\documentclass[11pt,twoside]{amsart}

\usepackage{pstricks}
\usepackage{multido}

\usepackage{graphicx}

\usepackage{amssymb}
\usepackage{amsmath}
\usepackage{amstext}
\usepackage{amsbsy}
\usepackage{amsxtra}
\usepackage{amsfonts}
\usepackage{latexsym}
\usepackage{alltt}
\usepackage{bbm}

\usepackage{verbatim}
\usepackage{amscd}
\usepackage[all]{xy}

\theoremstyle{plain}
\newtheorem{theorem}{Theorem}[section]
\newtheorem{corollary}[theorem]{Corollary}
\newtheorem{lemma}[theorem]{Lemma}
\newtheorem{proposition}[theorem]{Proposition}
\newtheorem{remark}[theorem]{Remark}
\newtheorem{example}[theorem]{Example}
\newtheorem{definition}[theorem]{Definition}

\newcommand{\LSL}{\Lambda  \SL}
\newcommand{\LSU}{\Lambda  \SU}

\newcommand{\Lsu}{\Lambda \su}

\newcommand{\moduli}{\mathcal{M}}
\newcommand{\C}{\mathbb C}
\newcommand{\R}{\mathbb R}

\newcommand{\Sp}{\mathbb S}
\newcommand{\Z}{\mathbb Z}
\newcommand{\N}{\mathbb N}

\newcommand{\CP}{\mathbb{CP}^1}

\newcommand{\tr}{\mathrm{tr}}

\newcommand{\Order}{{\rm O}}

\newcommand{\lbcl}{c}
\newcommand{\ubdh}{D}
\newcommand{\ubcl}{c\ind{max}}
\newcommand{\II}{\mathfrak{h}}
\newcommand{\norm}{n}
\newcommand{\pk}{\Lambda_{-1}^g\Sl(\C)}
\newcommand{\pksg}{\Lambda_{-1}^g\Sl(\C)^{\times}}
\newcommand{\pkcmc}{\mathcal{C}_g}
\newcommand{\vn}{}

\DeclareMathOperator{\dist}{dist}

\newcommand{\ind}[1]{_{\mbox{\rm\scriptsize{#1}}}}
\newcommand{\Con}[1]{{\sf({#1})}}

\newcommand{\SL}{\mathrm{SL}_{\mbox{\tiny{$2$}}}}

\newcommand{\Sl}{\mathfrak{sl}_{\mbox{\tiny{$2$}}}}

\newcommand{\SU}{\mathrm{SU}_{\mbox{\tiny{$2$}}}}

\newcommand{\su}{\mathfrak{su}_{\mbox{\tiny{2}}}}

\numberwithin{equation}{section}

\title[{\sc{cmc}} cylinders in $\Sp^3$]
{On the moduli of constant mean curvature cylinders of finite type
in the 3-sphere}

\author{M. Kilian}
\address{Department of Mathematics, University College Cork, Cork,
Ireland}
\email{m.kilian@ucc.ie}

\author{M. U. Schmidt}
\address{Institut f\"ur Mathematik,
Universit\"at Mannheim, 68131 Mannheim, Germany}
\email{schmidt@math.uni-mannheim.de}

\begin{document}



\begin{abstract}
  We show that one-sided Alexandrov embedded constant mean
  curvature cylinders of finite type in the 3-sphere
  are surfaces of revolution. This confirms a conjecture by Pinkall
  and Sterling that the only embedded constant mean
  curvature tori in the 3-sphere are rotational.
\end{abstract}
\maketitle


\begin{center}
  \sc Introduction
\end{center}

Alexandrov \cite{Ale0} proved that there are no compact embedded
surfaces with constant mean curvature ({\sc{cmc}}) in Euclidean
3-space $\R^3$ other than round spheres. However, while there are no
compact minimal surfaces in $\R^3$, there is an abundance of such in
the 3-sphere $\Sp^3$. For instance 2-spheres in the 3-sphere are
minimal precisely when they are great 2-spheres, and Lawson proved
that compact embedded minimal surfaces in $\Sp^3$ exist for every
genus \cite{Law:compact,Law:S3}. Lawson further showed
\cite{Law:unknot} that any embedded minimal torus in $\Sp^3$ is
unknotted, and conjectured that up to isometry the Clifford torus is
the only embedded minimal torus in $\Sp^3$. Hsiang and Lawson
\cite{HsiL} proved that the only embedded minimal torus of
revolution is the Clifford torus. Further results suggest that an
embedded minimal torus indeed has additional symmetries: Montiel and
Ros \cite{MonR:min} showed that the only minimal torus immersed into
$\Sp^3$ by the first eigenfunctions is the Clifford torus, and Ros
\cite{Ros} proved that the normal surface of an embedded minimal
torus in $\Sp^3$ is also embedded. Various methods for obtaining
minimal surfaces in $\Sp^3$ have been employed to study specific
classes, as in Karcher, Pinkall and Sterling \cite{KarPS}, and more
recently by Kapouleas and Yang \cite{KapY}.

Wente's discovery \cite{Wen} of {\sc{cmc}} tori provided the first
compact examples other than spheres in Euclidean 3-space. The
studies of Abresch \cite{Abr, Abr:twi}, Wente \cite{Wen:twi} and
Walter \cite{Wal} on special classes of {\sc{cmc}} tori in $\R^3$
concluded in the classification by Pinkall and Sterling \cite{PinS},
and their algebro-geometric description by Bobenko
\cite{Bob:tor,Bob:cmc}. In fact, Bobenko gave explicit formulas for
{\sc{cmc}} tori in $\R^3,\,\Sp^3$ and hyperbolic 3-space
$\mathbb{H}^3$ in terms of theta--functions, and provided a unified
description of {\sc{cmc}} tori in the 3-spaceforms in terms of
algebraic curves and spectral data. Independently, Hitchin
\cite{Hit:tor} classified harmonic 2-tori in the 3-sphere, and thus
also as a special case the harmonic Gauss maps of {\sc{cmc}} tori.
These ideas culminated in the description of harmonic tori in
symmetric spaces by Burstall, Ferus, Pedit and Pinkall
\cite{BurFPP}, and the generalized Weierstra{\ss} representation by
Dorfmeister, Pedit and Wu \cite{DorPW}.

Associated to a {\sc{cmc}} torus in the 3-sphere is a hyperelliptic
Riemann surface, the so called spectral curve. The structure
equation for {\sc{cmc}} tori is the $\sinh$-Gordon equation. Hitchin
\cite{Hit:tor}, and Pinkall and Sterling \cite{PinS} independently
proved that all doubly periodic solutions of the $\sinh$-Gordon
equation correspond to spectral curves of finite genus. The genus of
the spectral curve is called spectral genus. Ercolani, Kn\"orrer and
Trubowitz \cite{ErcKT} proved that for every even spectral genus
$g\geq 2$ there exists a hyperelliptic curve which corresponds to an
immersed {\sc{cmc}} torus in $\R^3$. The remaining cases of odd
genera $g>1$ was settled by Jaggy \cite{Jag}. Adapting these
results, Carberry \cite{Car:thesis} showed that minimal tori in
$\Sp^3$ exist for every spectral genus. Note that while a {\sc{cmc}}
torus in $\R^3$ has at least spectral genus 2, there is no such
restriction for {\sc{cmc}} tori in $\Sp^3$. In particular,
{\sc{cmc}} tori of revolution in $\Sp^3$ have spectral genus $g \leq
1$. Pinkall and Sterling \cite{PinS} conjectured that the only
embedded {\sc{cmc}} tori in $\Sp^3$ are tori of revolution.

\begin{figure}
  \centering
  \includegraphics[scale=0.72]{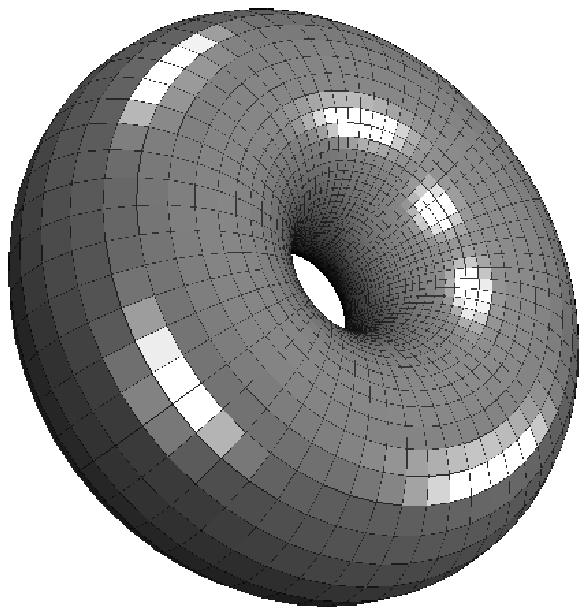}
  \includegraphics[scale=0.28]{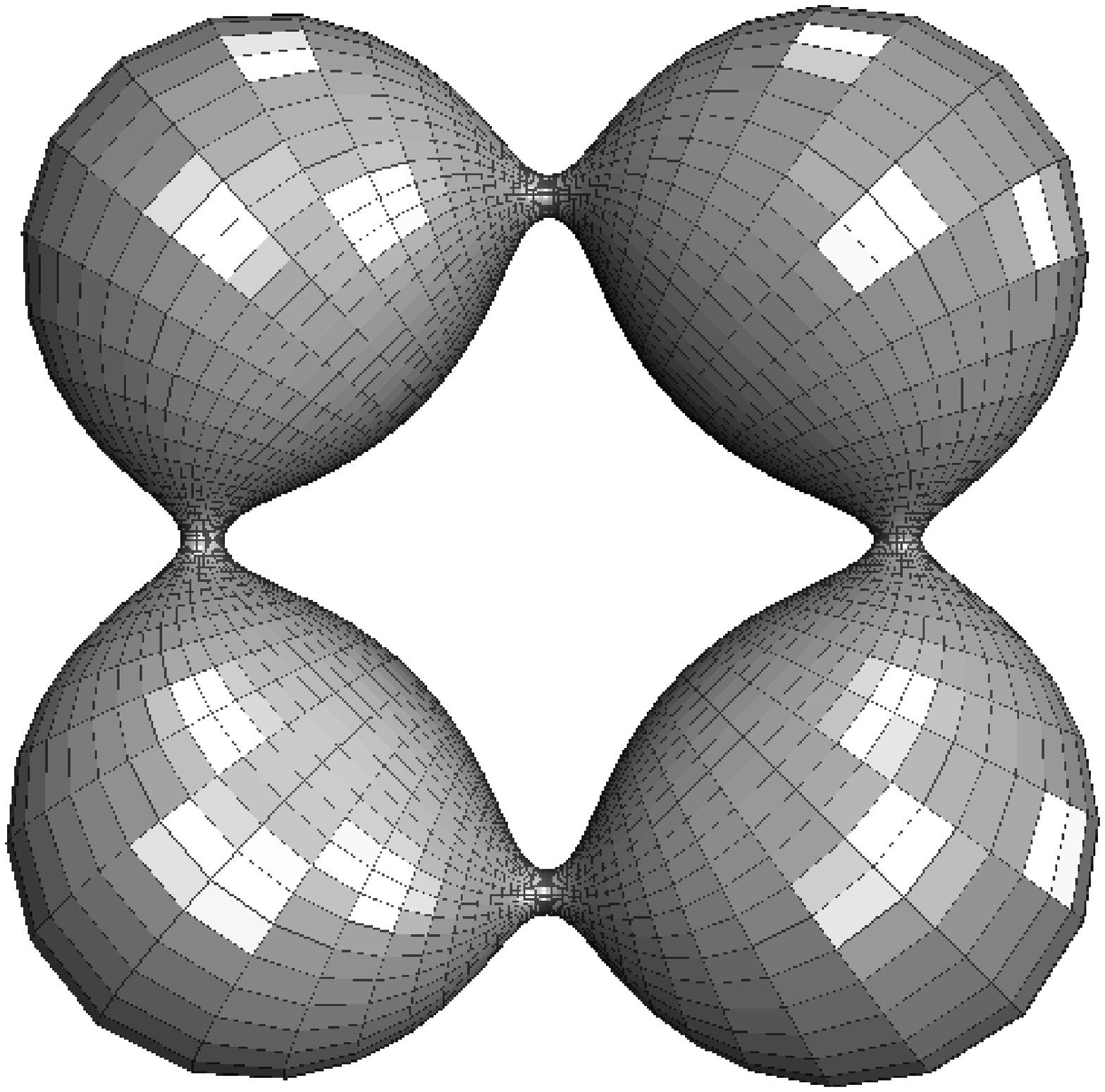}
  \includegraphics[scale=0.72]{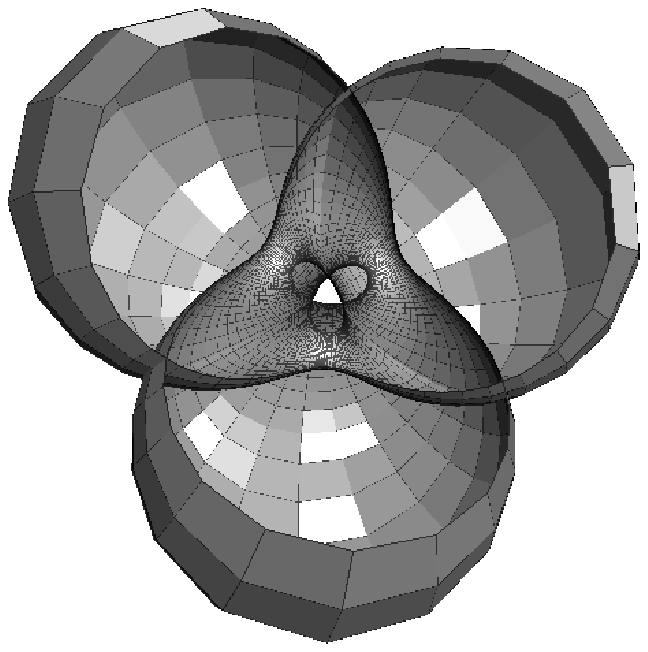}
  \caption{Stereographic projections of
    constant mean curvature tori of revolution in
  $\mathbb{S}^3$:
  On the left the Clifford torus,
  an embedded minimal torus with spectral genus 0.
  In the middle, an embedded non-minimal torus of spectral genus 1.
  On the right, a cutaway view of a
  non-embedded minimal torus of spectral genus 1.
  Images were created with cmclab \cite{Sch:cmclab}.}
  \label{fig:tori}
\end{figure}

In parallel the global theory of embedded {\sc{cmc}} (especially
minimal) surfaces in space forms was developed using geometric
{\sc{PDE}} methods. Meeks \cite{Mee} proved that a properly embedded
end of a {\sc{cmc}} surface is cylindrically bounded, which was used
by Korevaar, Kusner and Solomon \cite{KorKS} to prove that the only
embedded {\sc{cmc}} cylinders in $\R^3$ are surfaces of revolution -
either a standard round cylinder (spectral genus $g=0$), or a
Delaunay unduloid (spectral genus $g=1$). There are analogous
{\sc{cmc}} surfaces of revolution in $\Sp^3$, some of which close up
into tori, see Figure \ref{fig:tori} for some simple examples.

Kapouleas \cite{Kap1, Kap2, Kap3} proved the existence of compact
{\sc{cmc}} surfaces in $\R^3$ for any genus greater than 1, as well
as many new classes of non-compact {\sc{cmc}} surfaces in $\R^3$,
but not much is known about the moduli of {\sc{cmc}} surfaces in
general. Progress on understanding the moduli of {\sc{cmc}}
immersions of punctured spheres has recently been made in the case
of three punctures by Grosse-Brauckmann, Kusner and Sullivan
\cite{GroKS:Tri} and by Schmitt et. al. \cite{SKKR}. Kusner, Mazzeo
and Pollack \cite{KusMP} show that the moduli space of {\sc{cmc}}
surfaces is an analytic variety. The local linearization of the
moduli space is described by Jacobi fields which correspond to a
normal variation of the surface which preserve the constant mean
curvature property. Recently Korevaar, Kusner and Ratzkin
\cite{KorKR} studied Jacobi fields on a class of {\sc{cmc}} surfaces
with the additional property of being Alexandrov immersed. An
Alexandrov immersed surface in $\R^3$ is a complete noncompact
properly immersed surface $f:\Sigma \to \R^3$ that is the boundary
of a 3-manifold $M$ with two additional features: The mean curvature
normal of $\Sigma$ points into $M$, and $f$ extends to a proper
immersion of $M$ into $\R^3$. When the target is the 3-sphere, we
replace properness by completeness, and as Lawson \cite{Law:unknot}
we consider in analogy a smooth immersion $f:\Sigma \to \Sp^3$ that
we call a one-sided Alexandrov embedding if $\Sigma$ is the boundary
of a connected 3-manifold $M$ and the following two conditions hold:
The mean curvature of $\Sigma$ with respect to the inward normal is
non-negative. Secondly, the manifold $M$ is complete with respect to
the metric induced by $f$. We prove that the property of one-sided
Alexandrov embeddedness is stable under continuous deformation,
which allows us to study continuous families of one-sided Alexandrov
embedded surfaces.

In this paper we consider {\sc{cmc}} cylinders which have constant
Hopf differential, and whose metric is a periodic solution of the
$\sinh$-Gordon equation of finite type. Such {\sc{cmc}} cylinders
are said to be of finite type. We describe such finite type
{\sc{cmc}} cylinders in Section~\ref{sec:deformation} and
\ref{sec:moduli} by spectral data, and show that the spectral data
can be deformed in such a way that the corresponding family of
{\sc{cmc}} surfaces are all topologically cylinders. It turns out
that the corresponding moduli space of spectral data of genus $g$ is
$g+1$-dimensional. Furthermore, we can control the spectral genus
under the deformation, and by successively coalescing branch points
of the spectral curve, we continuously deform the spectral curve in
Lemma~\ref{thm:connected} into a curve of genus zero. In
Section~\ref{sec:AE space} we show that one-sided Alexandrov
embedded surfaces with constant mean curvature have collars with
depths uniformly bounded from below. For this purpose we use a
'maximum principle at infinity' which was communicated to us by
Harold Rosenberg \cite{Ros:com}. This allows us to show in
Theorem~\ref{thm:spec_AE} that a large class of continuous
deformations of {\sc{cmc}} cylinders of finite type preserve the
one-sided Alexandrov embeddedness. In Lemma~\ref{thm:continuous
genus 1} we continuously deform any one-sided finite type {\sc{cmc}}
cylinder in $\Sp^3$ into a one-sided Alexandrov embedded flat
cylinder in $\Sp^3$ with spectral genus zero. These are classified
in Theorem~\ref{thm:onesided revolution}. Finally this
classification is extended to all possible deformations of these
flat {\sc{cmc}} cylinders in $\Sp^3$ in Theorem~\ref{thm:main1}.
Since an embedded {\sc{cmc}} torus in the 3-sphere is covered by a
one-sided Alexandrov embedded cylinder, our result confirms the
conjecture by Pinkall and Sterling, and since the only embedded
minimal torus of revolution is the Clifford torus, also affirms
Lawson's conjecture.

\textbf{Acknowledgments. } We thank Fran Burstall, Karsten
Grosse-Brauckmann, Ian McIntosh, Rob Kusner, Franz Pedit and Ulrich
Pinkall for useful discussions. This work was mostly carried out
whilst M Kilian was a research assistant at the University of
Mannheim, and he would like to thank the Institute of Mathematics
there for providing excellent research conditions.

We had several beneficial conversations with Laurent Hauswirth, and
we are especially grateful to Antonio Ros and Harold Rosenberg who
helped us close a gap in a first draft of this paper, and
Univerit\'{e} Paris 7 for its hospitality during these discussions.
%
%
\section{Conformal cmc immersions into $\Sp^3$}
\label{sec:CONF}
This preliminary section recalls the relationship between {\sc{cmc}}
immersed surfaces in $\Sp^3$ and solutions of the $\sinh$-Gordon
equation, before considering the special case of {\sc{cmc}}
cylinders in $\Sp^3$, the notion of monodromy and the period
problem.
%
%
\subsection{The $\sinh$-Gordon equation}
We identify the 3-sphere $\Sp^3 \subset \R^4$ with $\Sp^3 \cong
\SU$. The Lie algebra of the matrix Lie group $\SU$ is $\su$,
equipped with the commutator $[\,\cdot,\,\cdot\,]$. For
$\alpha,\,\beta \in \Omega^1(\R^2,\su )$ smooth $1$--forms on $\R^2$
with values in $\su$, we define the $\su$--valued $2$--form
\begin{equation*}
  [ \alpha \wedge \beta ] (X,\,Y) =
  [\alpha(X),\,\beta(Y)] - [\alpha(Y),\,\beta(X)],
\end{equation*}
for vector fields $X,\,Y$ on $\R^2$. Let $L_g:h \mapsto gh$ be left
multiplication in $\SU$. Then by left translation, the tangent
bundle is $T\SU \cong \SU \times \su$ and $\theta : T\SU \to \su,\,
v_g \mapsto dL_{g^{-1}}(v_g)$ is the (left) Maurer--Cartan form. It
satisfies the {\bf{Maurer-Cartan-equations}}
\begin{equation} \label{eq:MC_equation}
 2\, d\theta + [ \theta \wedge \theta ] = 0.
\end{equation}
For a map $F:\R^2 \to \SU$, the pullback $\alpha = F^{\ast}\theta$
also satisfies \eqref{eq:MC_equation}, and conversely, every
solution $\alpha \in \Omega^1(\R^2,\su)$ of \eqref{eq:MC_equation}
integrates to a smooth map $F:\R^2 \to \SU$ with $\alpha =
F^{\ast}\theta$.

Complexifying the tangent bundle $T\R^2\cong T\C$ and decomposing
into $(1,0)$ and $(0,1)$ tangent spaces, and writing $d=\partial +
\bar{\partial}$, we may split $\omega \in \Omega^1(M,\su)$ into the
$(1,0)$ part $\omega^\prime$, the $(0,1)$ part $\omega^{\prime
\prime}$ and write $\omega = \omega^\prime + \omega^{\prime\prime}$.
We set the $*$--operator on $\Omega^1(M,\Sl)$ to $*\omega = - i
\omega^\prime + i \omega^{\prime \prime}$.

We denote by $\langle \cdot \, , \cdot \rangle$ the bilinear
extension of the Ad--invariant inner product $(X,\,Y) \mapsto
-\tfrac{1}{2} \tr(XY)$ of $\su$ to $\su^{\mbox{\tiny{$\C$}}} =
\Sl(\C)$. The double cover of the isometry group $\mathrm{SO}(4)$ is
$\SU \times \SU$ via the action $((F,G),X) \mapsto FXG^{-1}$.

Now let $g:\R^2 \to \SU$ be an immersion and $\omega = g^{-1}dg =
\omega' + \omega''$. Then $g$ is conformal if and only if the
$(1,0)$-part of $\omega$ is isotropic
\begin{equation} \label{eq:conformal}
  \langle \omega' ,\, \omega' \rangle = 0\,.
\end{equation}
If $g$ is a conformal immersion then there exists a smooth function
$v:\R^2 \to \R^\times = \R \setminus \{ 0 \}$, called the
{\bf{conformal factor}} of $g$ such that
\begin{equation} \label{eq:conformal_factor}
   v^2 = 2\,\langle \omega',\,\omega'' \rangle\,.
\end{equation}
The mean curvature function $H$ of $g$ (see e.g. \cite{SKKR}) is
given by
\begin{equation} \label{eq:H_S3}
  2\,d*\omega = H \,[ \omega \wedge \omega ].
\end{equation}
Recall the following observation of Uhlenbeck \cite{Uhl}, based on
an earlier result by Pohlmeyer \cite{Poh}, and suppose that in the
following $g:\R^2 \to \SU$ is a conformal immersion with non-zero
constant mean curvature $H_0$ and conformal factor $v$. Then
\eqref{eq:H_S3} and $2d\omega + [\,\omega \wedge \omega \,]=0$
combined give $d\omega+ H_0^{-1} d*\omega =0$, or alternatively
\begin{equation} \label{eq:w'}
  (1-iH_0^{-1})\,d\omega' + (1+iH_0^{-1})\,d\omega'' = 0.
\end{equation}
Inserting $d\omega'' = -d\omega' - [\,\omega' \wedge \omega'' \,]$
respectively $d\omega' = -d\omega'' - [\,\omega' \wedge \omega''
\,]$ into \eqref{eq:w'} gives $2d\omega' = (iH_0-1)[\,\omega' \wedge
\omega''\,]$ and $2d\omega'' = -(1+iH_0)[\,\omega' \wedge
\omega''\,]$. Then an easy computation shows that
\begin{equation*}
  \alpha_\lambda = \tfrac{1}{2}(1+\lambda^{-1})(1+iH_0)\,\omega'
  + \tfrac{1}{2}(1+\lambda)(1-iH_0)\,\omega''
\end{equation*}
satisfies the Maurer-Cartan-equations
$$
    2\,d\alpha_\lambda + [\,\alpha_\lambda \wedge \alpha_\lambda \,] =
    0 \mbox{ for all $\lambda \in \C^{\times} = \C \setminus \{ 0\}$. }
$$
The Maurer-Cartan-equations are an integrability condition, so we
can integrate and obtain a corresponding {\bf{extended frame}}
$F_\lambda : \R^2 \times \C^\times \to \SL(\C)$ with $dF_\lambda =
F_\lambda\,\alpha_\lambda$ and $F_\lambda(0) = \mathbbm{1}$. Since
$\omega$ takes values in $\su$, we conclude that $F_\lambda$ takes
values in $\SU$ when $\lambda \in \Sp^1$. Now define for
$\lambda_0,\,\lambda_1 \in \Sp^1,\, \lambda_0 \neq \lambda_1$, the
following map $f:\R^2 \to \SU$ by the {\bf{Sym-Bobenko-formula}}
\begin{equation}\label{eq:Sym_S3}
  f = F_{\lambda_1} F_{\lambda_0}^{-1}\,.
\end{equation}
Then for $\Omega = f^{-1}df = \mathrm{Ad}\,F_{\lambda_0} \left(
    \alpha_{\lambda_1} - \alpha_{\lambda_0} \right)$ we obtain
    $\Omega' = \tfrac{1}{2}(\lambda_1^{-1} - \lambda_0^{-1})(1+iH_0)\,
    \mathrm{Ad}\,F_{\lambda_0} \,\omega'$
so the conformality of $f$ follows from the conformality of $g$,
since
\begin{align}
    \langle \Omega',\,\Omega' \rangle &=
        \tfrac{1}{4}(\lambda_1^{-1} - \lambda_0^{-1})^2(1+iH_0)^2\,
    \mathrm{Ad}\,F_{\lambda_0} \,\langle \omega',\,\omega' \rangle =
    0\\
    2\,\langle \Omega',\,\Omega'' \rangle &=
    \sin^2(t_1-t_0)(1+H_0^2)\,v^2\,. \label{eq:conf_Om}
\end{align}
Here we have written $\lambda_{0,1} = e^{2it_{0,1}}$, and $v$ is the
conformal factor of $g$. Furthermore
\begin{align*}
    d\ast\Omega &= \tfrac{i}{4}(1+H_0^2)
    \left( \lambda_0^{-1}\lambda_1 - \lambda_0\lambda_1^{-1}
    \right) \mathrm{Ad}\,F_{\lambda_0} [\omega' \wedge \omega'']\,, \\
    [\Omega \wedge \Omega] &=
    \tfrac{1}{2}(1+H_0^2)(\lambda_1^{-1}-\lambda_0^{-1})(\lambda_1 -
    \lambda_0)\,\mathrm{Ad}\,F_{\lambda_0}\,[\omega' \wedge
    \omega'']\,.
\end{align*}
Hence by \eqref{eq:H_S3}, the map $f:\R^2 \to \SU$ given by
\eqref{eq:Sym_S3} has constant mean curvature
\begin{equation} \label{eq:H_mu}
  H = i \,\frac{\lambda_0 + \lambda_1}{\lambda_0 -\lambda_1}\,.
\end{equation}
In summary, by starting with one non-minimal conformal {\sc{cmc}}
immersion $g$, we have just seen how to obtain a whole
$\C^\times$-family of solutions of the Maurer-Cartan-equations, and
from the corresponding extended frame we then obtained another
conformal {\sc{cmc}} immersion $f$. Since the mean curvature
\eqref{eq:H_mu} and the conformal factor of $f$ in
\eqref{eq:conf_Om} only depend on the angle between
$\lambda_0,\,\lambda_1$, we in fact get a whole $\Sp^1$-family of
isometric conformal {\sc{cmc}} immersions, called an {\bf{associated
family}}, which is obtained by simultaneously rotating
$\lambda_0,\,\lambda_1$ while keeping the angle between them fixed.

We next recall the following version of Theorem 14.1 in Bobenko
\cite{Bob:cmc}, which provides a correspondence between solutions of
the $\sinh$-Gordon equation and associated families of {\sc{cmc}}
surfaces in the 3-sphere.
\begin{theorem} \cite{Bob:cmc} \label{thm:sinh}
Let $u:\R^2 \to \R$ be a smooth function and define
\begin{equation} \label{eq:general_alpha}
  \alpha_\lambda = \frac{1}{2}\,\begin{pmatrix}
  u_z\,dz-u_{\bar{z}}\,d\bar{z} &
  i\,\lambda^{-1}e^u\,dz + i\,e^{-u}\,d\bar{z}\\
  i\,e^{-u}\,dz + i\,\lambda\,e^{u}\,d\bar{z}&
  -u_z\,dz+u_{\bar{z}}\,d\bar{z}
  \end{pmatrix}\,.
\end{equation}
Then $2\,d\alpha_\lambda + [\,\alpha_\lambda \wedge
\alpha_\lambda\,] = 0$ if and only if $u$ is a solution of the
$\sinh$-Gordon equation
\begin{equation}\label{eq:sinh-Gordon}
  \partial \bar{\partial}\,2u + \sinh (2u) = 0.
\end{equation}
For any solution $u$ of the $\sinh$-Gordon equation and
corresponding extended frame $F_\lambda$, and $\lambda_0,\,\lambda_1
\in \Sp^1,\, \lambda_0 \neq \lambda_1$, the map defined by the
Sym-Bobenko-formula \eqref{eq:Sym_S3} is a conformal immersion with
constant mean curvature $H$ \eqref{eq:H_mu}, conformal factor $v =
e^u/\sqrt{H^2 + 1}$, and constant Hopf differential $Q\,dz^2$ with
$Q = i\,(\lambda_1^{-1}-\lambda_0^{-1})/4$.
\end{theorem}
\begin{proof} Decomposing $\alpha_\lambda = \alpha'_\lambda\,dz +
\alpha''_\lambda\,d\bar{z}$ into $(1,\,0)$ and $(0,\,1)$ parts, we
compute
\begin{equation*} \begin{split}
  &\bar{\partial} \alpha'_\lambda = \frac{1}{2}\,
  \begin{pmatrix} u_{z\bar{z}} & i\lambda^{-1}u_{\bar{z}}e^u \\
    -iu_{\bar{z}}e^{-u} & -u_{z\bar{z}} \end{pmatrix}\,,
  \quad \partial \alpha''_\lambda = \frac{1}{2}\,
  \begin{pmatrix} -u_{z\bar{z}} & -iu_z e^{-u} \\
    i\lambda u_z e^u & u_{z\bar{z}} \end{pmatrix}\,, \\
  &\left[ \alpha'_\lambda,\,\alpha''_\lambda \right] =
  \frac{1}{4} \begin{pmatrix} -e^{2u} + e^{-2u} &
    2iu_{\bar{z}}\lambda^{-1}e^u + 2iu_ze^{-u} \\
    -2i\lambda u_ze^u - 2iu_{\bar{z}}e^{-u} &
    e^{2u} - e^{-2u} \end{pmatrix}\,. \end{split}
\end{equation*}
Now $2\,d\alpha_\lambda + [\,\alpha_\lambda \wedge \alpha_\lambda\,]
= 0$ is equivalent to $\bar{\partial}\alpha'_\lambda - \partial
\alpha''_\lambda = [ \alpha'_\lambda,\,\alpha''_\lambda ]$, which
holds if and only if $u$ solves the sinh-Gordon equation
\eqref{eq:sinh-Gordon}.

If $u$ is a solution of the sinh-Gordon equation, then we may
integrate $dF_\lambda = F_\lambda\,\alpha_\lambda$ to obtain a map
$F_\lambda:\R^2 \times \mathbb{S}^1 \to \SU$. Let
$\lambda_0,\,\lambda_1 \in \mathbb{S}^1,\, \lambda_0 \neq
\lambda_1$, and consider the map $f:\R^2 \to \SU$ defined by the
Sym-Bobenko-formula \eqref{eq:Sym_S3}. Conformality
\eqref{eq:conformal} is a consequence of the fact that the
complexified tangent vector
\begin{equation*}
  f^{-1}\partial f = \mathrm{Ad}\,F_{\lambda_0} (\alpha'_{\lambda_1} -
  \alpha'_{\lambda_0}) =
  \frac{i}{2}\,e^u \,(\lambda_1^{-1}-\lambda_0^{-1})\,
  \mathrm{Ad}\,F_{\lambda_0}\,\bigl( \begin{smallmatrix}0&1\\0&0\end{smallmatrix}\bigr)
\end{equation*}
is isotropic with respect to the bilinear extension of the Killing
form. The mean curvature can be computed using formula
\eqref{eq:H_S3}. The conformal factor is obtained from
\begin{equation*}
  v^2 = 2\,\langle f^{-1}\partial f,\,f^{-1}\bar{\partial} f \rangle =
  \tfrac{1}{4}\,e^{2u}\,(\lambda_1^{-1}-\lambda_0^{-1})
  (\lambda_1-\lambda_0).
\end{equation*}
From \eqref{eq:H_mu} we have $(H^2+1)(\lambda_1^{-1}-\lambda_0^{-1})
(\lambda_1-\lambda_0) = 4$, which proves the formula for the
conformal factor.

Define the normal $N = F_{\lambda_1}
\,\varepsilon\,F_{\lambda_0}^{-1}$ with $\varepsilon = \bigl(
\begin{smallmatrix} i&0\\0&-i \end{smallmatrix} \bigr)$.
Then $\partial N = F_{\lambda_1} (\alpha'_{\lambda_1}\, \varepsilon
- \varepsilon\,\alpha'_{\lambda_0})\, F_{\lambda_0}^{-1}$ and
\begin{equation*}
  \alpha'_{\lambda_1}\,\varepsilon -
  \varepsilon\,\alpha'_{\lambda_0} = \begin{pmatrix}
    0&\tfrac{1}{2}e^u(\lambda_1^{-1}+\lambda_0^{-1})\\-e^{-u}&0
    \end{pmatrix}\,.
\end{equation*}
Consequently, $Q := -\langle \partial \partial f,\, N \rangle =
  \langle \partial f,\, \partial N \rangle =
  \langle F_{\lambda_1}^{-1}\,\partial f \,F_{\lambda_0},\,
  F_{\lambda_1}^{-1}\,\partial N \,F_{\lambda_0} \rangle
  = \tfrac{i}{4}\,(\lambda_1^{-1} - \lambda_0^{-1})$,
which proves the formula for the Hopf differential, and concludes
the proof.
\end{proof}
There is an analogous but more general theorem (see e.g Bobenko
\cite{Bob:2x2}) than the one above, which asserts that if functions
$(u,\,Q,\,H \equiv const.)$ satisfy the Gauss-Codazzi equations,
then one obtains a $\C^\times$-family of solutions of the
Maurer-Cartan-equations, thus an extended frame and consequently an
associated family via the Sym-Bobenko-formula.
%
%
\subsection{Monodromy and periodicity condition}
The {\sc{cmc}} condition implies that the Hopf differential is a
holomorphic quadratic differential \cite{Hop}. On the cylinder
$\C^{\times}$ there is an infinite dimensional space of holomorphic
quadratic differentials, large classes of which can be realized as
Hopf differentials of {\sc{cmc}} cylinders \cite{KilMS}. On a
{\sc{cmc}} torus the Hopf differential is constant (and non-zero).
Since we are ultimately interested in tori, we restrict our
attention to {\sc{cmc}} cylinders considered via
Theorem~\ref{thm:sinh} which have constant non-zero Hopf
differentials on the universal covering $\C$ of $\C^{\times}$. Note
that for given solution $u$ of the sinh-Gordon equation an extended
frame $\lambda \mapsto F_\lambda$ is holomorphic on $\C^{\times}$
and has essential singularities at $\lambda = 0,\,\infty$.

Let $F_\lambda$ be an extended frame for a {\sc{cmc}} immersion
$f:\R^2 \to \Sp^3$ such that \eqref{eq:Sym_S3} holds for two
distinct unimodular numbers $\lambda_0,\,\lambda_1$. Let $\tau:\R^2
\to \R^2,\,z\mapsto z+\tau$ be a translation, and assume that
$\alpha_\lambda = F_\lambda^{-1}dF_\lambda$ has period $\tau$, so
that $\tau^* \alpha_\lambda = \alpha_\lambda \circ \tau =
\alpha_\lambda$. Then we define the {\bf{monodromy}} of $F_\lambda$
with respect to $\tau$ as
\begin{equation}\label{eq:monodromy}
    M_\lambda(\tau) = \tau^{\ast}(F_\lambda)\,F_\lambda ^{-1}\,.
\end{equation}
Periodicity $\tau^{\ast} f = f$ in terms of the monodromy is then
$\tau^{\ast} f = M_{\lambda_1}(\tau) F_{\lambda_1}
F_{\lambda_0}^{-1} M_{\lambda_0}^{-1}(\tau)$, so $\tau^{\ast} f = f$
if and only if
\begin{equation}\label{eq:periodicity}
    M_{\lambda_0}(\tau) = M_{\lambda_1}(\tau) = \pm \mathbbm{1}\,.
\end{equation}
If $\Delta(\lambda)$ is the trace of $M_{\lambda}(\tau)$ then
$\tau^{\ast} f = f$ if and only if
$\Delta(\lambda_0)=\Delta(\lambda_1) = \pm 2$.
%
%
\section{Finite type solutions of the $\sinh$-Gordon equation}
\label{sec:finite type}
In this section we introduce the solutions of the $\sinh$-Gordon
equation which are called finite type solutions. Finite type
solutions of the $\sinh$-Gordon equation are in one-to-one
correspondence with maps called polynomial Killing fields. These
polynomial Killing fields take values in certain $2\times 2$-matrix
polynomials, and solve a non-linear partial differential equation,
but they are uniquely determined by one of their values. We shall
call these values initial values of polynomial Killing fields or
just initial values. The Symes method calculates the solutions in
terms of the initial values with the help of a loop group splitting.
The eigenvalues of these matrix polynomials define a real
hyperelliptic algebraic curve, which is called {\bf{spectral
curve}}. One spectral curve corresponds to a whole family of finite
type solutions of the $\sinh$-Gordon equation. We call the sets of
finite type solutions (or their initial values), which belong to the
same spectral curve, isospectral sets. The eigenspaces of the matrix
polynomials define a holomorphic line bundle on the spectral curves
called eigenbundle. These holomorphic line bundles completely
determine the corresponding initial value and the corresponding
solution of the $\sinh$-Gordon equation. Consequently, the
isospectral sets can be identified with one connected component of
the real part of the Picard group. In case the spectral curve has
singularities, then the isospectral set can be identified with the
real part of the compactification of a generalized Jacobian. These
compactifications have stratifications, whose strata are the orbits
under the action of the generalized Jacobian. In our case the
spectral curves are hyperelliptic and we shall describe the
corresponding stratifications of the isospectral sets.
%
%
\subsection{Polynomial Killing fields} \label{sec:killing}
For some aspects of the theory untwisted loops are advantageous, and
avoiding the additional covering map $\lambda \mapsto \sqrt\lambda$
simplifies for example the description of Bianchi-B\"acklund
transformations by the simple factors \cite{TerU, KilSS}. For the
description of polynomial Killing fields on the other hand, the
twisted loop algebras as in
\cite{BurFPP,BurP_adl,BurP:dre,DorPW,McI:tor} are better suited, but
we remain consistent and continue working in our 'untwisted'
setting.

Let $\varepsilon_+ = \bigl( \begin{smallmatrix} 0 & 1 \\ 0 & 0
  \end{smallmatrix} \bigr)$
and $\varepsilon_- = \varepsilon_+^t$, and consider for $g \in \N_0$
the finite dimensional vector space
$$
    \pk = \left\{
    \xi=\sum_{d=-1}^{g}\xi_d\lambda^d
    \mid\xi_{-1}\in\mathbb{C}\varepsilon_+,\,
    \xi_d=-\bar{\xi}^t_{g-1-d}\in\Sl(\C)\mbox{ for }d=-1,\ldots,g
    \right\}.
$$
Clearly $\pk$ is a real $3g+2$-dimensional vector space and has up
to isomorphism a unique norm $\|\cdot\|$. These Laurent polynomials
define smooth mappings from $\lambda\in\Sp^1$ into $\Sl(\C)$. Note
that $\sqrt{\lambda}\mapsto\lambda^{\frac{1-g}{2}}\xi$
belongs to the loop Lie algebra $\Lsu$ of the loop Lie group
$\LSU$. For the resulting
solution of the $\sinh$-Gordon equation to be of finite type, we
need in addition the conditions $\tr(\xi_{-1}\varepsilon_-)\neq 0
\neq \tr(\xi_0\varepsilon_+)$. These conditions ensure that
$\xi_{-1}$ and the lower left entry of $\xi_0$ do not vanish, and
therefore that $\xi_{-1}+\tr(\xi_0\varepsilon_+)\,\epsilon_-$ is
semisimple. This is the same as semisimplicity of the leading order
term in the twisted setting. We thus define
$$
    \pksg=\left\{\xi\in\pk\mid
    \tr(\xi_{-1}\varepsilon_-)\neq 0\neq\tr(\xi_0\varepsilon_+)
    \right\}.
$$
By the Symes method \cite{Symes_80}, elucidated by Burstall and
Pedit \cite{BurP_adl,BurP:dre}, the extended framing $F_\lambda:\R^2
\to \Lambda \SU$ of a {\sc{cmc}} immersion of finite type is given
by the unitary factor of the Iwasawa decomposition of
\begin{equation} \label{eq:FB}
  \exp(z\,\xi) = F_\lambda\,B
\end{equation}
for some $\xi\in\pksg$ with $g\in\mathbb{N}_0$. Due to Pressly and
Segal \cite{PreS}, the Iwasawa decomposition is a diffeomorphism
between the loop group $\LSL(\C)$ of $\SL(\C)$ into point wise
products of elements of $\LSU$ with elements of the loop group
$\Lambda^+\SL(\C)$ of holomorphic maps from $\lambda\in\mathbb{D}$
to $\SL(\C)$, which take at $\lambda=0$ values in the subgroup of
$\SL(\C)$ of upper-triangular matrices with positive real diagonal
entries. For every $\xi\in\pk$ there exists a unique
$\alpha(\xi)\in\Omega^1(\R^2,\,\Lambda_{-1}^1\Sl(\C))$, such that
$\xi dz-\alpha(\xi)$ takes values in the Lie algebra of
$\Lambda^+\SL(\C)$ of the right hand factor in the Iwasawa
decomposition~\eqref{eq:FB}.

A {\bf{polynomial Killing field}} is a map $\zeta:\R^2\to\pk$ which
solves
\begin{equation} \label{eq:pKf}
    d\zeta =[\,\zeta,\,\alpha(\zeta)\,]\quad\mbox{ with }\quad
    \zeta(0)=\xi.
\end{equation}
For each {\bf initial value} $\xi\in\pk$, there
exists a unique polynomial Killing field given by
\begin{equation}\label{eq:solution pk}
\zeta=B\xi B^{-1}=F^{-1}_\lambda\xi F_\lambda\quad
\mbox{ with }F_\lambda\mbox{ and }B\mbox{ as in \eqref{eq:FB}.}
\end{equation}
For $\xi\in\pksg$ with $\tr(\xi_{-1}\varepsilon_-)\in \mathbb{R}^+i$
and $\tr(\xi_0\varepsilon_+)\in\mathbb{R}^+i$ the corresponding
$\Lambda_{-1}^1\Sl(\C)$-valued 1-form $\alpha(\zeta)$ is the
$\alpha$ as in \eqref{eq:general_alpha} for that particular solution
$u$ of the $\sinh$-Gordon equation corresponding to the extended
frame $F_\lambda$ of \eqref{eq:FB}. For general initial values
$\xi\in\pksg$ the leading term
$\zeta_{-1}+\tr(\zeta_0\varepsilon_+)\,\epsilon_-$ of the
corresponding polynomial Killing field does not depend on the
surface parameter $z$. The corresponding $\alpha(\zeta)$ differs
from \eqref{eq:general_alpha} by multiplication of $\lambda$ and
$dz$ with constant unimodular complex numbers. Given a polynomial
Killing field $\zeta$, we set the initial value $\xi=\zeta|_{z=0}$
in \eqref{eq:FB}. Thus $\zeta$, or the initial value $\xi$, gives
rise to an extended frame, and thus to an associated family.
\begin{definition}\label{def:finite type}
A solution of the $\sinh$-Gordon equation is called a {\bf finite
type} solution if and only if it corresponds to a polynomial Killing
field $\zeta:\R^2 \to \pksg$ with $g\in\N_0$.
\end{definition}
%
%
\subsection{Roots of polynomial Killing fields}
If an initial value $\xi$ has a root at some
$\lambda=\alpha\in\C^{\times}$, then the corresponding polynomial
Killing field has a root at the same $\lambda$ for all $z\in\C$.
In this case we may reduce the order of $\xi$ and $\zeta$ without
changing the corresponding extended frame $F$ \eqref{eq:FB}. The
following polynomials transform under
$\lambda\mapsto\bar{\lambda}^{-1}$ as
\begin{align}\label{eq:transformation}
  p(\lambda)&=\begin{cases}
  i(\sqrt{\bar{\alpha}}\lambda-\sqrt{\alpha})
  &\mbox{ for }\alpha\bar{\alpha}=1\\
  (\lambda-\alpha)(1-\bar{\alpha}\lambda)&
  \mbox{ for all }\alpha\in\C\end{cases}&
  \overline{\lambda^{\deg(p)}p\left(\bar{\lambda}^{-1}\right)}&=p(\lambda).
\end{align}
If the polynomial Killing field $\zeta$ with initial value
$\xi\in\pksg$ has a simple root at $\lambda=\alpha\in\C^{\times}$,
then $\zeta/p$ does not vanish at $\alpha$ and is the polynomial
Killing field with initial value
$\xi/p\in\Lambda_{-1}^{g-\deg(p)}\Sl(\C)^\times$. Furthermore,
obviously $\zeta$ and $\zeta/p$ commute, and we next show that both
polynomial Killing fields $\zeta$ and $\zeta/p$ give rise to the
same extended frame $F_\lambda$ \eqref{eq:FB}.
\begin{proposition} \label{th:pKf_min}
If a polynomial Killing field $\zeta$ with initial value $\xi\in\pk$
has zeroes in $\lambda\in\C^{\times}$, then there is a polynomial
$p(\lambda)$, such that the following two conditions hold:
\begin{enumerate}
    \item $\zeta/p$ is the polynomial Killing field with initial value
    $\xi/p\in\Lambda_{-1}^{g-\deg p}\Sl(\C)$, which gives rise to the same
    associated family as $\zeta$.
    \item $\zeta/p$ has no zeroes in $\lambda \in \C^{\times}$.
\end{enumerate}
\end{proposition}
\begin{proof}
An appropriate M\"obius transformation \eqref{eq:moebius} transforms
any root $\alpha\in\C^{\times}$ into a negative root. For such
negative roots the corresponding initial values $\xi$ and $\xi/p$
are related by multiplication with a polynomial with respect to
$\lambda$ with positive coefficients. In the Iwasawa decomposition
\eqref{eq:FB} this factor is absorbed in $B$. Hence the
corresponding extended frames coincide, which proves (i). Repeating
this procedure for every root $\lambda \in \C^\times$ ensures (ii).
\end{proof}
Hence amongst all polynomial Killing fields that give rise to a
particular {\sc{cmc}} surface of finite type there is one of
smallest possible degree (without adding further poles), and we say
that such a polynomial Killing field has \emph{minimal degree}. A
polynomial Killing field has minimal degree if and only if it has
neither roots nor poles in $\lambda \in \C^{\times}$. We summarize two
results by Burstall and Pedit \cite{BurP_adl,BurP:dre}. The first
part is a variant of Theorem 4.3 in \cite{BurP_adl}, the second part
follows immediately from results in \cite{BurP:dre}.
\begin{theorem} \label{thm:PKF}
(i) A {\sc{cmc}} immersion $f:\R^2 \to \mathbb{S}^3$ is of finite
type if and only if there exists a polynomial Killing field $\zeta$
with initial value $\xi\in\pksg$ such that the map $F_\lambda$
obtained from \eqref{eq:FB} is an extended frame of $f$.

(ii) In particular there exists a unique polynomial Killing field of
minimal degree that gives rise to $f$. Thus we have a smooth 1-1
correspondence between the set of {\sc{cmc}} immersions of finite
type and the set of polynomial Killing fields without zeroes.
\end{theorem}
\begin{proof}
Point (i) is a reformulation of Theorem 4.3 in \cite{BurP_adl}. (ii)
We briefly outline how to prove the existence and uniqueness of a
minimal element.

If the initial value $\xi$ gives rise to $f$, then the corresponding
polynomial Killing field $\zeta$ can be modified according to
Proposition \ref{th:pKf_min} so that $\tilde{\zeta}$ is of minimal
degree, and still giving rise to $f$. Hence there exists a
polynomial of least degree giving rise to $f$.

For the uniqueness, assume we have two initial values
$\xi,\,\tilde{\xi}$ of least degree $g$ both giving rise to $f$.
Putting Proposition~3.3 and Corollary~3.8 in \cite{BurP:dre}
together gives: Two finite type initial values give rise to the same
associated family if and only if they commute and have equal
residues. Since the residues coincide and both $\xi,\,\tilde{\xi}$
are of minimal degree, we conclude that $\xi \equiv \tilde{\xi}$.
The unique minimal polynomial Killing field is thus $\zeta =
F^{-1}_\lambda\xi \,F_\lambda$.

Since the Iwasawa factorization is a diffeomorphism, and all other
operations involved in obtaining an extended frame from $\zeta$ are
smooth, the resulting {\sc{cmc}} surface depends smoothly on the
entries of $\zeta$.
\end{proof}
%
%
\subsection{Spectral curves I}
Due to \eqref{eq:solution pk} the characteristic equation
\begin{equation}\label{eq:characteristic0}
\det\left(\nu\,\mathbbm{1}-\zeta\right)=\nu^2+\det(\zeta)=0
\end{equation}
of a polynomial Killing field $\zeta$ with initial value $\xi\in\pk$
does not depend on $z\in\C$ and agrees with the characteristic
equation of the initial value $\xi$. If $\xi\in\pk$ then we may
write $-\det\xi = \lambda^{-1}a$ for a polynomial $a$ of degree at
most $2g$ which satisfies the reality condition
\begin{equation} \label{eq:a_reality}
  \lambda^{2g} \overline{a(\bar{\lambda}^{-1})} =
  -a(\lambda).
\end{equation}
Consequently the hyperelliptic curve has three involutions
\begin{align}\label{eq:involutions}
\sigma&:(\lambda,\nu)\mapsto(\lambda,-\nu)&
\rho&:(\lambda,\nu)\mapsto(\bar{\lambda}^{-1},\bar{\lambda}^{-g}\bar{\nu})&
\eta&:(\lambda,\nu)\mapsto(\bar{\lambda}^{-1},-\bar{\lambda}^{-g}\bar{\nu})&
\end{align}
If $a$ has $2g$ pairwise distinct roots, then $\nu^2 =
\lambda^{-1}a(\lambda)$ is a spectral curve of genus $g$ of a not
necessarily periodic solution of the $\sinh$-Gordon equation. The
genus $g$ is called the {\bf{spectral genus}}.
\begin{lemma}\label{thm:compact}
Let $a$ be a polynomial of degree $2g$ satisfying
\eqref{eq:a_reality}. Then the {\bf isospectral set}
$$
  \mathcal{K}_a = \left\{ \xi\in\pk\mid
  \det \xi(\lambda) = -\lambda^{-1}a(\lambda) \right\}
$$
is compact. Furthermore, if the $2g$ roots of $a$ are pair wise
distinct, then $\mathcal{K}_a \cong \left(\Sp^1\right)^g$.
\end{lemma}
\begin{proof}
For the compactness it suffices to show that all Laurent
coefficients of a $\xi\in\mathcal{K}_a$ are bounded, since
$\mathcal{K}_a$ is a closed subset of the $(3g+2)$-dimensional
vector space $\pk$. For $d=(1-g)/2$ the product $\lambda^d\xi$ is
skew hermitian on $|\lambda|=1$. The negative determinant of
traceless skew hermitian $2 \times 2$ matrices is the square of a
norm. Hence for all $\xi\in\mathcal{K}_a$ the Laurent polynomial of
$\lambda^d\xi(\lambda)$ with respect to $\sqrt{\lambda}$ is bounded
on $|\lambda |=1$. Thus the Laurent coefficients are bounded.

If $a$ has $2g$ pairwise distinct roots, then $\xi$ has no roots
since at all roots of $\xi$, the determinant $\det \xi$ has a root
of order two. If $\alpha$ is a root of $a$, then $\det\xi(\alpha)$
vanishes and $\xi(\alpha)$ is nilpotent. For a nonzero nilpotent
$2\times 2$-matrix $\xi(\alpha)$ there exists a $2\times 2$-matrix
$Q$ such that $\xi(\alpha) =[Q,\,\xi(\alpha)]$. Hence for every root
$\alpha$ of $a$, there exists a $2\times 2$-matrix $Q$, such that
$$
    \dot{\xi}(\lambda) = \frac{\xi^{-1}\det
    \xi+[Q,\,\xi]}{\lambda-\alpha}=
    \frac{-\xi+[Q,\,\xi]}{\lambda-\alpha}
$$
has no pole at $\lambda=\alpha$. The corresponding derivative of
$a=-\lambda\det(\xi)$ is equal to
$\dot{a}=\frac{2a}{\lambda-\alpha}$. Furthermore, $\lambda\dot{\xi}$
is polynomial with respect to $\lambda$ of degree $g$. Two
appropriate linear combinations with the analogous tangent element
at the root $\bar{\alpha}^{-1}$ of $\xi$ change the roots $\alpha$
and $\bar{\alpha}^{-1}$ and fixes all other roots of $a$ and
respects the reality condition of $\pk$. These two linear
combinations belong to the tangent space of $\pksg$. Hence the
derivatives of all the coefficients of $a$ as functions on $\pk$ are
non-zero at all $\xi\in\mathcal{K}_a$. By the implicit function
theorem this set is therefore a $g$-dimensional submanifold. The
corresponding eigenspaces of $\xi$ depend holomorphically on the
solutions $(\lambda,\nu)$ of \eqref{eq:characteristic0} and define a
holomorphic line bundle on the spectral curve. These
{\bf{eigenbundles}} have degree $g+1$, they are non-special in the
sense that they have no holomorphic sections vanishing at one of the
points at $\lambda=0$ or $\lambda=\infty$, and finally they obey
some reality condition. Vice versa, all holomorphic line bundles
obeying these three conditions correspond to one
$\xi\in\mathcal{K}_a$ (see McIntosh \cite[Section~1.4]{McI:tor}).
Hitchin has shown in \cite{Hit:tor}, that the third condition
implies the second condition. Therefore $\mathcal{K}_a$ can be
identified with the real part of one connected component of the
Picard group of the spectral curve, which is a $g$-dimensional
torus.
\end{proof}
If $a$ has multiple roots, then the real part of the Jacobian of the
corresponding hyperelliptic curve still acts on $\mathcal{K}_a$, but
not transitively. More precisely, in case of non-unimodular multiple
roots of $a$ the set $\mathcal{K}_a$ has a stratification, whose
strata are the orbits of the action of the real part of the
generalized Jacobian of the singular hyperelliptic curve defined by
$\nu^2=\lambda\,a(\lambda)$. The elements of different strata have
different orders of zeroes at the multiple roots of $a$. In
Proposition~\ref{th:pKf_min} we have seen that all $\xi\in\pk$ are
products of $\tilde{\xi}\in\Lambda_{-1}^{\tilde{g}}\Sl(\C)$ of lower
degree $\tilde{g}<g$ with polynomials of the form
\eqref{eq:transformation}.
\begin{definition} \label{def:bubbletons}
Every finite type solution of the $\sinh$-Gordon equation
corresponds to a unique polynomial Killing field $\zeta$ without
zeroes and initial value $\xi\in\pksg$. The curve defined by
$\nu^2=-\det \xi$ has a unique compactification to a projective
curve without singularities at $\lambda=0$ and $\lambda=\infty$. If
$\det(\xi)$ has multiple roots, then we say that the solution
contains {\bf{bubbletons}}. The arithmetic genus of this
hyperelliptic curve is equal to $g$.
\end{definition}
%
%
\subsection{Bubbletons}
We briefly motivate Definition \ref{def:bubbletons} above, and refer
the reader to \cite{Bur:iso,BurP:dre,KilSS,McI:tor,SteW:bub} for
further details. If an initial value $\xi\in\pksg$ gives rise to a
{\sc{cmc}} cylinder with extended frame $F_\lambda$ and monodromy
$M_\lambda$, and $\beta \in \{\lambda \in \C : 0 < |\lambda | < 1
\}$ is a point at which $M_\beta = \pm \mathbbm{1}$, we define a
simple factor
$$
    g = \begin{pmatrix}
    \sqrt{\frac{\lambda-\beta}{1-\bar{\beta}\,\lambda}} & 0 \\
    0 & \sqrt{\frac{1-\bar{\beta}\,\lambda}{\lambda-\beta}}
        \end{pmatrix}\,.
$$
Then for any $0 < r < |\beta |$ the dressed extended frame, obtained
from the $r$-Iwasawa factorization \cite{BurP:dre} of $g F_\lambda$,
is an extended frame of a {\sc{cmc}} cylinder with a bubbleton. On
the initial value level this dressing action corresponds to
$g\,\xi\,g^{-1}$ which obviously has singularities at
$\beta,\,1/\bar{\beta}$. To eliminate these, consider $\tilde{\xi} =
(\lambda-\beta)(1-\bar{\beta}\,\lambda)\,g\,\xi\,g^{-1}$. If $a =
-\lambda \det \xi$, then $\tilde{a} = -\lambda \det \tilde{\xi} =
(\lambda-\beta)^2(1-\bar{\beta}\,\lambda)^2\,a$, so the polynomial
$\tilde{a}$ of a bubbleton has double zeroes.
\begin{lemma}\label{thm:dense stratum}
If $a$ has multiple roots, then
$\mathcal{K}_a^\circ=\{\xi\in\mathcal{K}_a \,\mid
\mbox{ all roots of $\xi$ are unimodular }\}$
is open and dense in $\mathcal{K}_a$. If $a$
has no unimodular zeroes then it is a $g$-dimensional submanifold.
If $a$ has unimodular zeroes, then let $\tilde{a}$ denote the
quotient of $a$ by all real zeroes. Then $\mathcal{K}_a^\circ$ is
the image of the multiplication with an appropriate rational
function $p$ from $\mathcal{K}_{\tilde{a}}^\circ $ to
$\mathcal{K}_a^\circ$.
\end{lemma}
\begin{proof} Similar arguments as in the proof of
Lemma \ref{thm:compact} carry over to this situation.
\end{proof}
\begin{corollary} \label{thm:hom_S}
Suppose $a$ is polynomial of degree $2g$ satisfying the reality
condition \eqref{eq:a_reality}. Assume $a$ has precisely
$2\tilde{g}$ pairwise distinct non-unimodular roots and
$g-\tilde{g}$ pairs of unimodular roots of order 2. Then
$\mathcal{K}_a \cong (\Sp^1)^{\tilde{g}}$.
\end{corollary}
\begin{proof}
Since $-\det$ is the square of a norm on all skew-hermitian $2\times
2$ matrices, all $\xi\in\mathcal{K}_a$ have a zero at the unimodular
double roots of $a$. Let $a(\lambda)=\tilde{a}(\lambda)p^2(\lambda)$
be the corresponding decomposition of $a$ into an $\tilde{a}$ with
pairwise distinct roots and the corresponding factors
\eqref{eq:transformation}. Due to Proposition~\ref{th:pKf_min} the
one-to-one correspondence $\zeta\leftrightarrow\tilde{\zeta}$
between polynomial Killing fields $\zeta$ with roots and polynomial
Killing fields without roots $\tilde{\zeta}$ induces an isomorphism
$\mathcal{K}_a\simeq\mathcal{K}_{\tilde{a}}$. The assertion now
follows from Lemma \ref{thm:compact}.
\end{proof}
%
%
\subsection{Spectral curves II}
We also utilize the description of finite type {\sc{cmc}} surfaces
in $\Sp^3$ via spectral curves due to Hitchin \cite{Hit:tor}, and
relate this to our previous definition of spectral curves due to
Bobenko \cite{Bob:cmc}. While Hitchin defines the spectral curve as
the characteristic equation for the holonomy of a loop of flat
connections, Bobenko defines the spectral curve as the
characteristic equation of a polynomial Killing field. We shall use
both of these descriptions, and briefly recall their equivalence:
Due to \eqref{eq:monodromy}, the monodromy $\C^\ast \to \SL
(\C),\,\lambda \mapsto M_\lambda$ is a holomorphic map with
essential singularities at $\lambda = 0,\,\infty$. By construction
the monodromy takes values in $\SU$ for $|\lambda |=1$. The
monodromy depends on the choice of base point, but its conjugacy
class and hence eigenvalues $\mu_\lambda,\,\mu_\lambda^{-1}$ do not.
With $\Delta(\lambda)=\tr(M_\lambda)$ the characteristic equation
reads
\begin{equation} \label{eq:characteristic1}
    \mu_\lambda ^2 - \Delta(\lambda) \,\mu_\lambda + 1 =0\,.
\end{equation}
The set of solutions $(\lambda,\,\mu)\in\C^2$ of
\eqref{eq:characteristic1} yields another definition of the {\bf
spectral curve} of periodic (not necessarily finite type) solutions
of the $\sinh$-Gordon equation. Moreover, the eigenspace of
$M_\lambda$ depends holomorphically on $(\lambda,\,\mu)$ and defines
the {\bf eigenbundle} on the spectral curve. Let us compare this
with the previous definition of a spectral curve of periodic finite
type solutions of the $\sinh$-Gordon equations. Let $\zeta$ be a
polynomial Killing field with initial value $\xi\in\pk$, with period
$\tau$ so that $\zeta(p+\tau) = \zeta(p)$ for all $p \in \R^2$. Then
also the corresponding $\alpha(\zeta)$ is $\tau$-periodic. Let
$dF_\lambda= F_\lambda \alpha(\zeta),\,F_\lambda(0)=\mathbbm{1}$ and
$M_\lambda = F_\lambda(\tau)$ be the monodromy with respect to
$\tau$. Then for $z=0$ we have $\xi= \zeta(0) = \zeta(\tau) =
F_\lambda^{-1}(\tau) \,\xi\,F_\lambda(\tau) = M_\lambda^{-1} \xi
\,M_\lambda$ and thus
$$
    [\,M_\lambda,\,\xi\,] = 0\,.
$$
All eigenvalues of holomorphic $2\times 2$ matrix valued functions
depending on $\lambda\in\CP$ and commuting point wise with
$M_\lambda$ or $\xi$ define the sheaf of holomorphic functions of
the spectral curve. Hence the eigenvalues of $\xi$ and $M_\lambda$
are different functions on the same Riemann surface. Furthermore, on
this common spectral curve the eigenspaces of $M_\lambda$ and $\xi$
coincide point-wise. Consequently the holomorphic eigenbundles of
$M_\lambda$ and $\xi$ coincide.
\begin{proposition}
A finite type solution of the $\sinh$-Gordon equation is periodic if
and only if
\begin{enumerate}
    \item There exists a meromorphic differential $d\ln\mu$ on the
    spectral curve with second order poles without residues at the two
    points $\lambda=0$ and $\lambda=\infty$.
    \item This differential is the logarithmic derivative of a function
    $\mu$ on the spectral curve which transforms under the involutions
    \eqref{eq:involutions} as
    $\sigma^\ast\mu=\mu^{-1}$, $\rho^\ast\mu=\bar{\mu}^{-1}$ and
    $\eta^\ast\mu=\bar{\mu}$.
\end{enumerate}
Conversely, a periodic solutions of the $\sinh$-Gordon equation is
of finite type if and only if the monodromy \eqref{eq:monodromy}
fails at only finitely many points $\lambda\in\C^\times$ to be
semisimple.
\end{proposition}
\begin{proof} Due to Krichever \cite{Kri_77}, the
translations by $z\in\C$ act on the eigenbundle by the tensor
product with a one-dimensional subgroup of the Picard group. In
Sections~1.4-1.7 McIntosh \cite{McI:tor} describes this Krichever
construction for finite type solutions of the $\sinh$-Gordon
equation. The line bundle corresponding to $\tau\in\C$ is trivial if
and only if there exists a non-vanishing holomorphic function $\mu$
on the compactified spectral curve with essential singularities at
$\lambda=0$ and $\lambda=\infty$, whose logarithm has a first order
pole at $\lambda=0$ and $\lambda=\infty$ with singular part equal to
$\tau/\sqrt{\lambda}$ and $\bar{\tau}\sqrt{\lambda}$. This implies
the characterization of periodic finite type solutions.

At all simple roots of $\Delta^2-4$ the monodromy
\eqref{eq:monodromy} cannot be semisimple. Furthermore, at a double
root of $\Delta^2-4$ the monodromy fails to be semisimple, if and
only if it is dressed by a simple factor and contains a
corresponding bubbleton. An asymptotic analysis shows that there can
exists at most finitely many roots of $\Delta^2-4$ of order larger
than two.
\end{proof}
Pinkall and Sterling \cite{PinS}, and independently Hitchin
\cite{Hit:tor} proved that doubly periodic solutions of the
sinh-Gordon are of {\bf{finite type}}. Thus all metrics of
{\sc{cmc}} tori are of finite type. We enlarge this class by
relaxing one period, and make the following
\begin{definition}
The {\sc{cmc}} cylinders with constant Hopf differential and whose
metric is a periodic solution of finite type of the $\sinh$-Gordon
will be called {\bf{{\sc{cmc}} cylinders of finite type}}.
\end{definition}
\subsection{Examples} We compute some examples of initial values,
polynomial Killing fields and extended frames for spheres, and
spectral genus $g=0,\,1$ surfaces. Formulas for all finite type
surfaces in terms of theta-functions are given by Bobenko
\cite{Bob:tor}.
\subsubsection{Spheres} We start with a discussion of spheres. Since
the Hopf differential vanishes identically, spheres constitute a
degenerate case since their conformal factor is a solution to the
Liouville equation rather than the $\sinh$-Gordon equation, a fact
also reflected in the initial value which does not satisfy the
semi-simplicity condition. Consider
\begin{equation} \label{eq:sphere_alpha}
  \alpha_\lambda = \frac{1}{2}\begin{pmatrix}
  u_z\,dz-u_{\bar{z}}\,d\bar{z} &
  2\lambda^{-1}e^u\,dz \\
  -2\lambda\,e^ud\bar{z} &
  -u_z\,dz+u_{\bar{z}}\,d\bar{z}
  \end{pmatrix}\,.
\end{equation}
Then $2\,d\alpha_\lambda + [\,\alpha_\lambda \wedge
\alpha_\lambda\,] = 0$ if and only if $u$ solves the Liouville
equation $\partial \bar{\partial}\,u + e^{2u} = 0$. The solution is
$u(z,\,\bar{z}) = - \log(1+z\bar{z})$. Plugging this into the
$\alpha_\lambda$ in \eqref{eq:sphere_alpha}, and solving $dF_\lambda
= F_\lambda \,\alpha_\lambda,\,F_\lambda(0) = \mathbbm{1}$ gives the
extended framing
\begin{equation*}
  F_\lambda = \frac{1}{\sqrt{1+z\bar{z}}}\,
  \begin{pmatrix} 1 & z\lambda^{-1} \\ -\lambda\bar{z} & 1
    \end{pmatrix}\,.
\end{equation*}
Using the initial value $\xi = \lambda^{-1}\varepsilon_+ -
\lambda\,\varepsilon_-$, the corresponding polynomial Killing field
is
$$
\zeta = \frac{1}{1+z\bar{z}}\,\begin{pmatrix} z-\bar{z} & \lambda^{-1}(1+z^2) \\
-\lambda\,(1+\bar{z}^2) & \bar{z}-z \end{pmatrix}\,.
$$
%
%
\subsubsection{Flat cylinders} \label{sec:clifford}
We next discuss flat {\sc{cmc}} surfaces of revolution in
$\mathbb{S}^3$, and compute the closing conditions for the Clifford
torus.
\begin{proposition} \label{th:vaccuum}
After a possible isometry, for any flat surface $f: \R^2 \to
\mathbb{S}^3$ with constant mean curvature $H$, there exists a $t_0
\in \R$ such that $H= \cot(2t_0)$, and with $\lambda_0 = e^{i t_0}$
we have
\begin{equation} \label{eq:vaccuumF}
  f = F_{\lambda_0^{-1}}\,F_{\lambda_0}^{-1}
  \quad \mbox{ with }\quad
  F_\lambda = \exp \left( \,\frac{i}{2}
  \begin{pmatrix} 0 & z \lambda^{-1} + \bar{z}\\
    z  +  \bar{z}\lambda  & 0 \end{pmatrix} \,\right)\,.
\end{equation}
\end{proposition}
\begin{proof}
Solving $dF_\lambda = F_\lambda \alpha_\lambda$ with $u \equiv 0$ in
$\alpha_\lambda$ of \eqref{eq:general_alpha} proves that after a
possible isometry, any flat constant mean curvature immersion is
framed by such $F_\lambda$ as in \eqref{eq:vaccuumF}. Hence there
exist distinct $\lambda_0,\,\lambda_1 \in \mathbb{S}^1$, and a frame
of the given form such that after a possible isometry $f =
F_{\lambda_1}\,F_{\lambda_0}^{-1}$. By \eqref{eq:H_mu}, the two
distinct unimodular numbers $\lambda_j= e^{2i t_j},\,j=0,\,1$ must
be chosen so that for the constant mean curvature $H$ of $f$ we have
$H= \cot(t_0-t_1)$. Now $\lambda_0,\,\lambda_1$ are determined only
up to a phase, and since rotations in the $\lambda$-plane correspond
to rotations in the $z$-plane via a unitary gauge, we may adjust the
phase so that $\lambda_1 = \lambda_0^{-1}$.
\end{proof}
For such spectral genus zero surfaces, the polynomial Killing field
is constant and equal to the initial value $\xi = \tfrac{i}{2}\left(
(\lambda^{-1}+1)\,\varepsilon_+ + (1+\lambda)\,\varepsilon_-
\right)$. The factor $1/2$ is a consequence of the choices made in
Theorem \ref{thm:sinh}. For $\lambda_0 = i$ and $\lambda_1 = -i$ we
obtain a minimal flat surface. We determine the simple periods such
that the restriction to a fundamental domain gives the {\bf{Clifford
torus}}: The eigenvalues of $F_\lambda(z)$ are of the form $\exp(\pm
\mu)$ with $\mu(z,\,\lambda) = \tfrac{i}{2}(z\,\lambda^{-1/2} +
  \bar{z}\,\lambda^{1/2})$. Simple periods are numbers $\omega_1,\,\omega_2 \in \C$ with
smallest possible modulus satisfying $\mu(\omega_1,\,\pm i) = \pi i$
and $\mu(\omega_2,\,\pm i) = \pm \pi i$, and compute to $\omega_1 =
\pi \sqrt{2} $ and $\omega_2 =\pi i \sqrt{2} $.
\subsubsection{Delaunay surfaces}

For $a,\,b \in \R$ the 1-parameter family of conformal metrics of
Delaunay surfaces $v^2(x) (dx^2 + dy^2)$ is given by the Jacobian
elliptic functions $v(x) = 2b\, \mathrm{dn}(2bx\,|\,1-a^2/b^2)$. For
the derivation of this conformal factor from the Gauss equation in
the rotational case we refer to \cite{BurK,SKKR}. Note that in the
limiting cases $a=\pm b$ we have $v \equiv 1$, which is the flat
case, while when $a=0$ we have $v(x) = \mathrm{sech}(x)$. The
initial value \cite{Kil:del, SKKR} is $\xi = (a\lambda^{-1} +
b)\,i\varepsilon_+ + (b + a\lambda)\,i\varepsilon_-$, and the
polynomial Killing field computes to
$$
    \zeta = i\begin{pmatrix} -\tfrac{v'(x)}{2v(x)} &
    \tfrac{2ab\lambda^{-1}}{v(x)} +\tfrac{v(x)}{2}\\
    \tfrac{2ab\lambda}{v(x)} +\tfrac{v(x)}{2} & \tfrac{v'(x)}{2v(x)}
    \end{pmatrix}\,.
$$
%
%
\section{Deformation of spectral data}
\label{sec:deformation} In this section we describe the spectral
curves of periodic finite type solutions of the $\sinh$-Gordon
equation by two polynomials $a$ and $b$. The first polynomial
defines the hyperelliptic curve and the second polynomial the
meromorphic differential $d\ln\mu$ on this curve. Not all
polynomials $a$ and $b$ correspond to spectral curves of periodic
solutions of the $\sinh$-Gordon equation. In order to describe the
subsets of all such $a$ and $b$, we derive vector fields on the
space of coefficients of $a$ and $b$, which leave these subsets
invariant. We will be using the usual spectral parameter $\lambda$
in which $\lambda = 0,\,\infty$ are singularities, as well as the
transformed spectral parameter
$$
\kappa = i\,\frac{1-\lambda}{1+\lambda}\,.
$$
(Hence $\lambda = (i-\kappa)/(i+\kappa)$.) The part of the spectral
curve over $\{\lambda \in \C : |\lambda |=1 \}$ corresponds to real
$\kappa$. These parameters are fixed only up to M\"obius
transformations
\begin{align}\label{eq:moebius}
    \lambda &\mapsto e^{2i\varphi}\lambda&
    \kappa \mapsto \frac{\sin
    \varphi + \kappa \cos \varphi}{\cos \varphi - \kappa \sin \varphi}
\end{align}
This degree of freedom allows us to assume that no branch point or
otherwise significant point, which we shall introduce later (for
example a zero of $d\ln\mu$), lies at $\kappa = \infty$. The
spectral curve is then a hyperelliptic surface which we describe
with the equation
\begin{equation} \label{eq:def_a}
  \nu^2 = (\kappa^2 +1)\,a(\kappa)\,.
\end{equation}
Here $a$ is a real polynomial of degree $2g$ which has highest
coefficient equal to one, and which is non-negative for $\kappa \in
\R$. Thus $a$ only possesses real roots of even order.

In the following we will consider finite type {\sc{cmc}} cylinders
in $\Sp^3$. These possess a monodromy~\eqref{eq:monodromy}
whose eigenvalue $\mu$ is a
holomorphic function on the spectral curve, and has essential
singularities at $\kappa = \pm i$. Then $d\ln \mu$ is an abelian
differential of the second kind of the form
\begin{equation} \label{eq:def_b}
d\ln\mu = 2\pi i \frac{b(\kappa)\,d\kappa}{(\kappa^2+1)\,\nu}\,,
\end{equation}
where $b$ is a real polynomial of degree $g+1$.

As a consequence of the work of Bobenko \cite{Bob:tor, Bob:cmc}, our
starting point is the definition of what we call the spectral data
of a {\sc{cmc}} cylinder of finite type in the 3-sphere.
\begin{definition} \label{thm:spec_bobenko}
Let $a$ be a real polynomial of degree $2g$ with highest coefficient
equal to one, and let $b$ be a real polynomial of degree $g+1$, and
$\kappa_0,\,\kappa_1 \in \R$ two {\bf{marked points}}.

The {\bf{spectral data}} of a {\sc{cmc}} cylinder of finite type in
$\Sp^3$ with mean curvature
\begin{equation}\label{eq:H kappa}
H = \frac{1 + \kappa_0 \kappa_1}{\kappa_0 - \kappa_1}
\end{equation}
consists of a quadruple $(a,\,b,\,\kappa_0,\,\kappa_1)$ with the
following properties:
\begin{enumerate}
\item[\Con{A}] $a(\kappa) \geq 0$ for $\kappa \in \R$.
\item[\Con{B}] On the hyperelliptic surface $\nu^2 = (\kappa^2+1)\,a(\kappa)$
  there is a single valued holomorphic function $\mu$ with essential
  singularities at $\kappa = \pm i$ with logarithmic differential
  \eqref{eq:def_b}, that transforms under the three involutions
\begin{align*}
  \sigma:(\kappa,\,\nu)&\mapsto (\kappa,\,-\nu),&
  \rho:(\kappa,\,\nu)&\mapsto (\bar{\kappa},\,\bar{\nu}),&
  \eta:(\kappa,\,\nu)&\mapsto (\bar{\kappa},\,-\bar{\nu}),
\end{align*}
as
  $\sigma^{\ast}\mu = \mu^{-1}$, $\rho^{\ast}\mu=\bar{\mu}^{-1}$ and
  $\eta^{\ast}\mu=\bar{\mu}$.
\item[\Con{C}] $\mu(\kappa_0) = \mu(\kappa_1) = \pm 1$.
\end{enumerate}
\end{definition}
We call the fixed point set of $\rho$ the {\bf real part}. While a
M\"obius transformation \eqref{eq:moebius} of the parameter $\kappa$
changes the spectral data $(a,\,b,\,\kappa_0,\,\kappa_1)$, it
changes neither the corresponding periodic solutions of the
$\sinh$-Gordon equation nor the corresponding {\sc{cmc}} cylinders
in $\Sp^3$. Hence the moduli space of spectral data is the set of
equivalence classes of spectral data up to the action on the
spectral data induced by \eqref{eq:moebius}, prompting the following
\begin{definition}
{\rm{(i)}} For all $g\in\N_0$ let $\hat{\moduli}_g$ be the space of
equivalence classes of spectral data $(a,\,b)$ obeying
conditions~\Con{A}-\Con{B} up to the action of \eqref{eq:moebius} on
$(a,b)$.

{\rm{(ii)}} For all $g\in\N_0$ let $\moduli_g$ be the space of
equivalence classes of spectral data $(a,\,b,\,\kappa_0,\,\kappa_1)$
obeying conditions~\Con{A}-\Con{C} up to the action of
\eqref{eq:moebius} on $(a,b,\kappa_0,\kappa_1)$.
\end{definition}
Thus $\hat{\moduli}_g$ is the moduli space of spectral data of
periodic solutions of the $\sinh$-Gordon equation of arithmetic
genus $g$, and $\moduli_g$ is the moduli space of spectral data of
finite type {\sc{cmc}} cylinders in $\Sp^3$.

We now derive vector fields on open sets of spectral data $\left\{
(a,\,b,\,\kappa_0,\,\kappa_1)\right\}$ and show that their integral
curves are differentiable families of spectral data of periodic
finite type solutions of the $\sinh$-Gordon equation. We
parameterize such families by one or more real parameters, which we
will denote by $t$. From condition~\Con{B} in Definition
\ref{thm:spec_bobenko} we conclude that $\partial_t\ln \mu$ is
meromorphic on the corresponding family of spectral curves. If we
view these functions locally in terms of $\kappa$ and $t$, then
$\partial_t\ln \mu$ can only have poles at the branch points, or
equivalently at the zeroes of $a$, and at $\kappa = \pm i$. If we
assume that for such a family of spectral curves the genus $g$ is
constant, then $\partial_t \ln \mu$ can at most have poles of first
order at simple roots of $a$. In general we have
\begin{equation}\label{eq:def_c}
  \partial_t\ln \mu =  \tfrac{2 \pi i}{\nu} \,c(\kappa)
\end{equation}
with a real polynomial $c$ of degree at most $g+1$.

To compute the corresponding vector field on the space of
spectral data we view $\mu$ locally as a function of the parameters
$\kappa$ and $t$. Differentiating \eqref{eq:def_b} and
\eqref{eq:def_c} gives
\begin{equation*}
  \partial^2_{t\kappa}\ln \mu = 2\pi i
  \frac{2\dot{b}a - b\dot{a}}{2\nu^3}\,,\quad
  \partial^2_{\kappa t}\ln \mu = 2\pi
  i\frac{2(\kappa^2+1)ac' - 2\kappa ac -
  (\kappa^2+1)a'c}{2\nu^3}\,.
\end{equation*}
Second partial derivatives commute if and only if
\begin{equation} \label{eq:integrability_1}
  2\dot{b}a - b\dot{a} = 2(\kappa^2+1)ac' - 2\kappa ac -
  (\kappa^2+1)a'c\,.
\end{equation}
The highest coefficient on the right hand side vanishes, so both
sides are polynomials of at most degree $3g+1$. As the highest
coefficient of $a$ does not depend on $t$ we conclude that $\dot{a}$
is a real polynomial of degree $2g-1$, and $\dot{b}$ a real
polynomial of degree $g+1$. Thus we have to determine $3g+2$ real
coefficients. In case $a$ and $b$ have no common roots, equation
\eqref{eq:integrability_1} uniquely determines the values of
$\dot{a}$ at the roots of $a$ and the values of $\dot{b}$ at the
roots of $b$. Since the highest coefficient on the right hand side
depends only on the highest coefficient of $\dot{b}$, in this case
\eqref{eq:integrability_1} uniquely determines a tangent vector on
the space of spectral data of periodic finite type solutions of the
$\sinh$-Gordon equation. By defining such polynomials $c$ we obtain
vector fields on the space of real polynomials $a$ of degree $2g$
and highest coefficient one and real polynomials $b$ of degree
$g+1$.

For spectral data of {\sc{cmc}} cylinders in $\Sp^3$ we have to
deform in addition to the polynomials $a$ and $b$ the two marked
points, such that the closing condition~\Con{C} of
Definition~\ref{thm:spec_bobenko} is preserved. As long as $\kappa_0
\neq \kappa_1$, and thus $| H | <\infty$, we preserve the closing
condition if $\partial_t \ln \mu(\kappa_j(t),\,t) = 0$, which holds
precisely when $\partial_{\kappa_j} \ln
\mu(\kappa_j(t),\,t)\,\partial_t\kappa_j +
    \partial_t \ln \mu(\kappa_j(t),\,t) = 0$.
Using equations \eqref{eq:def_b} and \eqref{eq:def_c}, the closing
conditions are therefore preserved if and only if
\begin{equation}\label{eq:integrability_2}
    \dot{\kappa}_j = -\frac{(\kappa_j^2 +1)
    \,c(\kappa_j)}{b(\kappa_j)}\,.
\end{equation}
The equations \eqref{eq:integrability_1} and
\eqref{eq:integrability_2} define rational vector fields on the
space of spectral data $(a,b,\kappa_0,\kappa_1)$ not necessarily
obeying conditions~\Con{A}-\Con{C} of Definition~\ref{thm:spec_bobenko}
\begin{theorem}\label{thm:deformation}
{\rm{(i)}} Let $U$ be an open subset of spectral data $(a,\,b)$,
with $a,\,b$ having no common roots. Let $c$ be a smooth function
from $U$ to the real polynomials of degree $g+1$. Then equations
\eqref{eq:integrability_1} define a smooth vector field on $U$. The
corresponding flow leaves invariant the subset of spectral data
obeying conditions~\Con{A}-\Con{B} of
Definition~\ref{thm:spec_bobenko}.

{\rm{(ii)}} Let $U$ be an open subset of spectral data
$(a,\,b,\,\kappa_0,\,\kappa_1)$, with $a,\,b$ having no common
roots, and real $\kappa_0 \neq \kappa_1$ Let $c$ be a smooth
function from $U$ to the real polynomials of degree $g+1$. Then
equations \eqref{eq:integrability_1} and \eqref{eq:integrability_2}
define a smooth vector field on $U$. The corresponding flow leaves
invariant the subset of spectral data obeying
conditions~\Con{A}-\Con{C} of Definition~\ref{thm:spec_bobenko}.
\end{theorem}
\begin{proof}
We only prove (i). The proof of (ii) is similar.
The solutions $\dot{a}$ and $\dot{b}$ of
equation~\eqref{eq:integrability_1} are rational expressions of the
coefficients of $a$, $b$ and $c$. If $a$ and $b$ have no common roots,
then the Taylor coefficients of $\dot{a}$ and $\dot{b}$ at the roots
of $a$ and $b$ up to the order of the roots minus one, respectively,
and the highest coefficient of $\dot{b}$
are uniquely determined by equation~\eqref{eq:integrability_1}.
Hence, in this case the denominators of the rational expressions for
$\dot{a}$ and $\dot{b}$ do not vanish. Hence for smooth $c$ the
corresponding $\dot{a}$ and $\dot{b}$ are smooth too.

Due to \eqref{eq:def_c} $\partial_t\ln \mu$ is a meromorphic
function on the hyperelliptic curve. Hence the periods of the
meromorphic differential $d\ln\mu$ do not depend on $t$. The
transformation rules of $\mu$ under $\sigma$, $\rho$ and $\eta$ are
preserved under the flows of the vector field corresponding to $c$.
Hence the integrals of $d\ln\mu$ along any smooth path from one root
of $a$ to another root of $a$ is preserved too. This implies that
the subset of spectral data $(a,b)$, which determine by
\eqref{eq:def_b} a single valued function $\mu$ with
$\sigma^{\ast}\mu=\mu^{-1}$, $\rho^{\ast}\mu=\bar{\mu}^{-1}$ and
$\eta^{\ast}\mu=\bar{\mu}$ is preserved under this flow.
\end{proof}
\begin{remark}
The space of real hyperelliptic curves of genus $g$ is up to M\"obius
transformations is $2g-1$ dimensional. All such curves correspond to
real solutions of the $\sinh$-Gordon equation. The subset
$\hat{\moduli}_g$ of curves
corresponding to periodic solutions has codimension growing with $g$.
The space of real polynomials $c$ of degree at most $g+1$ is
$g+2$-dimensional. In case $c$ is proportional to $b$, then the
deformation corresponds to an infinitesimal M\"obius
transformation~\eqref{eq:moebius}. Hence $\hat{\moduli}_g$ is
$g+1$-dimensional.
\end{remark}
%
%
\section{Moduli of spectral data of periodic solutions of the $\sinh$-Gordon
equation}\label{sec:moduli}
In this section we switch to another description of the spectral
data of Definition~\ref{thm:spec_bobenko}. We shall use the values
of the trace $\Delta$ of the monodromy \eqref{eq:characteristic1} at
the roots of the derivative of $\Delta$, as in Grinevich and Schmidt
\cite{GriS1}, to determine local parameters on the moduli of
spectral curves of genus $g$. In doing so, and switching to the
parameter $\kappa$ of the previous section,  we describe the
twice-punctured Riemann sphere $\CP\setminus\{\pm i\}$, which is the
domain of definition of the covering map $\Delta$, as the result of
gluing infinitely many copies of $\C$ along cuts to be specified by
the branch points and combinatorial data. The combinatorial data
specifies which sheets are joined by which branch points and branch
cuts. We call this combinatorial data the {\bf{gluing rules}}.

The meromorphic function $\kappa$ and the holomorphic function $\mu$
with essential singularities at $\kappa=\pm i$ fulfill an equation
of the form
\begin{equation} \label{eq:characteristic}
  \mu^2 - \Delta(\kappa)\,\mu + 1=0\,.
\end{equation}
Hence $\Delta$ is a holomorphic function
\begin{equation}\label{eq:covering}
    \Delta\,:\,\CP\setminus\{\pm i\}\rightarrow \C\,, \qquad \kappa
    \mapsto\Delta(\kappa)
\end{equation}
The curve \eqref{eq:characteristic} is hyperelliptic with
hyperelliptic involution
$\sigma:(\kappa,\mu)\mapsto(\kappa,\mu^{-1})$. The branch points are
the odd-ordered roots of $\Delta^2 -4$. We next characterize those
maps $\Delta$ which correspond to spectral data $(a,\,b)$ described
in Definition~\ref{thm:spec_bobenko}.

Let us indicate to what extent the function $\Delta$ is determined
by its values at the branch points, and the gluing rules. We first
recall the simpler, but essentially comparable situation of finitely
sheeted coverings investigated by Hurwitz \cite{Hur}. Finitely
sheeted covering maps $\CP\rightarrow\CP$ are determined by their
branch points and their gluing rules up to M\"obius transformations
of both copies of $\CP$. In most cases the parametrization of the
image of the covering map is fixed by the values of this parameter
at the branch points. By fixing the values of the parameter of the
domain at some marked points, we can fix also the parameter of the
domain. For some fixed parametrization of the domain and the image
such a finitely sheeted covering map is a rational function.
Consequently a rational function can be characterized by the values
of the function at the zeroes of the derivative, some gluing rules,
and some conditions on the parametrization of the domain.

Now we return to our infinitely sheeted covering map
\eqref{eq:covering}. In this article we shall be concerned only with
spectral curves of finite geometric genus. In this case $\ln\mu$
extends to a meromorphic function on two neighbourhoods of
$\kappa=\pm i$ with first order poles at these two points. In
particular $(\ln\mu)^{-2}=-(\arccos(\Delta/2))^{-2}$ is a local
parameter on $\kappa\in\CP\setminus\{\pm i\}$ at $\kappa=\pm i$.
With this local parameter the domain of \eqref{eq:covering} can be
compactified to $\CP$. Since $\kappa$ takes at the two marked points
the values $\pm i$ and transforms under $\rho$ and $\eta$ as
$\kappa\mapsto\bar{\kappa}$, the parameter $\kappa$ is determined up
to M\"obius transformations \eqref{eq:moebius}.
\begin{theorem}\label{thm:Delta}
The functions $(\kappa,\,\mu)$ correspond to spectral data $(a,\,b)$
obeying conditions~\Con{A} and \Con{B} in
Definition~\ref{thm:spec_bobenko} if and only if the function
$\Delta$ satisfies the following conditions:
\begin{enumerate}
\item[\Con{D}] All but finitely many roots of $\Delta^2 -4$
  are roots of even order.
\item[\Con{E}] A branch of the function
  $(\ln\mu)^{-2}=-(\arccos(\Delta/2))^{-2}$ extends to a holomorphic
  function on two neighbourhoods of $\kappa=\pm i$ with simple zeroes
  at $\kappa=\pm i$.
\item[\Con{F}] The function $\Delta$ transforms as
  $\Delta(\bar{\kappa})=\bar{\Delta}(\kappa)$ and obeys
  $\Delta^2(\kappa)\leq 4$ for all $\kappa\in\mathbb{R}$.
\end{enumerate}
The roots of the corresponding polynomial $a$ have to be a subset of
the roots of $\Delta^2-4$, such that $(\Delta^2-4)/a$ has only roots
of even order. All such real $a$ are possible choices. Conversely,
the function $\Delta$ is uniquely determined by the spectral data
$(a,\,b)$.
\end{theorem}
\begin{proof}
The odd order roots of $\Delta^2-4$ are the branch points of the
two-sheeted covering \eqref{eq:characteristic} over $\kappa\in\CP$.
Hence only $\Delta$ obeying condition~\Con{D} correspond to
polynomials $a$. Since $\kappa=\pm i$ is a branch point,
condition~\Con{E} follows from condition~\Con{B}. Moreover, the
involutions $\rho$ and $\eta$ in condition~\Con{B} induce the
involution in Condition~\Con{F}.

Vice verse, Condition~\Con{F} implies that all real roots of
$\Delta^2-4$ are roots of
even order. Locally on $\kappa\in\CP\setminus\{\pm i\}$ the function
$\mu$ is holomorphic in terms of $\Delta$ and a square root of
$\Delta^2-4$. Hence for all real $a$, such that $(\Delta^2-4)/a$ has
only roots of even order, the differential
$d\ln\mu$ is a meromorphic differential
on the hyperelliptic curve defined by $\nu^2=(\kappa^2+1)a(\kappa)$.
This differential is antisymmetric with respect to the hyperelliptic
involution $\sigma:(\kappa,\nu)\mapsto(\kappa,-\nu)$. Furthermore,
due to condition~\Con{E} it has second order poles at $\kappa=\pm i$
without residues and no other poles. Hence it is of the form
\eqref{eq:def_b}. Due to condition~\Con{F} the involution
$\kappa\mapsto\bar{\kappa}$ induces on the hyperelliptic curve two
antilinear involutions $\rho$ and $\eta$, one of which denoted by
$\eta$ has no fixed points. The function $\mu$ transforms as
$\rho^{\ast}\mu=\bar{\mu}^{-1}$ and $\eta^{\ast}\mu=\bar{\mu}$.
\end{proof}
We shall see later that condition~\Con{E} implies \Con{D}. The
function $\Delta$ defines an infinitely sheeted covering map with
essential singularities at $\kappa=\pm i$. The branch points of this
covering map are the zeroes of $\Delta'(\kappa) =0$ and thus
precisely the roots of $d\ln \mu$ together with the set of
singularities of the spectral curve defined by equation
\eqref{eq:characteristic}.

We shall see that essentially we can move all the branch points
independently without destroying the periodicity. The reality
condition~\Con{F} in Theorem \ref{thm:Delta} imposes the only
restriction. Consequently the moduli space is a covering space over
the parameter space of the values of $\Delta$ at all branch points.
We shall see that for spectral curves of finite geometric genus all
of them with the exception of finitely many are fixed. Consequently
we can assume this parameter space to have finite dimension. The
gluing rules of this covering are completely determined by the
gluing rules of the covering \eqref{eq:covering}. In Hurwitz
\cite{Hur} the analogous deformations of finitely sheeted coverings
are investigated.
\begin{lemma} \label{thm:Delta 1}
For spectral data $(a,\,b)$ of genus $g$ obeying
conditions~\Con{A}-\Con{B} of Definition~\ref{thm:spec_bobenko}, for
which none of the pairwise distinct roots of $b$ are roots of $a$,
the element $[(a,\,b)]\in\hat{\moduli}_g$ possesses an open
neighbourhood in $\hat{\moduli}_g$, which is uniquely parameterized
by the values of $\Delta$ at the roots of $b$.
\end{lemma}
\begin{proof}
By assumption there exists for every root $\beta_i$ of $b$ a unique
polynomial $c_i$ that vanishes at every root of $b$ except at
$\beta_i$, where it attains a value such that
\begin{equation*}
  \partial_t \Delta(\beta_i)
  = 2\sinh(\ln(\mu(\beta_i)))\partial_t \ln\mu(\beta_i) = 1\,.
\end{equation*}
The corresponding vector fields in a neighbourhood $U$ of spectral
data in the moduli space commute: the map $U \to \C^{g+1},\,\beta_i
\mapsto \Delta(\beta_i)$ sends these vector fields to coordinate
vector fields. Since these vector fields are linearly independent on
$U$, they generate the quotient of all polynomials $c$ modulo the
polynomials $b$. Now the Lemma follows from
Theorem~\ref{thm:deformation}.
\end{proof}
In the following we shall investigate how the above parametrization
is affected if the roots of the polynomial $b$ are either not all
distinct, or some of the roots of $b$ coincide with the roots of
$a$. Recall that the local parameters are the values of $\Delta$ at
the branch points of the cover $\kappa \mapsto \Delta(\kappa)$.
Describing this cover by means of the branch cuts we also have a
description of the moduli space as a cover of the parameter space,
thus obtaining a global picture of the moduli space. Assume we are
given two branch points of the covering map $\kappa \mapsto
\Delta(\kappa)$ that connect a given sheet of the cover with two
different sheets. If in the process of a continuous deformation one
of these branch points circumvents the other branch point, then the
sheets which these branch points connect permute. Hence higher order
roots of $b$ are branch points of the moduli space.

If on the other hand, a root of $b$ coincides with a root of $a$,
then the derivative of $\kappa \mapsto \Delta(\kappa)$ has a higher
order root there and $\Delta^2(\kappa) =4$ there. Thus a root of $b$
coincides with a root of $\Delta'$ that arises from a singularity of
the spectral curve \eqref{eq:characteristic}. Hence also in this
case two branch points of the covering map $\kappa \mapsto
\Delta(\kappa)$ coalesce.

For any branch point of a covering map \eqref{eq:covering} we can
choose small open neighbourhoods $U$ and $V$ in the domain and in
the image respectively, such that the restriction $\Delta|_U$ of
$\Delta$ to $U$ is a finitely sheeted covering over $V$. Small
movements of the branch points of $\Delta|_U$ over a closed subset
$A\subset V$ do not change the restriction of $\Delta$ to the
pre-image of $V\setminus A$ under $\Delta|_U$. Moreover, if we
consider higher order branch points as coalescing first order branch
points, we can move all these first order branch points
independently over $A$, without changing the restrictions of
$\Delta$ to the pre-image of $V\setminus A$ under $\Delta|_U$. More
precisely, if $\beta$ is a branch point at $\Delta_0=\Delta(\beta)$
of order $k$, then $w=\sqrt[k+1]{\Delta-\Delta_0}$ is a local
parameter of the covering space on an open neighbourhood of $\beta$.
All polynomials $P(w)$ of degree $k+1$ with highest coefficient $1$
and small lower order coefficients describe small perturbations of
the covering map $w\mapsto\Delta=w^{k+1}+\Delta_0$. All of them are
covering maps with $k$ branch points, which are the zeroes of the
derivative $P'$. The parameter $w$ of such coverings is determined
up to M\"obius transformations $w\mapsto aw+b$ with
$a\in\C^{\times},b\in\C$. Hence the values of $\Delta$ at the $k$
branch points together with the gluing rules of these covering maps
$w\mapsto P(w)$ determine such $P$ uniquely up to a
reparametrization $\tilde{P}(w)=P(aw+b)$ with
$a\in\C^{\times},b\in\C$. Furthermore, all small values of $P$ at
the $k$ branch points are realized by such polynomials. By gluing
the deformed covering $\tilde{\Delta}|_U$ along the pre-image of
$V\setminus A$ under $\Delta|_U$ with the restriction of the
undeformed $\Delta$ to $\CP\setminus(\{i,-i\}\cup A)$ we obtain a
deformation $\tilde{\Delta}$ of the covering $\Delta$, as a covering
map from an abstract Riemann surface without fixed parameter
$\kappa$ onto $\C$ with fixed parameter $\Delta$.

Finite combinations of such deformations we call {\bf{local
deformations}}. If the movements of the branch points respect
condition~\Con{F} in Theorem~\ref{thm:Delta}, then we call the
corresponding deformations {\bf{real local deformations}}.

In Section~\ref{sec:deformation} we described deformations of
spectral data by real polynomials $c$ of degree $g+1$. We shall
calculate the corresponding deformation of \eqref{eq:covering}. If
we consider also $\Delta$ and $\mu$ as functions depending on
$\kappa$ and $t$, then we have due to \eqref{eq:def_b} and
\eqref{eq:def_c}
\begin{align*}
\Delta'=2\sinh(\ln\mu)(\ln\mu)'&= 4\pi
i\frac{\sinh(\ln\mu)b}{(\kappa^2+1)\nu}&
\dot{\Delta}=2\sinh(\ln\mu)\dot{(\ln\mu)}&= 4\pi
i\frac{\sinh(\ln\mu)c}{\nu}.
\end{align*}
Hence $\dot{\Delta}$ is given by
\begin{equation}\label{eq:dotDelta}
\dot{\Delta}=\frac{(\kappa^2+1)c}{b}\Delta'.
\end{equation}
\begin{lemma}\label{thm:Delta 2}
Local deformations preserve conditions~\Con{D}-\Con{E} in
Theorem~\ref{thm:Delta}. In particular, real local deformed
$\tilde{\Delta}$ corresponds to a parameterized covering map
\eqref{eq:covering}, but the parameter $\kappa$ is determined only
up to M\"obius transformations \eqref{eq:moebius}. Those real local
deformations $\tilde{\Delta}$ of the covering map $\Delta$
corresponding to spectral data $(a,b)$ obeying
conditions~\Con{A}-\Con{B} of Definition~\ref{thm:spec_bobenko},
which move only the branch points of $\Delta$ corresponding to roots
of $b$, correspond to unique
$[(\tilde{a},\,\tilde{b})]\in\hat{\moduli}_g$.
\end{lemma}
\begin{proof} Since local deformations move only finitely many branch
points, they preserve conditions~\Con{D}-\Con{E} in
Theorem~\ref{thm:Delta}. Due to condition~\Con{E} the covering space
corresponding to a local deformed $\tilde{\Delta}$ can be
compactified to $\CP$. A parameter $\kappa$, which takes at the two
added points the values $\pm i$ and transforms under
$\Delta\mapsto\bar{\Delta}$ as $\kappa\mapsto\bar{\kappa}$ is unique
up to \eqref{eq:moebius}. Due to this Theorem all real local
deformations $\tilde{\Delta}$ of $\Delta$ corresponding to spectral
data $(a,b)$ correspond also to spectral data
$(\tilde{a},\tilde{b})$. The branch points of $\tilde{\Delta}$ are
in one-to-one correspondence with the zeroes of $d\ln\mu$ and the
singularities of \eqref{eq:characteristic}. Hence a real local
deformation, which moves only the branch points corresponding to
zeroes of $b$ corresponds to deformed spectral data
$(\tilde{a},\tilde{b})$, which are unique up to the M\"obius
transformations~\eqref{eq:moebius}.
\end{proof}
Since the space of local deformations of the covering map
\eqref{eq:covering} are manifolds, this Lemma can be used to make
$\hat{\moduli}_g$ into a real $(g+1)$-dimensional manifold. In the
sequel we shall call those real local deformations, which can be
realized as combinations of continuous movements of single branch
points, {\bf{continuous deformations}}. It is not difficult to
introduce a topology on $\moduli_g$, such that these deformations
corresponds to continuous paths. In fact the space of real
polynomials of fixed degree have a natural topology. As a quotient
space of such spaces, $\moduli_g$ also has a natural topology.

We shall first compute the map $\kappa \mapsto \Delta(\kappa)$ for
the periodic solutions of the $\sinh$-Gordon equation of spectral
genus zero: Then $a \equiv 1$ and
\begin{equation}\label{eq:a_vac}
    \ln \mu = 2\pi i \frac{b_0\kappa -b_1}{\nu} \,.
\end{equation}
For the anti-linear fix point free involution $\eta$ of Definition
~\ref{thm:spec_bobenko} we have $\eta^{\ast} \overline{d\ln\mu} =
d\ln\mu$. As $b_0$ and $b_1$ are real we obtain
\begin{equation}\label{eq:b_vac}
    d\ln\mu = 2\pi i \tfrac{b_0 + b_1
    \kappa}{(\kappa^2 +1)\,\nu}\,.
\end{equation}
Then
\begin{equation}\label{eq:genus zero}
  \Delta(\kappa) = 2\cosh(\ln\mu) = 2 \cos\left( 2\pi \tfrac{b_0\kappa
  -b_1}{\nu}\right) = 2\cos\left(\sqrt{4\pi^2 \tfrac{b_0^2\kappa^2
  - 2b_0b_1\kappa +b_1^2}{\kappa^2 +1}} \right)
\end{equation}
which is the composition of the two maps $\kappa \mapsto\delta=
4\pi^2 \tfrac{b_0^2\kappa^2
  - 2b_0b_1\kappa +b_1^2}{\kappa^2 +1}$ and
$\delta \mapsto\Delta= 2\cos\sqrt{\delta}$. Of the second of these
maps we choose the cuts along the lines
\begin{equation*}
  \sqrt{\delta} \in \pi n +i\,\R \Longleftrightarrow
  \left\{ \begin{array}{ll}
  \Delta \in (-\infty,\,-2] & \mbox{ for $n$ odd}\,,\\
  \Delta \in [2,\,\infty) & \mbox{ for $n$ even}\,. \end{array} \right.
\end{equation*}
Besides the branch points that arise for each $n \in \Z^{\times}$ we
have two additional branch points of $\kappa \mapsto
\delta(\kappa)$, situated at $\delta =0$ and
$$
\kappa = -\tfrac{b_0}{b_1} \Longleftrightarrow \delta = 4\pi^2
\tfrac{b_0^4+2b_0^2b_1^2 + b_1^4}{b_0^2 + b_1^2} = 4\pi^2 (b_0^2 +
b_1^2)\,.
$$
Between these latter two points we cut along the line segment
joining them. The covering map $\kappa\mapsto
2\cos\sqrt{\pi^2\kappa^2/(\kappa^2+1)}$ has sheets numbered by two
copies of $\N$. We denote these sheets by $l^+$ and $l^-$ with $l\in
\N$. In order to describe the movements of the branch points it is
convenient to choose only parallel cuts on each sheet. Besides the
branch cut along the real part we choose on the sheets $l^\pm$
\begin{equation}\label{eq:branch point}\begin{aligned}
\mbox{for $l\geq 1$ a cut along }&\Delta\in(-1)^l2+ i\mathbb{R}^+_0
&\mbox{ connecting the sheets }&l^+\mbox{ and }(l+1)^+\\
\mbox{for $l\geq 1$ a cut along }&\Delta\in(-1)^l2- i\mathbb{R}^+_0
&\mbox{ connecting the sheets }&l^-\mbox{ and }(l+1)^-\\
\mbox{for $l> 1$ a cut along }&\Delta\in-(-1)^l2+ i\mathbb{R}^+_0
&\mbox{ connecting the sheets }&l^+\mbox{ and }(l-1)^+\\
\mbox{for $l> 1$ a cut along }&\Delta\in-(-1)^l2- i\mathbb{R}^+_0
&\mbox{ connecting the sheets }&l^-\mbox{ and }(l-1)^-\\
\end{aligned}\end{equation}
Hence besides the sheets $1^\pm$ each sheet has at $\Delta=2$ and
$\Delta=-2$ a branch point. Those sheets, whose labels have
exponents $+$, do not have branch points and cuts at small imaginary
values of $\Delta$. Those sheets, whose labels have exponent $-$, do
not have branch points and cuts at large imaginary values of
$\Delta$. Each sheet has exactly one branch point at $\Delta=-2$.
This branch point connects the sheet with another sheet, to which we
pass along large circles in the $\Delta$-plane from small imaginary
parts to small imaginary parts, if the label has exponent $+$ and
from large imaginary parts of $\Delta$ to large imaginary parts, if
the label has exponent $-$, respectively. More precisely, for all
$l\in\mathbb{N}$ the branch point at $\Delta=-2$ connects the sheet
with label
\begin{equation}\label{eq:label sheets}\begin{aligned}
(2l-1)^+\mbox{ with the sheet reached by traversing large circles }&2l-1
\mbox{ times anti-clockwise,}\\
(2l)^+\mbox{ with the sheet reached by traversing large circles }&2l-1
\mbox{ times clockwise,}\\
(2l-1)^-\mbox{ with the sheet reached by traversing large circles }&2l-1
\mbox{ times clockwise,}\\
(2l)^-\mbox{ with the sheet reached by traversing large circles }&2l-1
\mbox{ times anti-clockwise.}\\
\end{aligned}\end{equation}
The labels of the sheets are completely determined by
this rule. If we choose $b_0=1/2$ and $b_1 =0$, then we have an
additional real branch cut that joins the sheets $1^-$ and $1^+$
with two additional real branch points at $\Delta=\pm 2$. This is
the spectral data of the standard round {\sc{cmc}} cylinder in
$\mathbb{R}^3$. The spectral data \eqref{eq:genus zero} of the standard
round {\sc{cmc}} cylinders in $\Sp^3$ are obtained by moving the
real branch point at $\Delta=-2$ in $\Delta\in[-2,2]$ along the sheets
$2^\pm, 3^\pm,\ldots$ to and fro. In doing so, the real part contains
branch cuts connecting the sheets $1^+$ with $1^-$, $2^+$ with
$2^-,\ldots$ and $l^+$ with $l^-$. We introduce the following
class of covering maps $\kappa \mapsto \Delta(\kappa)$ connecting
the sheets $\left(l^-\right)_{l\in\N}$ and $\left(l^+\right)_{l\in\N}$
at branch cuts such that the following hold:
\begin{enumerate}
  \item[\Con{G}] There exists an $L \in \N$
      such that for all $l>L$ the sheets $l^{\pm}$ only have the
      branch points and branch cuts \eqref{eq:branch point}.
  \item[\Con{H}] Along large circles in the $\Delta \in \C$-plane
    we get the following sequence of sheets:
\begin{align*}
  \ldots,\,(2l+1)^+,\,(2l-1)^+,\,\ldots,
  \,3^+,\,&1^+,\,2^+,\,4^+,\,\dots,\,(2l)^+,\,(2l+2)^+,\,\ldots
\end{align*}
traversing anti-clockwise from small imaginary parts to small imaginary
parts, and
\begin{align*}
  \ldots,\,(2l+1)^-,\,(2l-1)^-,\,\ldots,
  \,3^-,\,&1^-,\,2^-,\,4^-,\,\dots,\,(2l)^-,\,(2l+2)^-,\,\ldots
\end{align*}
traversing clockwise from large imaginary parts to large imaginary
parts. Reversing the anti-clockwise and clockwise order gives the
same sequences in reverse order.
  \item[\Con{I}] There exists an $L\in\mathbb{N}$ such that for $l\geq
    L$ the total number of all branch points that join two sheets from
    the set $\{1^-,1^+,\,\ldots,\,l^-,l^+\}$ is equal to
    $2l$. Furthermore, $2l-2$ of the corresponding branch cuts
    terminate at $2l-2$ additional branch points at infinity.
  \item[\Con{J}] The configuration is invariant under
    $l^+ \mapsto l^-,\,l^-\mapsto l^+,\,\Delta \mapsto \bar{\Delta}$.
    Furthermore, the branch order of real branch points at
    $\Delta=\pm2$ is odd.
    The fixed points of the corresponding anti-linear involution is called
    real part and is a branch cut along $\Delta\in[-2,2]$ between sheets
    $l^+$ and $l^-$
  \item[\Con{K}] Besides the branch cuts along the real part all
    branch cuts run along unbounded lines parallel to the imaginary axis
    in direction to very large or small imaginary parts according to
    the superscript $\pm$ of the corresponding sheets. These branch cuts
    start either at branch points or at the real part.
\end{enumerate}
Due to \eqref{eq:label sheets} the second condition~\Con{H} fixes
the labeling of the sheets. Together with condition~\Con{G} it
ensures condition~\Con{E} in Theorem~\ref{thm:Delta}. With the
holomorphic coordinates of this setting the corresponding covering
can be compactified to a compact Riemann surface. Hence
condition~\Con{G} is the finite type condition. Condition~\Con{I}
ensures that the compactified covering space has genus zero. The
analogous condition for finite sheeted coverings $\CP\rightarrow\CP$
is that the branching order is two times the number of sheets minus
two \cite{Hur}. Condition~\Con{J} ensures the reality
condition~\Con{F} in Theorem~\ref{thm:Delta}. This endows the
spectral curve~\eqref{eq:characteristic} with the involutions
$\sigma$, $\rho$ and $\eta$, the last of which is without fixed
points. Finally condition~\Con{K} describes a choice of the
corresponding branch cuts. All branch points away from the real part
are the starting point of a unique branch cut running to very large or
very small imaginary parts. In general, the branch cuts, which start
at the real part can be moved along the real part without changing
$\Delta$. In fact, if $\Delta_0\in(-2,2)$ is the starting point of a
branch cut connecting the sheets $l^+$ and $m^+$, and if the sheets
$l^+$ and $l^-$ and the sheets $m^+$ and $m^-$ are connected along a
real branch cut running along some neighbourhood of
$\Delta_0\in(-2,2)$, then this $\Delta_0$ is no branch point.
Moreover the corresponding branch cut can be moved along the real
part until it reaches either on the sheets $l^\pm$ or on the sheets
$m^\pm$ a real branch point of $\Delta$ without changing $\Delta$.
\begin{lemma}\label{thm:characterization}
All $\Delta$ obeying conditions~\Con{D}-\Con{F} in
Theorem~\ref{thm:Delta} fulfill conditions~\Con{G}-\Con{J}.
Moreover, the branch cuts may be chosen as described in
condition~\Con{K}.
\end{lemma}
\begin{proof} For all $r>0$ only finitely many of the infinite
branch points of $\lambda\mapsto 2\cos(\sqrt{\lambda})$ belong to
$B(0,r)\subset\C$. Hence condition~\Con{E} of
Theorem~\ref{thm:Delta} implies condition~\Con{G}. For all
$n\in\mathbb{N}$ the absolute value of the derivative $\sin$ of
$\cos$ is not smaller than $1$ on the circles
$|\ln\mu|=(n+\frac{1}{2})\pi$. Hence, due to Rouche's Theorem, the
number of zeroes of $\Delta'$ in the complement of small discs
around $\kappa=\pm i$ with appropriate radius is the same as the
corresponding number of the map \eqref{eq:genus zero}. This argument
is a slight variation of the Counting Lemma~2 in Chapter~2 of
P\"oschel and Trubowitz \cite{PoeT}. Therefore condition~\Con{E} of
Theorem~\ref{thm:Delta} implies also condition~\Con{I}. In
particular, the values of $\Delta$ at the branch points are bounded.
Now condition~\Con{H} just fixes the labeling of the sheets.
Condition~\Con{F} in Theorem~\ref{thm:Delta} implies
condition~\Con{J}. Due to condition~\Con{I} at all non-real branch
points there ends an unbounded branch cut, which may be chosen
according to condition~\Con{K}. Moreover with the exception of two
real branch points all real branch points correspond to an
additional unbounded branch cut starting at the real part.
Condition~\Con{K} describes a choice of these unbounded branch cuts.
\end{proof}
In order to show the converse we shall show that we can deform all
$\Delta$ obeying conditions~\Con{G}-\Con{K} into the covering
\eqref{eq:genus zero} and use Lemma~\ref{thm:Delta 2}. It would be
natural to do this with decreasing geometric genus. In
Theorem~\ref{thm:spec_AE} we shall concentrate on deformations
decreasing the following number
\begin{equation}\label{eq:G}\begin{aligned}
G&=\mbox{\rm geometric genus}\quad + \sum\limits\ind{real
singularities}\delta-\mbox{invariant
(see e.g. \cite{Ser})}\\
&=\frac{1}{2}\#\{\mbox{non real roots of }d\ln\mu\}+ \#\{\mbox{\rm
real branch points of }\Delta\}-1.
\end{aligned}\end{equation}
\begin{lemma} \label{thm:connected}
Any covering map \eqref{eq:covering} obeying
conditions~\Con{G}-\Con{K} may be continuously deformed with
decreasing $G$ \eqref{eq:G} within this class into a covering map
\eqref{eq:genus zero}.
\end{lemma}
\begin{proof}
In a first step we deform without changing the genus any covering
map \eqref{eq:covering} obeying the five conditions~\Con{G}-\Con{K}
into a covering map, which has only branch points of first order.
Nearby branch points of order $k$ at $\Delta=\pm 2$ any deformation
of $\Delta\mp2=z^{k+1}$ into a polynomial $p(z)$ with at least
$\frac{k}{2}$ distinct double roots deforms the branch point at
$\Delta=\pm 2$ into simple branch points at $\Delta=\pm2$. A
combination of such {\bf{continuous deformations}} deforms all
branch points at $\Delta=\pm 2$ into simple branch points at
$\Delta=\pm 2$. The movement of all other branch points at
$\Delta\not=\pm 2$ in $\Delta\in\C\setminus\{-2,2\}$ does
not change the genus at all.

In a second step we increase the values of the real branch points at
$\Delta=-2$ and decrease the values of the real branch points at
$\Delta=2$ by moving them into $\Delta\in(-2,2)$. As a result all
real singularities of \eqref{eq:characteristic} are deformed into
real zeroes of $d\ln\mu$ and the geometric genus becomes equal to
$G$. In particular, $\Delta$ takes on the real part only values in
$(-2,2)$. In the remaining steps we shall decrease the geometric
genus.

In a third step we increase the values at those real branch points
on $\Delta\in(-2,\,2)$, which are on the real part local minima, and
decrease the values of $\Delta$ at those real branch points which
are on the real part local maxima. If two real branch points
coalesce, on $\Delta\in(-2,\,2)$, then we move them away from the
real part. We may continue to shrink the real part, until only two
real branch points of $\kappa\mapsto\Delta(\kappa)$ remain and
converge against each other. In the limit we would obtain a spectral
curve with two connected components. But we stop shortly before this
happens. Consequently the sheets of the covering \eqref{eq:covering}
divide into the two groups labeled by two copies of $\N$, which are
joined only by a small circle between the two remaining real branch
points. Since all other simple branch points do not belong to the
real part, they occur in complex conjugate pairs, which can be moved
in complex conjugate directions.

In a fourth step we show that we can move the complex conjugate
branch points over $\Delta\in\C\setminus\{\pm 2\}$ in such a way
that the genus becomes at most equal to one. Let $L$ denote the
minimum of all $l\in\mathbb{N}$, such that for all $k>l$ the sheets
$l^{\pm}$ contain the branch points \eqref{eq:branch point} and
possibly the two remaining real branch points, but no other branch
point. In the subsequent discussion we neglect these real branch
points, which in this step are not moved at all. Due to
condition~\Con{G} the number $L$ is finite. Now we claim that we may
inductively decrease this number $L$ until it is equal to $0$. In
order to avoid branch points on vertical branch cuts we move all
branch points at $\Delta\not=\pm 2$ to places with pairwise
different real values of $\Delta$. Moreover, we can achieve that the
real part of $\Delta$ takes at these simple branch points at
$\Delta\not=\pm 2$ pairwise different values in $(-2,\,2)$, which
simplifies the subsequent argument. Furthermore if the branch points
on sheets with exponents $+$ cross the lines with real part of
$\Delta$ equals to $\pm 2$ along negative imaginary values of
$\Delta$ and on sheets with exponents $-$ along positive imaginary
values of $\Delta$, the branch points do not cross the branch cuts
described in \eqref{eq:branch point}. Consequently the number $L$ is
preserved under this deformation. At $\Delta=(-1)^L2$ the sheets
with labels $(L+1)^{\pm}$ join a unique branch point with sheets
with labels not larger than $L$. Let $k^\pm$ denote the labels of
these sheets. These sheets have besides these branch points at
$\Delta=(-1)^L2$ only branch points connecting with sheets in
$\{1^-,1^+,\ldots,L^-,L^+\}$. Due to condition~\Con{H} they have at
least one other branch point. If $L$ is even and the sheets
$k^{\pm}$ contain more than two branch points, then we move all
branch points with the exception of those with minimal and maximal
real parts of $\Delta$ starting with the smaller real parts through
the vertical branch cut with the lowest real parts of $\Delta$. If
$L$ is odd, then we start with the branch points with larger real
parts of $\Delta$ and move these branch points through the branch
cut with the largest real part of $\Delta$. Finally the sheets with
label $k^{\pm}$ contain besides the branch point at $\Delta=(-1)^L2$
exactly one other branch point. If we move this branch point to
$\Delta=-(-1)^L 2$ the number $L$ decreases. This proves the claim.
Hence we can decrease the number $L$, until it is equal to $0$ and
arrive at a covering map \eqref{eq:covering} that corresponds to a
spectral curve of genus one.

In a fifth step we finally move one of the two remaining real branch
points away from the other real branch point several times along
$\Delta\in[-2,\,2]$ to and fro until it eventually reaches
$\Delta=2$ at the sheet with label one. In order to do so we
distinguish two cases: If the real branch points are located on
sheets with odd label, then we move the branch point with larger
value of $\Delta$. If both real branch points are located  on a
sheet with even label, then we move the branch point with smaller
value of $\Delta$. We thus obtain a covering map \eqref{eq:covering}
of the form~\eqref{eq:genus zero} corresponding to a spectral curve
of geometric genus zero.
\end{proof}
We remark, that it is also possible to deform all $\Delta$ obeying
conditions~\Con{G}-\Con{K} with decreasing geometric genus into
$\Delta$ of geometric genus zero, but the proof is more complicated.
Moreover, in Lemma~\ref{thm:2 real zeroes} and
Lemma~\ref{thm:continuous genus 1} we shall deform those $\Delta$
corresponding to finite type {\sc{cmc}} cylinders, and again
decrease $G$ \eqref{eq:G} instead of the geometric genus.
\begin{theorem} \label{thm:weights}
All spectral curves of real finite type periodic solutions of the
$\sinh$-Gordon equation correspond uniquely to covering maps that
satisfy conditions~\Con{G}-\Con{K} above. The geometric genus $g+1$
of the corresponding spectral curve is equal to the sum of weighted
branch orders with the following weights:
$$
\mbox{weight} = \left\{ \begin{array}{cll} \mbox{branch order}
&\mbox{
    at } \Delta \neq \pm 2, & \\
\tfrac{1}{2}(\mbox{branch order}) &\mbox{ at } \Delta = \pm 2
&\mbox{for
  even branch order}, \\
\tfrac{1}{2}(\mbox{branch order}-1) &\mbox{ at } \Delta = \pm 2
&\mbox{for
  odd branch order}. \end{array} \right.
$$
\end{theorem}
\begin{proof}
Lemma~\ref{thm:characterization} shows that all $\Delta$
\eqref{eq:covering} of periodic solutions of the $\sinh$-Gordon
equation of finite type fulfill \Con{G}-\Con{K}. Due to
Lemma~\ref{thm:connected} we can by suitable movements continuously
deform an arbitrary $\Delta$ \eqref{eq:covering} obeying
conditions~\Con{G}-\Con{K} into the family of covering maps
\eqref{eq:genus zero}, which corresponds to spectral curves of
geometric genus zero. In Lemma~\ref{thm:Delta 2} it is shown that
these deformations preserve those $\Delta$, which correspond to
periodic solutions of the $\sinh$-Gordon equation of finite type. The
form $d\ln\mu$ has $2g+2$ zeroes on the spectral curve. The formula
for the genus is obtained by computing the order of
the roots of $d\ln\mu$ on $\mu^2-\mu\Delta(\kappa) + 1 =0$.
\end{proof}
%
%
\section{One-sided Alexandrov embeddings in $\Sp^3$}
\label{sec:AE space}
In this section we consider one-sided Alexandrov embeddings of
general manifolds $N$ and $M$. Some statements apply only to
one-sided Alexandrov embeddings with constant mean curvature, but we
do not use the special properties of {\sc{cmc}} cylinders of finite
type. We provide sufficient conditions which allow us to perturb
surfaces which are one-sided Alexandrov embedded inside of a collar
of the unperturbed surfaces into surfaces which remain one-sided
Alexandrov embedded. For this purpose we have to ensure that the
surfaces have collars with depths uniformly bounded from below.
Making use of the fact that there are no complete stable minimal
surfaces in $\Sp^3$, we present in Lemma \ref{th:Rosenberg_Lemma} a
crucial technical result communicated to us by Harold Rosenberg
\cite{Ros:com}: If both principal curvatures of an Alexandrov
embedded {\sc{cmc}} surface are uniformly bounded, then the cut
locus function is bounded from below by a positive number.

We consider one-sided Alexandrov embedded cylinders in $\Sp^3$. In
the literature we only found the notion of Alexandrov embeddings for
compact domains on the one hand, and the concept of properly
Alexandrov embedded immersions from open manifolds into open
Riemannian manifolds on the other hand. Since we are interested in
immersions of open manifolds into the compact Riemannian manifold
$\Sp^3$, we make the following
\begin{definition} A {\bf{one-sided Alexandrov embedding}} in
$\Sp^3$ is a smooth immersion $f$ from a connected 3-manifold $N$
with connected boundary $M=\partial N$ to $\Sp^3$ with the following
properties:
\begin{enumerate}
  \item The mean curvature of $M$ in $\Sp^3$ with respect to the
  inward normal is non-negative everywhere.
  \item The manifold $N$ is complete with respect to the metric
  induced by $f$.
\end{enumerate}
An immersion $f:N\rightarrow\Sp^3$ obeying condition~(ii) is called
an {\bf{Alexandrov embedding}}.
\end{definition}
A fixed orientation of $\Sp^3$ induces on $N$ and $M=\partial N$ an
orientation. Conversely, if $M$ is endowed with an orientation, then
there exists a unique normal, which points inward to the side of $M$
in $\Sp^3$, which induces on the boundary $M$ the given orientation of
$M$. In this sense the orientation of $M$ determines the inner normal
of $N$.

For each point $p \in M$ of a hypersurface of a Riemannian manifold
$N$ there exists a unique arc length parameterized geodesic
$\gamma(p,\cdot)$ emanating from $p = \gamma(p,0)$ and going in the
direction of the inward normal at $p$. Such geodesics are called
{\bf{inward $M$-geodesics}} \cite{Heb}.

Let $\gamma(p,\cdot)$ be an inward $M$-geodesic. Points $q \in N$ in
the ambient manifold that are 'close to one side' of $M$ can thus be
uniquely parameterized by $(p,\,t)$ where $p \in M$ and $q =
\gamma(p,t)$ for some inward $M$-geodesic $\gamma(p,\cdot)$ and some
$t \in \R_0^+$. The value of $t$ is the geodesic distance of $q$ to
$M$. Extending the geodesic further into $N$ it might eventually
encounter a point past which $\gamma(p,t)$ no longer minimizes the
distance to $M$. Such a point is called a cut point. The cut locus
of $M$ in $N$ consists of the set of cut points along all inward
$M$-geodesics. We define the {\bf{cut locus function}} as the
geodesic distance of the cut point to $M$:
\begin{equation}\label{eq:cutlocus}
  c:M\rightarrow\R^+,\quad p\mapsto c(p),\quad\mbox{ such that }
\gamma(p,c(p))\mbox{ is the cut point.}
\end{equation}
If we want to stress the dependence on $f$ we decorate $\gamma$ and
$c$ with index $f$. A known fact from Riemannian Geometry
(\cite[Lemma~2.1]{Heb}) asserts that a cut point is either the first
focal point on an inward $M$-geodesic,  or is the intersection point
of two shortest inward $M$-geodesics of equal length.

For a one-sided Alexandrov embedding $f:N \to \Sp^3$ of a 3-manifold
$N$ with boundary $\partial N = M$, the inward $M$-geodesics give us
a parametrization of $N$, which we call {\bf{generalized cylinder
coordinates}}:
\begin{equation}\label{eq:cylinder}
  \gamma_f : \{(p,t)\in M\times\R\mid 0\leq t\leq c_f(p)\}
  \rightarrow N\,.
\end{equation}
The restriction of $\gamma_f$ to $\{(p,t)\in M\times\R\mid 0\leq t<
c_f(p)\}$ is a diffeomorphism onto the complement of the cut locus.
The cut locus is homeomorphic to the quotient space $M/\sim_f$ with
the following equivalence relation on $M$:
\begin{align*}
    p&\sim_f q&&\Longleftrightarrow&
    c_f(p)&=c_f(q)\quad\mbox{ and }&\gamma_f(p,c_f(p))&=\gamma_f(q,c_f(q))\\
    &&&\Longleftrightarrow&&&\gamma_f(p,c_f(p))&=\gamma_f(q,c_f(q)).
\end{align*}
For all $p\in M$ we denote the corresponding equivalence classes by
\begin{align}\label{eq:ec cutlocus}
    [p]_f&=\{q\in M\mid \gamma_f(p,c_f(p))=\gamma_f(q,c_f(q))\}.
\end{align}
\begin{lemma}\label{thm:unique embedded 1}
Let $f:N\rightarrow\Sp^3$ and $\tilde{f}:\tilde{N}\rightarrow\Sp^3$
be two one-sided Alexandrov embeddings. The orientation of $\Sp^3$
induces orientations on $N$, $\partial N$, $\tilde{N}$ and
$\partial\tilde{N}$. If the two oriented boundaries $\partial
N=M=\partial\tilde{N}$ and the two restrictions to the boundaries
$f|_M=\tilde{f}|_M$ coincide, then there exists a diffeomorphism
$\Psi:\tilde{N}\rightarrow N$, whose restriction to $M$ is the
identity map of $M$, such that $\tilde{f}=f\circ \Psi$.
\end{lemma}
\begin{proof}
The immersions $f$ and $\tilde{f}$ induce on $N$ and $\tilde{N}$ two
Riemannian metrics $g$ and $\tilde{g}$, which coincide on $M$. The
Riemannian metrics induce on $N$ and $\tilde{N}$ two metrics $d$ and
$\tilde{d}$. Let $c_f$ and $c_{\tilde{f}}$ denote the cut locus
functions of the submanifold $M$ in the two Riemannian manifolds
$(N,g)$ and $(\tilde{N},\tilde{g})$. The generalized cylinder
coordinates \eqref{eq:cylinder} define diffeomorphisms $\gamma_f$
and $\gamma_{\tilde{f}}$ of
$$
  L=\{(p,t)\in M\times\R\mid 0\leq t<\min\{c_f(p),c_{\tilde{f}}(p)\}\}
$$
onto open subsets of $N$ and $\tilde{N}$. First we claim that these
diffeomorphisms $\gamma_f$ and $\gamma_{\tilde{f}}$ together with
the metrics $d$ and $\tilde{d}$ induce on $L$ the same metrics. In
fact, the subset of $L \times L$ on which the metrics coincide is
open and closed and therefore all of $L\times L$. Consequently both
cut locus functions $c_f$ and $c_{\tilde{f}}$ coincide, and for all
points $p\in M$ the corresponding classes $[p]_f=[p]_{\tilde{f}}$
\eqref{eq:ec cutlocus} coincide too. Hence the diffeomeorphism
$\gamma_f\circ\gamma^{-1}_{\tilde{f}}$ from the complement of the
cut locus in $\tilde{N}$ onto the complement of the cut locus in $N$
extends to a homeomorphism $\Psi$ from $\Tilde{N}$ to $N$. By
definition of the cylinder coordinates the immersions
$f\circ\gamma_f$ and $\tilde{f}\circ\gamma_{\tilde{f}}$ from $L$
into $\Sp^3$ coincide. Hence we have $\tilde{f}=f\circ\Psi$. Since
$f$ and $\tilde{f}$ are immersions, $\Psi$ is a diffeomorphism form
$\tilde{N}$ onto $N$.
\end{proof}
\begin{definition}
We call an immersion $f:W\rightarrow\Sp^3$ of a connected 3-manifold
$W$ with connected boundary $V=\partial W$ a {\bf{local one-sided
Alexandrov embedding}}, if the following hold:
\begin{enumerate}
    \item The mean curvature of $V$ in $\Sp^3$ with respect to the
  inward normal is non-negative everywhere.
    \item All inward $V$-geodesics exist in $W$ until they reach
  the cut locus \eqref{eq:cutlocus}.
    \item The generalized cylinder coordinates
    $\gamma_f$ \eqref{eq:cylinder} are surjective.
\end{enumerate}
An immersion obeying conditions~(ii)-(iii) is called a {\bf{local
Alexandrov embedding}}.
\end{definition}
If $f:N\rightarrow\Sp^3$ is an Alexandrov embedding, and $V\subset
M=\partial N$ is an open subset, which contains for all $p\in V$ the
classes $[p]_f$ \eqref{eq:ec cutlocus}, then the restriction $f|_W$
of $f$ to
$$
    W=\{\gamma_f(p,t)\in N\mid p\in V\mbox{ and }0\leq t\leq c_f(p)\}
$$
is a local Alexandrov embedding. The proof of Lemma \ref{thm:unique
embedded 1} carries over to the following situation:
\begin{corollary}\label{thm:unique embedded 2}
Let $f:W\rightarrow\Sp^3$ and $\tilde{f}:\tilde{W}\rightarrow\Sp^3$
be two local Alexandrov embeddings. If the oriented boundaries
$\partial W=V=\partial\tilde{W}$ and $f|_V=\tilde{f}|_V$ coincide,
then there exists a diffeomeorphism $\Psi:\tilde{W}\rightarrow W$,
whose restriction to $V$ is the identity map of $V$, such that
$\tilde{f}=f\circ \Psi$.\qed
\end{corollary}
We shall prove that 'one-sided Alexandrov embeddedness' is an open
condition, which will allow us to study deformation families of
one-sided Alexandrov embeddings. The main tool is a general
perturbation technique of Alexandrov embeddings, which we call
{\bf{collar perturbation}}. We consider perturbations $\tilde{f}$ of
a given smooth immersion $f:M\rightarrow\Sp^3$, which are `small'
with respect to the $C^1$-topology on the space of immersions from
$M$ into $\Sp^3$. For this purpose we use the trivialization of the
tangent bundle $T\Sp^3 \cong \SU \times \su$ by left invariant
vector fields.
\begin{remark}
For sake of simplicity we denote the derivative considered as a
smooth function into $\su$ by $f'$ instead of $f^{-1}df$, and endow
$\su$ with the norm $\|X\|=(-\frac{1}{2}\tr(X^2))^{\frac{1}{2}}$.
\end{remark}
\begin{lemma} \label{thm:injectivity radius}
Let $f:M\rightarrow\Sp^3$ be an immersion inducing on the 2-manifold
$M$ a complete Riemannian metric. If the absolute values of both
principal curvatures are bounded by $\kappa_{\max}>0$, then for some
$r>0$ depending only on $\kappa_{\max}$ the exponential map $\exp_p$
is at all points $p\in M$ a diffeomorphism from $B(0,r)\subset T_pM$
onto an open neighbourhood of $p$. Furthermore, for all $\epsilon>0$
there exists a $\delta>0$ depending only on $\epsilon$ and
$\kappa_{\max}$, such that
$$\|(f\circ\exp_p)'(v)-(f\circ\exp_p)'(0)\|<\epsilon\quad
\mbox{ for all }p\in M\mbox{ and all }v\in B(0,\delta)\subset
T_pM.$$
\end{lemma}
\begin{proof}
Choose two radii $r,R\in(0,\pi)$ such that the following two
conditions hold:

(i) $r+R\leq\arctan(\kappa_{\max}^{-1})$.

(ii) Let $q^\pm$ be the
    centers of two spheres of radius $R$ in $\Sp^3$ touching each other
    at the center $p$ of an open ball $B(p,r)\subset\Sp^3$. All
    geodesics in $\Sp^3$ emanating from $q^+$ or $q^-$ can intersect
    transversally the sphere of radius $R$ around $q^-$ and $q^+$
    respectively, inside of $B(p,r)$.

For $R<\arctan(\kappa_{\max}^{-1})$ both conditions can be fulfilled
for small $r$. For all $p\in M$ there are exactly two points $q^\pm$
on the $M$-geodesic through $p$, whose distance to $p$ are equal to
$R$. Both spheres of radius $R$ and centers $q^\pm$ in $\Sp^3$ touch
$M$ at $p$. The absolute values of both principal curvatures are not
larger than the principal curvatures of all spheres in $\Sp^3$ with
radius smaller than $\arctan(\kappa_{\max}^{-1})$. Hence on all
geodesics of $M$ emanating from $p$ the distances to both centers
$q^\pm$ are monotone increasing, as long as these distances are
smaller than $\arctan(\kappa_{\max}^{-1})$. Due to condition (i)
this is the case for all points on the geodesic, whose geodesic
distance to $p$ is not larger than $r$. Hence the ball
$B(p,r)\subset M$ is mapped by $f$ into the complement of both balls
$B(q^\pm,R)\subset\Sp^3$. We parameterize the points of
$B(p,r)\subset M$ by the intersection points of the shortest
geodesics connecting these points with $q^\pm$ with the spheres of
radius $R$ around $q^\pm$, respectively. Due to the second
condition~(ii) these parameters are smooth. Hence $\exp_p$ is on
$B(0,r)\subset T_pM$ a diffeomorphism onto an open neighbourhood of
$p$. Since on all $M$-geodesics through $p$ the distances to both
centers $q^\pm$ is monotone increasing, the uniform estimate
$\|(f\circ\exp_p)'(v)-(f\circ\exp_p)'(0)\|<\epsilon$ holds for all
$p\in M$ and all sufficiently small $v\in T_pM$.
\end{proof}
\begin{corollary}\label{thm: cut locus function bound}
Let $f:N\rightarrow\Sp^3$ be a one-sided Alexandrov embedding with
uniform upper bound $\kappa\ind{max}$ on both principal curvatures.
Then there exist constants $0<\ubcl<\arctan(\kappa\ind{max}^{-1})$
and $L>0$ depending only on $\kappa\ind{max}$ such that for all
$p\in M=\partial N$ with $c_f(p)<\ubcl$ the set $[p]_f$ \eqref{eq:ec
cutlocus} contains exactly two points. Furthermore, the angle
between both corresponding inward $M$-geodesics at the cut locus
$\gamma_f(p,c_f(p))$ is larger than $\pi-L\cdot c_f(p)$.
\end{corollary}
\begin{proof}
All 2-spheres of radius $t \in (0,\frac{\pi}{2}]$ in $\Sp^3$ have
principal curvatures both equal to $\kappa =\cot(t)$ (for the
computations see e.g. Montiel and Ros \cite{MonR:alex}). For every
$p \in M$ there is a unique inward $M$-geodesic $\gamma(p,\cdot)$.
Evaluating these geodesics at some $t \in \R^+_0$ gives a
hypersurface
\begin{equation}\label{eq:hypersurface}
M_t = \bigcup_{p \in M} \gamma(p,t)
\end{equation}
in $N$. This hypersurface is smooth as long as $\gamma(p,t)$ does
not pass through a focal point. If $\kappa_1$ and $\kappa_2$ are the
two principal curvatures of $M$ in $p$ with respect to the inner
normal, then the corresponding principal curvatures of $M_t$ at
$\gamma(p,t)$ are given by $\cot(\arctan(\kappa_1^{-1})-t)$ and
$\cot(\arctan(\kappa_2^{-1})-t)$. The focal point on $\gamma(p,\cdot)$
is at the value of $t$ given by
\begin{equation}\label{eq:focalpoint}
  t\ind{foc} =
  \arctan\left((\max\{\kappa_1,\,\kappa_2\})^{-1}\right) \geq
  \arctan(\kappa\ind{max}^{-1}).
\end{equation}
At all $p\in M$ with $c_f(p)<\arctan(\kappa\ind{max}^{-1})$ the
cut locus $\gamma_f(p,c_f(p))$ cannot be a focal point. Hence $[p]_f$
has to contain at least one other element $q\in
M\setminus\{p\}$. Choose two radii $r,R$ as in the proof of
Lemma~\ref{thm:injectivity radius}. The two spheres in $\Sp^3$ of
radius $R$, which touch $M$ at $p$ and $q$ outside of $N$, must not
intersect each other at distances $<r$ to $p$ and $q$,
respectively. Hence there exist $\ubcl>0$ and $L>0$, such that for
$c_f(p)<\ubcl$ the angle between the inward $M$-geodesics
$\gamma_f(p,\cdot)$ and $\gamma_f(q,\cdot)$ at the cut locus is larger
than $\pi-L\cdot c_f(p)$. If $\ubcl$ is small enough, this implies
that any element of $[p]_f$ belongs to one of the balls $B(p,r)$ or
$B(q,r)$ in $M$. If $c_f(p)$ is smaller than $R$, then due to
Lemma~\ref{thm:injectivity radius} $[p]_f\cap B(p,r)=\{p\}$ and
$[p]_f\cap B(q,r)=\{q\}$.
\end{proof}
We thank Harold Rosenberg for communicating the following result to
us. This arose through discussions with Antonio Ros and Harold
Rosenberg, and makes use of a technique for constructing a global
non-negative, non-trivial Jacobi field that was recently also
employed by Meeks, Perez and Ros \cite{MeePR:limit}.
\begin{lemma}[Rosenberg's Lemma]\label{th:Rosenberg_Lemma}
Let $f:N \to \Sp^3$ be a one-sided Alexandrov embedding with
constant mean curvature, and principal curvatures bounded by
$\kappa\ind{max}>0$. Then the cut locus function is bounded from below
by $\arctan(\kappa\ind{max}^{-1})$.
\end{lemma}
\begin{proof}
The mean curvature of the hypersurface \eqref{eq:hypersurface}
\begin{equation}\label{eq:mean curvature}
  H(t) = \tfrac{1}{2}\left(
  \cot\left(\arctan(\kappa_1^{-1})-t\right) +
  \cot\left(\arctan(\kappa_2^{-1})-t\right) \right)
\end{equation}
is positive for all $t \in (0,\,t\ind{foc})$, and strictly monotone
increasing, since
$$
    H'(t) = \tfrac{1}{2}\left(\sin^{-2}(\arctan(\kappa_1^{-1})-t) +
    \sin^{-2}(\arctan(\kappa_2^{-1})-t)\right) > 0\,.
$$
Let $c_f$ denote the cut locus function \eqref{eq:cutlocus}. If
there exists a point $p \in M$ for which $c_f(p) <
\arctan(\kappa\ind{max}^{-1})\leq t\ind{foc}$, then two inward
$M$-geodesics $\gamma(p,\cdot),\,\gamma(q,\cdot)$ through $p,\,q \in
M$ respectively, have to intersect at a distance of $c_f(p)$ from
$M$, and thus $c_f(p) = c_f(q)$. Hence, if there exists a point $p
\in M$ with $c_f(p) <\arctan(\kappa\ind{max}^{-1})$ then $M_t$
intersects itself for a value of $t < \arctan(\kappa\ind{max}^{-1})$
over two points $p,\,q \in M$. Let
$$
    \lbcl_0 = \inf\{t\,\left| \right. M_t \mbox{ intersects over two
    points of } M\}.
$$
Since over all points $p \in M$ the mean curvature of $M_t$ is
positive for all $0<t<\arctan(\kappa\ind{max}^{-1})$ with respect to the
inner normals, the surfaces $M_t$ cannot touch each other from
different sides over two points for
$0<t<\arctan(\kappa\ind{max}^{-1})$.

A normal graph in $\Sp^3$ over a domain in a geodesic sphere $\Sp^2$
is the graph of the composition of the exponential map with a
section of the normal bundle on this domain, see Fornari,
deLira and Ripoll \cite{ForLR}. By Lemma \ref{thm:injectivity
radius} there exists a $r>0$ depending only on $\kappa\ind{max}$
such that all $p \in M$ have open neighbourhoods in $M$, which are
normal {\sc{cmc}} graphs over the ball $B(p,\,r)$ inside the unique
geodesic 2-sphere $\Sp^2$, which touches $M$ at $p$. Due to
Arzel\`{a}-Ascoli, and the a priori gradient bound from Proposition
4.1 in \cite{ForLR}, every bounded sequence of normal {\sc{cmc}}
graphs over $B(p,\,r) \subset \Sp^2$ has a convergent subsequence.
Now let $(p_k)_{k\in\N}$ be a sequence in $M$ with
$$
    \lim_{k \to\infty}c_f(p_k) = c_0 = \inf\left\{c_f(p) \mid p \in
    M \right\}\,.
$$
Then there exists a sequence $\Theta_k$ of isometries of $\Sp^3$
which transform each point $p_k$ into a fixed reference point $p_0
\in \Sp^3$, and the tangent plane of $M$ at $p_k$ into the tangent
plane of a fixed geodesic sphere $\Sp^2_{p_0}\subset \Sp^3$ which
contains $p_0$. This sequence of isometries transforms
neighbourhoods $U_k$ of $p_k\in M$ into normal {\sc{cmc}} graphs
$\Theta_k[U_k]$ over $B(p_0,\,r)$. By passing to a subsequence we
may achieve that these graphs converge to a normal {\sc{cmc}} graph
$U$ over $B(p_0,\,r) \subset \Sp^2_{p_0}$, which is tangent to
$\Sp^2_{p_0}$ at $p_0$. Due to Corollary \ref{thm: cut locus
function bound}, the sets $[p_k]_f$ contain besides $p_k$ another
point $q_k$ for large $k$. Furthermore, the sequence of isometries
$\Theta_k$ transforms the sequence of geodesic 2-spheres tangent to
$M$ at $q_k$ into a converging sequence of spheres with limit
$\Sp^2_{q_0}$. This sphere contains the limit
$q_0=\lim\Theta_k(q_k)$ with distance $\mathrm{dist}(p_0,\,q_0) =
2c_0$. The shortest geodesic connecting $p_0$ and $q_0$ intersects
orthogonally both geodesics spheres $\Sp^2_{p_0}$ and $\Sp^2_{q_0}$.
By Corollary~\ref{thm: cut locus function bound}, for large $k$ the
points $q_k$ have neighbourhoods $V_k$, whose transforms
$\Theta_k[V_k]$ are normal {\sc{cmc}} graphs over
$B(q_0,\,r)\subset\Sp^2_{q_0}$. By passing again to a subsequence
the normal {\sc{cmc}} graphs $\Theta_k[V_k]$ converge to a normal
{\sc{cmc}} graph $V$ tangent to $Sp^2_{q_0}$ at $q_0$.

The transformed inward $M$-geodesics nearby $p_k$ and $q_k$ converge
to normal geodesics of these two limiting {\sc{cmc}} surfaces $U$
and $V$ in $\Sp^3$. If we shift both limiting {\sc{cmc}} surfaces by
$c_0$ along their normal geodesics as in \eqref{eq:hypersurface} we
obtain two surfaces touching each other from different sides at the
limit of the transformed cut points
$\lim\Theta_k(\gamma_f(p_k,\,c(p_k)))=
\lim\Theta_k(\gamma_f(q_k,\,c(q_k)))$. Hence the shifted surfaces
cannot have positive mean curvature with respect to the inner
normal. This implies $\lbcl_0=0$ and $H=0$. In this case
$\Sp^2_{p_0}=\Sp^2_{q_0}$ and the two minimal surfaces $U$ and $V$
are tangent at $p_0=q_0$ and are graphs over the same domain
$B(p_0,\,r)$. By Hopf's maximum principle both limiting minimal
normal graphs $U$ and $V$ coincide.

In this case, due to Corollary~\ref{thm: cut locus function bound},
for large $k$ neighbourhoods $W_k\subset M$ of $q_k$ are normal
graphs over $U_k\subset M$. Let $\norm_k$ be the corresponding
sequence of non-negative sections of the inward normal bundle of the
sequence $\Theta_k[U_k]$ of minimal graphs over $B(p_0,\,r)$. Then
$\Theta_k[W_k]$ are normal graphs over the sequence $\Theta_k[U_k]$
of normal graphs over $B(p_0,r)$. Then
$$
    \lim_{k\to\infty} \frac{\norm_k}{\sup
    \left\{|\norm_k(q)| \mid q \in \Theta_k[U_k] \right\}}
$$
converges to a non-trivial non-negative Jacobi field on the limiting
minimal graph $U$ over $B(p_0,\,r)$. This argument is a slight
variation of an observation by Meeks, Perez and Ros \cite{MeePR:limit}.

Since the limits of both sequences of normal minimal graphs
coincide, the function $|\norm_k|$ converges on $B(p_0,r)$ uniformly
to zero. Hence we may repeat the same line of argument at all points
of the boundary of these minimal normal graphs $\Theta_k[U_k]$ over
$B(p_0,\,r)$. By passing to a diagonal subsequence like in the
Arzel\`{a}-Ascoli theorem we may extend $U$ to a complete minimal
surface in $\Sp^3$ with a non-negative Jacobi field which does not
vanish identically. (Note that due to unique continuation of Jacobi
fields \cite{Wol:elliptic} the limit of the re-scaled difference of
both normal graphs stays bounded on all connected compact subsets of
the complete minimal surface, whenever it is bounded on one compact
subset.) A non-trivial non-negative Jacobi field implies that the
spectrum of the Jacobi operator is non-negative
(Allegretto-Piepenbrink Theorem, see Fischer-Colbrie and Schoen
\cite{Fis-ColS} Theorem 1 and Davies \cite{Dav} Lemma 4.1).
Therefore the limiting complete minimal surface is stable. But by a
result of Fischer-Colbrie \cite{Fis-ColS} ( see also Corollary 3 in
Schoen \cite{Sch:stable} ) there exist no complete stable minimal
surfaces in $\Sp^3$. This contradicts $c_0 <
\arctan(\kappa\ind{max}^{-1})$.
\end{proof}
Later we apply the collar deformation to bounded open subsets
$W\subset N$ of one-sided Alexandrov embeddings
$f:N\rightarrow\Sp^3$. Hence we shall find bounded open subsets
$V\subset M=\partial N$, which contain for all $p\in V$ the classes
$[p]_f$ \eqref{eq:ec cutlocus}. For this purpose we need a
{\bf{chord-arc}} bound:
\begin{lemma}\label{thm:ec diameter bound}
Let $f:N\rightarrow\Sp^3$ be a one-sided Alexandrov embedding with
second fundamental form $\II$ with respect to the inner normal
$\mathfrak{N}$. If $\lbcl$ is a lower bound on the cut locus
functions $c_f$ \eqref{eq:cutlocus} and $\ubdh$ a uniform bound on
the covariant derivative of the second fundamental form:
\begin{equation}\label{eq:upper bound derivative of h}
|(\nabla_X\II)(X,X)|=|\nabla_X(\II(X,X))-2\II(\nabla_XX,X)|\leq\ubdh|X|^3
\quad\mbox{ for all }X\in TM,
\end{equation}
then there exists $C>0$ depending only on $\lbcl$ and $\ubdh$, such
that
$$\dist_N(p,q)\leq\dist_M(p,q)\leq C\dist_N(p,q)\quad
\mbox{ for all }p,q\in M.$$
\end{lemma}
\begin{proof}
For all $p,q\in M$ we have $\dist_N(p,q)\leq\dist_M(p,q)$. In
general, these distances do not coincide. We shall construct a path
from $p$ to $q$ of length at most $C\dist_N(p,q)$. Due to ( Rinow
\cite{Rin}, pages 172 and 141) the points $p$ and $q$ are joined by
a shortest path in $N$. In case this shortest path touches at some
points the boundary \cite[Theorem~1.]{AA} we decompose it into
pieces. The boundary points of a shortest path might have
accumulation points. But any point of a shortest path, which is not
a boundary point, belongs to a unique geodesic piece in $N$, which
has only two boundary points at both ends. Hence it suffices to
construct such a path for two points $p$ and $q$, which are
connected by a geodesic in $N$ with only two boundary points at both
ends.

Due to \cite[Lemma~2.1]{Heb} the cut locus function $c_f(p)$
\eqref{eq:cutlocus} is for all $p\in M$ not larger than the first
focal point $t\ind{foc}$ \eqref{eq:focalpoint}. Hence both principal
curvatures are uniformly bounded by $\kappa\ind{max}=\cot(\lbcl)$.
If $\dist_N(p,q)\leq\lbcl$, then
\begin{align*}
  \dist_N(\gamma_f(p,\dist_N(p,q)),q)&\geq\dist_N(p,q)\quad\mbox{ and }&
  \dist_N(\gamma_f(q,\dist_N(p,q)),p)&\geq\dist_N(p,q).
\end{align*}
Hence the angles between the geodesic connecting $p$ and $q$, and
the inward $M$-geodesics at $p$ and $q$ are not smaller than the
angles of the triangle in $\Sp^2$ with three sides of length
$\dist_N(p,q)$, which is larger than $\frac{\pi}{3}$. If we reduce
the length of $\dist_N(p,q)$ along the gradient flow of $\dist_N$ on
$M\times M$, then we obtain a path in $M\times M$ from $(p,q)$ to
the diagonal in $M\times M$ of length smaller than $2\dist_N(p,q)$.
This shows the claim for short distances $\dist_N(p,q)\leq\lbcl$. On
the other hand, all geodesics $\gamma$ in $N$ starting at $p$, which
do not meet $M$ for distances $d\in(0,\pi)$ meet each other at the
antipode of $p$. The exponential map of $N$ maps a unique half space
of $T_pN$ into $N$. If the pre-image of $B(p,\pi)\subset N$ with
respect to $\exp_p$ contains this half space, then due to Hopf's
maximum principle (see e.g. \cite{Esc}) $B(p,\pi)\subset M$ is a
geodesic sphere in $\Sp^3$ and the statement is obvious. Otherwise
there exists a geodesic starting at $p$, which touches $M$ for some
$t\in(0,\pi)$. In particular, if $\dist_N(p,q)\geq\pi$, then there
exists $\tilde{q}\in M$ with
$$
    \dist_N(p,\tilde{q}) < \pi\quad\mbox{ and }\quad
    \dist_N(p,q)\geq\dist_N(p,\tilde{q})+\dist_N(\tilde{q},q).
$$
Hence we may assume $\lbcl<\dist_N(p,q)<\pi$. Let $\chi_p,\,\chi_q
\in [0,\frac{\pi}{2}]$ denote the angles in $T_pN$ and $T_qN$
between the inward geodesic $\gamma$ connecting $p$ and $q$ and the
inward normal to $M$, respectively. In this proof we shall consider
smooth families of geodesics $\gamma$ connecting two smooth paths
$s\mapsto p(s)$ and $s\mapsto q(s)$ in $M$ parameterized by a real
parameter $s$. For fixed $s$ the geodesic is parameterized by the
real parameter $t$. The derivatives with respect to $s$ are denoted
by prime and the derivatives with respect to $t$ by dot. Let $|p'|$
and $|q'|$ denote the lengths of the tangent vectors $p'$ and $q'$
with respect to the Riemannian metric. The function
\begin{equation}\label{eq:gradient}
|\nabla\dist_N(p,q)|=
\sup\left\{\tfrac{|\dist_N'(p,q)|}{|p'|+|q'|}\mid p'\in T_pM,q'\in
T_qM\right\}= \max\{|\tan(\chi_p)|,|\tan(\chi_q)|\}
\end{equation}
is the length of the gradient of the function $\dist_N$ on $M\times
M$. The geodesic $\gamma$ extends in $\Sp^3$ to a closed geodesic.
For any tangent vector $(p',q')\in T_pM\times T_qM$ there exists a
Killing field $\vartheta$ on $\Sp^3$, which moves the closed
geodesic $\gamma$ in such a way, that the intersection points at $p$
and $q$ moves along $p'$ and $q'$, respectively. Conversely, all
Killing fields $\vartheta$ generate a one-dimensional group of
isometries of $\Sp^3$. Let $s\mapsto\gamma_\vartheta(s,\cdot)$
denote the corresponding family of geodesics and
$s\mapsto(p_{\vartheta}(s),q_{\vartheta}(s))$ the corresponding
intersection points with $M$. In order to proceed we need
\begin{lemma}\label{thm:second derivative}
For all $\kappa\ind{max}>0$ and $0<\lbcl<\pi$ there exist
$\epsilon,\delta>0$ and $0<s_0<\frac{3}{2}$ with the following
property: For all $(p,q)\in M\times M$ with $\lbcl\leq\dist_N(p,q)<\pi$ and
$\max\{\tan(\chi_p),\,\tan(\chi_q)\}\leq\epsilon$ there exists a non
trivial Killing field $\vartheta$, such that $d:s\mapsto
d(s)=\dist_N(p_{\vartheta}(s),q_{\vartheta}(s))$ obeys
\begin{align}\label{eq:second derivatives}
d'(s)&\leq 0& d''(s)&\leq-\delta\cos\left(\tfrac{d(s)}{2}\right)&
|p'_\vartheta(s)|+|q'_\vartheta(s)|&\leq 3\cos\left(\tfrac{d(s)}{2}\right)&
\mbox{ for all }&s\in[0,s_0].
\end{align}
\end{lemma}
\begin{proof} We shall construct a Killing field $\vartheta$ with the
desired properties, which rotates $\gamma$ around two antipodes of
$\gamma$. The corresponding rotated geodesics
$\gamma_\vartheta(s,\cdot)$ belong to a unique geodesic 2-sphere
$\Sp^2\subset\Sp^3$. The corresponding paths $s\mapsto p(s)$ and
$s\mapsto q(s)$ move along the intersection of this sphere $\Sp^2$
with $M$. Hence we can calculate all derivatives on this sphere.

We parameterize this sphere by the real parameter $s$ of the family
$s\mapsto\gamma_{\vartheta}(s,\cdot)$ of rotated geodesics, and the
real parameter $t$ of these geodesics. We choose the equator as the start
points $\gamma_{\vartheta}(s,0)$ with distance $\frac{\pi}{2}$ to
the rotation axis. The vector fields $\vartheta$ and the geodesic
vector field $\dot{\gamma}$ along the geodesics
$\gamma_\vartheta(s,\cdot)$ form an orthogonal base of the tangent
space $T\Sp^2$ of this sphere away from the zeroes of $\vartheta$.
The vector fields $\vartheta$ and $\dot{\gamma}$ have at $(s,t)$ the
scalar products
\begin{align*}
g(\vartheta,\vartheta)&=\cos^2(t)&
g(\vartheta,\dot{\gamma})&=0&
g(\dot{\gamma},\dot{\gamma})&=1.
\end{align*}
Since $\dot{\gamma}$ is a geodesic vector field the derivative
$\nabla_{\dot{\gamma}}\dot{\gamma}$ vanishes. Moreover, the mean
curvature of the integral curve of $\vartheta$ starting at $(s,t)$
is equal to $\tan(t)$. Therefore at $(s,t)$ we have
\begin{align*}
    \nabla_{\vartheta}\vartheta&
    =\cos^2(t)\tan(t)\dot{\gamma}
    =\cos(t)\sin(t)\dot{\gamma}\,,&
    \nabla_{\vartheta}\dot{\gamma}&
    =-\tan(t)\vartheta \,,\\
    \nabla_{\dot{\gamma}}\vartheta&
    =-\tan(t)\vartheta \,,&
    \nabla_{\dot{\gamma}}\dot{\gamma}& =0\,.
\end{align*}
We parameterize a neighbourhood of the geodesic from $p$ to $q$ in
such a way that the corresponding vector field $\dot{\gamma}$ points
inward to $N$ at $p$ and outward of $N$ at $q$, respectively. Let
$(s_p,t_p)$ and $(s_q,t_q)$ be the coordinates of $p$ and $q$,
respectively. The Killing field $\vartheta$ induces along the paths
$s\mapsto p(s)$ and $s\mapsto q(s)$ the vector fields
\begin{align*}
    p'&=\vartheta(p)-\dot{\gamma}(p)
    \frac{g(\mathfrak{N}(p),\vartheta)(p)}
     {g(\mathfrak{N}(p),\dot{\gamma}(p))}\quad\mbox{ and }&
    q'&=\vartheta(q)-\dot{\gamma}(q)
    \frac{g(\mathfrak{N}(q),\vartheta(q))}
     {g(\mathfrak{N}(q),\dot{\gamma}(q))}
\end{align*}
with lengths $|p'|=\vn\frac{|\cos(t_p)|}{\cos(\chi_p)}$ and
$|q'|=\vn\frac{|\cos(t_q)|}{\cos(\chi_q)}$. The derivative $d'$ is
equal to
$$
    d'=
    \frac{g(\mathfrak{N}(p),\,\vartheta(p))}
         {g(\mathfrak{N}(p),\,\dot{\gamma}(p))}-
    \frac{g(\mathfrak{N}(q),\,\vartheta(q))}
         {g(\mathfrak{N}(p),\,\dot{\gamma}(q))}=
    \frac{g(\mathfrak{N}(p),\,\vartheta(p))}
         {\cos(\chi_p)}-
    \frac{g(\mathfrak{N}(q),\,\vartheta(q))}
         {\cos(\chi_q)}.
$$
Along the paths $p$ and $q$ with $X=p'$ and $X=q'$, respectively, we
have at $(s,t)$
\begin{align*}
    \nabla_X\frac{g(\mathfrak{N},\,\vartheta)}
    {g(\mathfrak{N},\,\dot{\gamma})}&=
    \frac{g(\nabla_X\mathfrak{N},\,\vartheta)+
    g(\mathfrak{N},\,\nabla_X\vartheta)}
    {g(\mathfrak{N},\,\dot{\gamma})}-
    \frac{g(\mathfrak{N},\vartheta)(
    g(\nabla_X\mathfrak{N},\,\dot{\gamma})+
    g(\mathfrak{N},\,\nabla_X\dot{\gamma}))}
    {(g(\mathfrak{N},\,\dot{\gamma}))^2}\\
    &=\frac{g(\nabla_X\mathfrak{N},X)+g(\mathfrak{N},\nabla_X\vartheta)}
           {g(\mathfrak{N},\dot{\gamma})}-
      \frac{g(\mathfrak{N},\theta)
                g(\mathfrak{N},\nabla_X\dot{\gamma})}
           {g(\mathfrak{N},\dot{\gamma})^2}\\
    &=-\frac{\II(X,\,X)}{g(\mathfrak{N},\dot{\gamma})}
    +\cos(t)\sin(t)
    +2\tan(t)\left(\frac{g(\mathfrak{N},\,\vartheta)}
                        {g(\mathfrak{N},\,\dot{\gamma})}\right)^2.
\end{align*}
Hence the second derivative is equal to
\begin{align*}
d''&=-\frac{\II(p',\,p')}{\cos(\chi_p)}-\frac{\II(q',\,q')}{\cos(\chi_q)}
+\frac{\sin(2t_p)-\sin(2t_q)}{2}\\
   &+2\tan(t_p)\left(\frac{g(\mathfrak{N}(p),\,\vartheta(p))}
         {g(\mathfrak{N}(p),\,\dot{\gamma}(p))}\right)^2
-2\tan(t_q)\left(\frac{g(\mathfrak{N}(q),\,\vartheta(q))}
         {g(\mathfrak{N}(q),\,\dot{\gamma}(q))}\right)^2.
\end{align*}
If along the rotation of the geodesic for $s\in[0,s_0]$ the following
inequalities are satisfied
\begin{align}\label{eq:assumption 1}
-\tfrac{\pi}{2}\leq t_p&\leq 0,&
0\leq t_q&\leq\tfrac{\pi}{2},&
\tfrac{\lbcl}{2}&\leq t_q-t_p,\quad\mbox{ and }&
\max\left\{\tan(\chi_p),\tan(\chi_q)\right\}&\leq\tfrac{1}{\sqrt{3}},
\end{align}
then $\frac{\lbcl}{2}\leq d=t_q-t_p$,
$\min\{\cos(\chi_p),\cos(\chi_q)\}\geq\tfrac{\sqrt{3}}{2}$ and the
last two terms of $d''$ are bounded by
\begin{align*}
\left|\tan(t_p)\left(\frac{g(\mathfrak{N}(p),\,\vartheta(p))}
         {g(\mathfrak{N}(p),\,\dot{\gamma}(p))}\right)^2\right|&\leq
\sin(|t_p|)\cos(t_p)\tan^2(\chi_q)
\leq\frac{\sin(2|t_p|)}{2\cdot 3}\\
\left|\tan(t_q)\left(\frac{g(\mathfrak{N}(q),\,\vartheta(q))}
         {g(\mathfrak{N}(q),\,\dot{\gamma}(q))}\right)^2\right|&\leq
\sin(|t_p|)\cos(t_p)\tan^2(\chi_q)
\leq\frac{\sin(2|t_q|)}{2\cdot 3}.
\end{align*}
Due to
$\sin(2t_q)-\sin(2t_p) = 2\sin(t_q-t_p)\cos(t_p+t_q)$ and
$\cos(t_p)+\cos(t_q) = 2\cos(\tfrac{t_q-t_p}{2})\cos(\tfrac{t_p+t_q}{2})$
we arrive at
\begin{align*}
    d''(s)&\leq
    -\frac{\II(p',\,p')}{\cos(\chi_p)}-\frac{\II(q',\,q')}{\cos(\chi_q)}
    -\sin(d)\cos(t_p+t_q)\left(1-\tfrac{2}{3}\right),\\
    |p'|+|q'|&=\vn\frac{\cos(t_p)}{\cos(\chi_p)}
              +\vn\frac{\cos(t_p)}{\cos(\chi_p)}
    \leq\vn\tfrac{4}{\sqrt{3}}\cos\left(\tfrac{d}{2}\right)
    \cos\left(\tfrac{t_p+t_q}{2}\right)\leq3\vn\cos(\tfrac{d}{2}).
\end{align*}
The assumption \eqref{eq:assumption 1} implies $t_p\leq-d+\frac{\pi}{2}$ and
$d-\frac{\pi}{2}\leq t_q$. For $d\in[\frac{\pi}{2},\pi)$ we get
$\cos(t_p)\leq\sin(d)$ and $\cos(t_q)\leq\sin(d)$. For
$d\in[\frac{\lbcl}{2},\frac{\pi}{2})$ we use $\cos(t_p)\leq 1$ and
$\cos(t_q)\leq 1$ and obtain
\begin{align*}
\sin\left(\tfrac{\lbcl}{4}\right)\leq\tfrac{1}{2}\sin\left(\tfrac{\lbcl}{2}\right)
\leq\tfrac{3\sqrt{3}}{8}\sin\left(\tfrac{\lbcl}{2}\right)&
\leq\sin(d)\min\left\{\tfrac{\cos(\chi_p)}{|p'|^2},
                      \tfrac{\cos(\chi_q)}{|q'|^2}\right\}.
\end{align*}
The third inequality of \eqref{eq:second derivatives} implies
$2\sin(\frac{\lbcl}{4})\cos(\frac{d}{2})\leq\sin(d)$. To sum up, the
second and the third inequalities of \eqref{eq:second derivatives}
are implied by \eqref{eq:assumption 1} and
\begin{align}\label{eq:assumption 2}
    \delta&\leq\frac{1}{3}\sin\left(\tfrac{\lbcl}{4}\right)\cos(t_p+t_q),&
    -\frac{\II(p',\,p')}{|p'|^2}&
    \leq\tfrac{1}{2}\delta\quad\mbox{ and }&
    -\frac{\II(q',\,q')}{|q'|^2}&
    \leq\tfrac{1}{2}\delta.
\end{align}
We shall show first that there exists a vector field $\vartheta$
obeying at $s=0$
\begin{align*}
    \delta&\leq\frac{1}{6}\sin\left(\tfrac{\lbcl}{4}\right)\cos(t_p+t_q),&
    -\frac{\II(p',\,p')}{|p'|^2}&
    \leq\tfrac{1}{4}\delta\quad\mbox{ and }&
    -\frac{\II(q',\,q')}{|q'|^2}&
    \leq\tfrac{1}{4}\delta.
\end{align*}
The Killing field $\vartheta$ is uniquely determined by two choices:

Firstly the choice of a geodesic sphere $\Sp^2\subset\Sp^3$, which
contains the closed geodesic from $p$ to $q$.

Secondly a choice of the zeroes of $\vartheta$, or equivalently a
choice of the coordinates $t_p$ and $t_q$ with $t_q-t_p=d\mod \pi$.
We start with $t_q=-t_p=\frac{d}{2}$.

We first choose a 2-sphere $\Sp^2\subset\Sp^3$. This 2-sphere is
uniquely determined by a choice of a line in $T_p\Sp^2 \cap T_pM$
which includes $p'$, or equivalently a choice of a line in $T_q\Sp^2
\cap T_qM$ which includes $q'$. Since $f$ is a one-sided Alexandrov
embedding and both principal curvatures are uniformly bounded by
$\kappa\ind{max}$, the cone angles of those double cones in $T_pM$
and $T_qM$, on which $\II(X,X)\geq-\frac{1}{4}\delta|X|^2$, are not
smaller than $\tfrac{\pi}{2}+\Order(\delta)$. For sufficiently small
$\epsilon\geq\max\{\tan(\chi_p)\tan(,\chi_q)\}$ these double cones
correspond in the plane orthogonal to $\dot{\gamma}(p)$ in $T_pN$
and in the plane orthogonal to $\dot{\gamma}(q)$ in $T_qN$ to double
cones with cone angles not smaller than $\tfrac{\pi}{2}$. Hence the
intersection of both double cones is non empty.

In a second step we shall show that the inequalities
\eqref{eq:assumption 1} and \eqref{eq:assumption 2} are satisfied
for $s\in[0,\,s_0]$ with some $s_0>0$. Due to
Lemma~\ref{thm:injectivity radius} and the assumption
$\max\{\tan(\chi_p),\tan(\chi_q)\}\leq\epsilon$ there exists $s_0$
such that $t_p$ and $t_q$ do not reach the roots of $\vartheta$ for
$s\in[0,\,s_0]$. Since the derivatives of $\cos(\chi_p)$,
$\cos(\chi_q)$, $t_p$ and $t_q$ with respect to $s$ are uniformly
bounded, there exists $s_0>0$ such that the inequalities
\eqref{eq:assumption 1} and the first inequality of
\eqref{eq:assumption 2} are satisfied for $s\in[0,\,s_0]$. Due to
\eqref{eq:upper bound derivative of h} also the derivatives of
$\II(p',p')$ and $\II(q',q')$ are uniformly bounded. Hence there
exists $s_0>0$ only depending on $\lbcl$ and $\kappa\ind{max}$, such
that the second and the third inequality of \eqref{eq:second
derivatives} are satisfied for $s\in[0,\,s_0]$.

Finally we have to satisfy the first inequality of \eqref{eq:second
derivatives}. At the start point $s=0$ this is always the case for
one choice of the sign of $\vartheta$. Now the second inequality of
\eqref{eq:second derivatives} implies the first.
\end{proof}
\emph{Continuation of the proof of Lemma~\ref{thm:ec diameter bound}.}
By Lemma \ref{thm:second derivative} there exists
$\epsilon,\delta>0$ and $0<s_0<\frac{2}{3}$, such that for all all
$p,q\in M$ with $\lbcl\leq\dist_N(p,q)<\pi$ and
$\left|\nabla\dist_N(p,q)\right|\leq\epsilon$ there exists a Killing
field $\vartheta$, along which the length $d$ is reduced for $0\leq
s\leq s_0$. The inequality $-d'(s)\leq|p'(s)|+|q'(s)|\leq
3\cos(\frac{d(s)}{2})\leq 3\cos(\frac{d(s_0)}{2})$ implies
$\cos(\frac{d(0)}{2})\geq(1-\frac{3}{2}s_0)\cos(\frac{d(s_0)}{2})$.
Hence we get
\begin{align*}
\frac{d(0)-d(s_0)}{\dist_M(p(s_0),p(0))+\dist_M(q(s_0),q(0))}&\geq
\frac{\delta\frac{s_0^2}{2}\cos\left(\frac{d(0)}{2}\right)}
     {3s_0\cos\left(\frac{d(s_0)}{2}\right)}
\geq\delta\frac{s_0}{6}\left(1-\tfrac{3}{2}s_0\right).
\end{align*}
We apply this procedure again and again, until either
$\dist_N(p,q)<\lbcl$ or $\left|\nabla\dist_N(p,q)\right|>\epsilon$.
As long as the gradient~\eqref{eq:gradient}
$\left|\nabla\dist_N(p,q)\right|>\epsilon$, the same estimate holds
for the corresponding gradient on the Riemannian manifold $M\times
M$. Furthermore, the corresponding gradient flow reduces $\dist_N$
with monotonic decreasing
$\frac{\epsilon}{\sqrt{2}}\dist_M(p(s),p(0))+
\frac{\epsilon}{\sqrt{2}}\dist_M(q(s),q(0))+\dist_N(p(s),q(s))$. In
summary, all points $p,q\in M$ obey
$\dist_M(p,q)\leq\max\{2,\,\frac{\sqrt{2}}{\epsilon},\,
\frac{12}{\delta(2s_0-3s_0^2)}\}\dist_N(p,q)$.
\end{proof}
\begin{lemma}(Collar perturbation)\label{thm:collar deformation 1}
For given $\kappa\ind{max}>0$ and
$0\leq H\ind{max}\leq\arctan(\kappa\ind{max}^{-1})$ there exist
$\epsilon>0$ and $R>0$ with the following property:
If $f:N\rightarrow\Sp^3$ is a one-sided Alexandrov embedding with
constant mean curvature $0\leq H\leq H\ind{max}$
and principal curvatures bounded by $\kappa\ind{max}>0$,
$p\in M$ is some point, and
$\tilde{f}:M\rightarrow\Sp^3$
is an immersion with constant mean curvature $\tilde{H}\geq 0$ obeying
\begin{align}\label{eq:immersion bound}
    \dist(f(q),\,\tilde{f}(q)) &< \epsilon\mbox{ and}&
    \|f'(q)-\tilde{f}'(q)\| &< \epsilon
&\mbox{for all }q&\in B(p,R)\mbox{ and}&
    |H-\tilde{H}|&< \epsilon.
\end{align}
Then $\tilde{f}$ extends to a local one-sided Alexandrov embedding
$\tilde{f}:W\rightarrow\Sp^3$, whose boundary
$V=\partial W\subset B(p,R)$ is an open neighbourhood of $p$.
\end{lemma}
\begin{proof} Due to Rosenberg's Lemma~\ref{th:Rosenberg_Lemma}
$\lbcl=\arctan(\kappa\ind{max}^{-1})$ is a lower bound of the cut
locus function~\eqref{eq:cutlocus}. The generalized cylinder
coordinates \eqref{eq:cylinder} define a diffeomorphism $\gamma_f$
of $M\times[0,\lbcl)$ onto an open subset of $N$, which is a collar.
Any lower bound on the cut locus function~\eqref{eq:cutlocus} is
also a lower bound on the distance to the first focal point
\eqref{eq:focalpoint} on the inward $M$-geodesics. Since $f$ is a
one-sided Alexandrov embedding, the absolute values of the negative
principal curvatures are smaller than the positive principal
curvatures. Consequently, due to the formula~\eqref{eq:mean
curvature}, the distances to the first focal points on the outward
$M$-geodesics are not smaller than the distances to the first focal
points on the inward $M$-geodesics. Hence the normal variation
defines an immersion of $(-\lbcl,\lbcl)\times M$ into $\Sp^3$. In
particular, the induced metric makes $(-\lbcl,\lbcl)\times M$ into a
Riemannian manifold with constant sectional curvature equal to one.
For all elements of this manifold the cylinder coordinates, i.e. the
distances to $M\simeq\{0\}\times M$ and the nearest point in $M$,
are uniquely defined. Hence we can glue $(-\lbcl,\lbcl)\times M$
along $\gamma_f([0,\lbcl)\times M)$ to $N$, and obtain a larger
3-manifold $\hat{N}\supset N$ without boundary, such that the
generalized cylinder coordinates \eqref{eq:cylinder} extend to a
diffeomorphism $\hat{\gamma}_f: (-\lbcl,\lbcl)\times
M\rightarrow\hat{N}$ and the immersion $f:N\rightarrow\Sp^3$ extends
to an immersion $\hat{f}:\hat{N}\rightarrow\Sp^3$.

The balls with radius $\lbcl$ around all elements of $M\subset\hat{N}$
are isometric to open balls in $\Sp^3$ of radius $\lbcl$. For
$\epsilon<\lbcl$ all immersions $\tilde{f}$ obeying the first
inequality of \eqref{eq:immersion bound} can be represented uniquely
as the composition of a smooth section $B(p,R) \to T\hat{N}|_{B(p,R)}$
with $\hat{f}\circ\exp_{\hat{N}}$. Here $\exp_{\hat{N}}$ denotes the
exponential map of $\hat{N}$. Hence for all immersions $\tilde{f}$
obeying the first inequality of \eqref{eq:immersion bound} there
exists a smooth embedding $i:B(p,R)\hookrightarrow\hat{N}$, such
that $\tilde{f}$ is equal to $\hat{f}\circ i$. Due to
\eqref{eq:immersion bound} $\Phi$ induces on $B(p,R)$ a Riemannian
metric denoted by $\dist_O$ nearby the original metric. Hence there
exists a constant $C_1$ depending only on $\epsilon$ such that
$$
\dist_O(i(q),i(q'))\leq C_1\dist_M(q,q')
\quad\mbox{ for all }q,q'\in B(p,R).
$$
We denote the image of this embedding $i$ by $O$ as an oriented
submanifold of $\hat{N}$. We shall identify $\tilde{f}|_{B(p,R)}$ with
the immersion $\hat{f}|_O$. For all $q\in O$ let
$t\mapsto\gamma_{\tilde{f}}(q,t)$ denote the inward $O$-geodesics in
$\hat{N}$. Since $O$ is not complete we have to be careful with the
cut locus.

Due to Lemma~\ref{thm:injectivity radius} and Theorem~1.2 in
\cite{ForLR} there exists $r>0$ depending only on $\kappa\ind{max}$
and $H\ind{max}$, such that all $q\in M$ belong to open domains in
$M$, which are normal {\sc{cmc}} graphs over the ball $B(q,r)$ in
the unique geodesic sphere $\Sp^2\subset\Sp^3$, which touches $M$ at
$q$. Furthermore, for sufficiently small $\epsilon$ the same is true
for all $q\in\overline{B(p,R-r)}$ and the perturbed immersion
$\tilde{f}|_{B(p,R)}$. We conclude that at all $q\in M$
the second fundamental form $\II$ of $f$ obeys \eqref{eq:upper bound
derivative of h} with $\ubdh$ depending only on $\kappa\ind{max}$
and $H\ind{max}$. Moreover, for sufficiently small $\epsilon$ the
principal curvatures of the perturbation $\tilde{f}$ are bounded at
all $q\in\overline{B(p,R-r)}$ by a constant depending only on
$\kappa\ind{max}$ and $H\ind{max}$. Finally, again at all
$q\in\overline{B(p,R-r)}$ the second fundamental form $\tilde{\II}$
of $\tilde{f}$ obeys \eqref{eq:upper bound derivative of h} with a
constant $\ubdh$ depending only on $\kappa\ind{max}$ and
$H\ind{max}$. Due to Lemma~\ref{thm:ec diameter bound} there exists
a corresponding bound on the chord-arc ratio $C_2>0$ depending only
on $\kappa\ind{max}$ and $H\ind{max}$.

All cut locus functions of one-sided Alexandrov embeddings are
uniformly bounded from above by $\frac{\pi}{2}$, since otherwise a
sphere with negative principal curvatures touches $M$ inside of $N$
contradicting Hopf's maximum principle. Now we choose $R=C_1(C_2
2\pi+\epsilon+r)$ and denote $U=B(i(p),C_2 \pi+\epsilon)\subset O$.
For all $q\in U\subset O$, we have
$$
\left\{ q' \in O \mid \exists t \in [0,\,\tfrac{\pi}{2}] \mbox{ with
} \dist_{\hat{N}} (\gamma_{\tilde{f}}(q,\,t),\,q') \leq t \right\}
\subset \left\{ q' \in O \mid \dist_{\hat{N}}(q,\,q') \leq \pi
\right\} \subset B(q,\,C_2\pi)\,.
$$
Therefore, for all $q\in U$ the cut locus function $c_{\tilde{f}}$
is well defined. For all such $q\in U$ let $[q]_{\tilde{f}}$ denote
the set
$$
[q]_{\tilde{f}}=\left\{q'\in O\mid
\dist_{\hat{N}}\left(\gamma_{\tilde{f}}(q,c_{\tilde{f}}(q)),\,q' \right)=
c_{\tilde{f}}(q)\right\}.
$$
For any closed subset $A\subset O$ the set
$\{q\in U\mid [q]_{\tilde{f}}\cap A\not=\emptyset\}$ is a closed
subset of $U$. Hence $V=\{q\in U\mid [q]_{\tilde{f}}\subset U\}$ is an
open subset of $O$. Furthermore $W=\{\gamma_{\tilde{f}}(q,t)\mid q\in V
\mbox{ and }0\leq t\leq c_{\tilde{f}}(q)\}$ is a submanifold of
$\hat{N}$ with boundary $V$. By construction $\hat{f}_W$ is a local
one-sided Alexandrov embedding. By choice of $R$, $V$ is contained in
the image of $B(p,R)$ under $i$.
\end{proof}
\begin{corollary} \label{thm:collar deformation 2}
Let $f:N\rightarrow\Sp^3$ be a one-sided Alexandrov embedding with
constant mean curvature $H\leq H\ind{max}$ and principal curvatures
bounded by $\kappa\ind{max}$. Assume that
$\tilde{f}:M\rightarrow\Sp^3$ is an immersion of the boundary
$M=\partial N$ with constant mean curvature $\tilde{H}\geq 0$
obeying~\eqref{eq:immersion bound} for all $q\in M$ instead of all
$q\in B(p,R)$ with $\epsilon$ as in Lemma~\ref{thm:collar deformation 1}.
Then $\tilde{f}$ extends to a one-sided Alexandrov embedding from
$N$ to $\Sp^3$.
\end{corollary}
\begin{proof}
We apply Lemma~\ref{thm:collar deformation 1} to all $p\in M$ and
obtain a covering of $M$ by open subsets $V$, which are boundaries
of local one-sided Alexandrov embeddings. From \eqref{eq:immersion
bound} we deduced in the proof of Lemma~\ref{thm:collar deformation
1} a bound on both principal curvatures of $f$ and $\tilde{f}$. Due
to Lemma~\ref{thm:ec diameter bound} this implies a uniform bound
$C_2$ on the chord-arc ratio. The choice of the radius $R$ in
Lemma~\ref{thm:collar deformation 1} ensures that for any $p\in M$
the constructed local one-sided Alexandrov embedding nearby $p$ is
not affected by the immersion $\tilde{f}$ restricted to $M\setminus
B(p,R)$. The corresponding local one-sided Alexandrov embeddings fit
together to an Alexandrov embedding
$\tilde{f}:\tilde{N}\rightarrow\Sp^3$ of a submanifold
$\tilde{N}\subset\hat{N}$ with boundary $\partial\tilde{N}=M$. It
remains to show that $\tilde{N}$ is complete with respect to the
Riemannian metric induced by $\tilde{f}$. Due to Rosenberg's
Lemma~\ref{th:Rosenberg_Lemma} the cut locus function is uniformly
bounded by $\lbcl$ from below. For all $r<\lbcl$, the submanifolds
$$
    \{\gamma_f(p,t)\in \Hat{N}\mid p\in M\mbox{ and }
    -r\leq t\leq c_f(p)\}\subset \hat{N}
$$
are complete with respect to the Riemannian metric induced by
$\hat{f}$. By construction, the Riemannian manifold $\tilde{N}$ with
the metric induced by $\tilde{f}$ is a closed submanifold
of one of these complete submanifolds of $\hat{N}$, and therefore
also complete.
\end{proof}
%
%
\section{Deformation of {\sc{cmc}} cylinders in $\Sp^3$}
\label{sec:deformation of cylinders}
In this section we consider the subspace in the moduli space of
periodic finite type solutions of the $\sinh$-Gordon equation that
contain the spectral data of {\sc{cmc}} cylinders in $\Sp^3$, and
begin with an investigation of the subset of spectral data of
one-sided Alexandrov embedded {\sc{cmc}} cylinders. We shall
determine all continuous deformations of spectral data, which
preserve the one-sided Alexandrov embeddedness. In order to do this,
we combine the description of finite type {\sc{cmc}} cylinders in
$\Sp^3$ in terms of polynomial Killing fields with the investigation
of one-sided Alexandrov embeddings. From
Definition~\ref{thm:spec_bobenko} and \ref{thm:weights} we
immediately conclude
\begin{corollary}
A covering map $\kappa \mapsto \Delta(\kappa)$ satisfying conditions
\Con{G}-\Con{K} corresponds to the spectral data of
a {\sc{cmc}} cylinder if and only if
\begin{enumerate}
\item[\Con{L}] there are two branch points
$\kappa_0,\,\kappa_1$ on the real part (i.e the fix point set of
the involution $\rho$) at which $\Delta(\kappa_0) = \Delta(\kappa_1) =
\pm 2$.\qed
\end{enumerate}
\end{corollary}
From the point of view of the covering map \eqref{eq:covering} the
two distinguished points are two real branch points at $\Delta=\pm
2$, which are specified by the sheets they connect. A continuous
deformation of a covering map \eqref{eq:covering} preserves
condition~\Con{L}, if the two distinguished real branch points at
$\Delta=\pm 2$ are not moved and stay inside the real part. The
coordinates of the two distinguished points $\kappa_0,\,\kappa_1$
are given by a parametrization of the cover $\kappa \mapsto
\Delta(\kappa)$, which is unique up to M\"obius transformations
\eqref{eq:moebius}. As a consequence of conditions~\Con{G}-\Con{I}
the covering map is biholomorphic to $\CP$, and condition~\Con{J}
equips it with an antiholomorphic involution. In general it is
difficult to read off the coordinates from the parameters of the
moduli space, that is the values of $\Delta$ at the branch points.
An interesting parameter is the value of the mean curvature
$H=(1+\kappa_0\kappa_1)/(\kappa_0-\kappa_1)$.
\begin{lemma}\label{thm:H continuous}
For continuous deformations of covering maps \eqref{eq:covering}
preserving conditions~\Con{G}-\Con{L} the mean curvature depends
continuously on the deformation parameter.
\end{lemma}
\begin{proof}
We describe the compactification of the covering space of $\Delta$
by an open cover $U_0,\,\ldots,\,U_n$. Here $U_1,\,\ldots,\,U_n$ are
small disjoint discs around all branch point of $\Delta$ which are to
be moved, and $U_0$ is the complement of the union of smaller closed
discs $A_j\subset U_j$ such that each $A_j$ is still a neighbourhhod
of the corresponding branch point. Let $w_j$ be a holomorphic
coordinate on $U_j$, such that $\Delta-\Delta_j$ is a polynomial
$P_j(w_j)$ with respect to $w_j$. The coefficients of $P_j$ are the
parameters of the continuous deformations. For sufficiently small
values of these parameters the covering space is a compact Riemann
surface biholomorphic to $\CP$. We describe such a biholomorphic map
$\Phi$ by a meromorphic function of degree one. This is unique up to
M\"obius transformations.

We need to show that if we eliminate the freedom of M\"obius
transformations by imposing three additional conditions, that $\Phi$
depends continuously on the deformation parameter. It suffices to
show that a sequence of meromorphic functions $\Phi_k$ corresponding
to a convergent sequence of parameters converges to the limit $\Phi$
that corresponds to the limit of the parameters. For this it
suffices to show that any such sequence of meromorphic functions
contains a convergent subsequence, and that the limit is the unique
limit that corresponds to the limit of the parameters. Pick on each
open set $U_j$ a M\"obius transformation $M_j$ whose composition
with $\Phi$ maps the branch point of $\Delta$ contained in $U_j$ to
the origin, and with derivative at the branch point is equal to one,
and such that each one of the two points $\kappa=\pm i$ is mapped to
$\infty$. Since the space of schlicht functions is compact (see
Proposition~7.15 in Section~14.7 of Conway \cite{Con2}), a
subsequence of $M_j \circ \Phi_k$ converges to a biholomorphic map
$U_j \to \C$. Since the parameters converge, the limit is equal to
$M_j \circ \Phi$.
\end{proof}
The special situation
$H=\infty$ corresponding to {\sc{cmc}} cylinders in $\R^3$ has a
simple description. In the case of {\sc{cmc}} cylinders in $\R^3$,
the above correspondence via spectral curves also holds with one
notable exception. The closing condition for {\sc{cmc}} cylinders in
$\R^3$ require one value $\kappa_0 \in \R$ at which the monodromy is
equal to $\pm \mathbbm{1}$ and its derivative vanishes there
\cite{DorH:per}. Assume we have the spectral data for a {\sc{cmc}}
cylinder in $\Sp^3$ in which $\kappa_0$ and $\kappa_1$ coalesce.
Then the mean curvature is infinite and the cylinder is shrunk to a
point. A blow-up then yields a {\sc{cmc}} cylinder in $\R^3$, see
Umehara and Yamada \cite{UmeY:tori}.
\begin{corollary}
A covering map satisfying conditions~\Con{G}-\Con{K} corresponds to a
{\sc{cmc}} cylinder in $\R^3$ if and only if there is a real branch
point of order at least 3.\qed
\end{corollary}
Let us first show that all one-sided Alexandrov embedded {\sc{cmc}}
cylinders in $\Sp^3$ of finite type obey the assumptions of
Rosenberg's Lemma~\ref{th:Rosenberg_Lemma} and the chord-arc bound
in Lemma~\ref{thm:ec diameter bound}.
\begin{lemma}\label{thm:killing bound}
For every compact subset $\mathcal{K}\subset\pksg$ and
$H\ind{max}>0$ there exist constants $\kappa\ind{max}>0$ and
$\ubdh>0$ with the following property: If
$f:\mathbb{R}^2\rightarrow\Sp^3$ is an immersion of finite type with
constant mean curvature $|H|\leq H\ind{max}$, whose polynomial
Killing field $\zeta$ takes at some point $p\in\mathbb{R}^2$ a value
$\zeta(p)\in \mathcal{K}$. Then the absolute values of both
principal curvatures of $f$ are uniformly bounded by
$\kappa\ind{max}$, and the second fundamental form $\II$ of $f$
satisfies \eqref{eq:upper bound derivative of h}.
\end{lemma}
\begin{proof}
The coefficients of $a(\lambda)=-\lambda\det(\xi(\lambda))$ depend
continuously on $\xi\in\pk$. Hence on every compact subset
$\mathcal{K}\subset\pksg$ these coefficients are uniformly bounded.
In particular, the union $\hat{\mathcal{K}}$ of all isospectral sets
$\mathcal{K}_a$ of Lemma~\ref{thm:compact} not disjoint from
$\mathcal{K}$ is compact too. If $f:\mathbb{R}^2\rightarrow\Sp^3$ is
a {\sc{cmc}} immersion of finite type, whose polynomial Killing
field $\zeta$ takes at some $p\in\mathbb{R}^2$ a value in
$\mathcal{K}$, then the other values belong to the corresponding
isospectral set. Since all derivatives of $f$ at some
$p\in\mathbb{R}^2$ depend continuously on $\zeta(p)$, the lemma
follows.
\end{proof}
\begin{lemma}\label{thm:immersion bound}
For all compact subsets $\mathcal{K}\subset\pksg$, $\epsilon>0$, $R>0$
and $H\ind{max}>0$ there exists $\delta>0$ with
the following property:
If $(\zeta,\,\lambda_0,\,\lambda_1)$ and
$(\tilde{\zeta},\,\tilde{\lambda}_0,\,\tilde{\lambda}_1)$ are
polynomial Killing fields and marked points of two
{\sc{cmc}} cylinders $f,\,\tilde{f}:M\rightarrow\Sp^3$ obeying
at some point $p\in M$
\begin{align}
\label{eq:killing bound 1}
\zeta(p)\in\mathcal{K}\mbox{ and }
|H|=\left|\frac{\lambda_0+\lambda_1}
    {\lambda_0-\lambda_1}\right| &\leq H\ind{max}
\quad\mbox{ or }&
\tilde{\zeta}(p)\in\mathcal{K}\mbox{ and }
|\tilde{H}|=\left|\frac{\tilde{\lambda}_0+\tilde{\lambda}_1}
    {\tilde{\lambda}_0-\tilde{\lambda}_1}\right| &\leq H\ind{max}\\
\label{eq:killing bound 2}
    \|\zeta(p)-\tilde{\zeta}(p)\|&<\delta\,,&
    |\lambda_0-\tilde{\lambda}_0|<\delta\mbox{ and }
    |\lambda_1-\tilde{\lambda}_1|&<\delta.
\end{align}
Then the corresponding immersions $f$ and $\tilde{f}$ satisfy
\eqref{eq:immersion bound} on $B(p,\,R)\subset M$.
\end{lemma}
\begin{proof}
Since $\pksg$ is an open subset of $\pk$, there exists a $\delta>0$
and a compact subset of $\pksg$, containing all balls
$B(\xi,\delta)$ with $\xi\in\mathcal{K}$. Furthermore, due to
Lemma~\ref{thm:killing bound} there exists also a compact subset
$\hat{\mathcal{K}}\subset\pksg$ containing all isospectral sets
$\mathcal{K}_a$ not disjoint from
$\cup_{\xi\in\mathcal{K}}B(\xi,\delta)$. Due to
Lemma~\ref{thm:killing bound} for given $\mathcal{K}\subset\pksg$ and
$H\ind{max}>0$ there exists $\delta>0$, $\kappa\ind{max}>0$ and
$\ubdh>0$, such that the conditions~\eqref{eq:killing bound 1} and
\eqref{eq:killing bound 2} imply that the corresponding immersions
$f$ and $\tilde{f}$ have principal curvatures bounded by
$\kappa\ind{max}$ and second fundamental forms with covariant
derivatives bounded by \eqref{eq:upper bound derivative of h}. In
particular there exists $C>0$ such that the Riemannian metrics
induced by $f$ and $\tilde{f}$ on $M$ are uniformly bounded by $C$
times the Euclidean metric on the conformal parametrization
$\mathbb{R}^2$. The abelian group $\mathbb{R}^2$ of spatial
translations  acts continuously on all elements of $\xi\in\pksg$.
Therefore all elements of the closed ball
$B(0,CR)\subset\mathbb{R}^2$ act uniformly continuous on $\hat{K}$.
Hence for all $\tilde{\epsilon}>0$ there exists $\delta>0$, such
that the inequalities \eqref{eq:killing bound 2} imply
$\|\zeta(q)-\tilde{\zeta}(q)\|<\tilde{\epsilon}$ for all $q\in
B(p,R)$. Since all derivatives of $f(q)$ and $\tilde{f}(q)$ depend
continuously on $\zeta(q)$ and $\tilde{\zeta}(q)$, respectively, the
lemma follows.
\end{proof}
\begin{lemma}\label{thm:convergency}
For all compact subsets $\mathcal{K}\subset\pksg$ and $H\ind{max}>0$
there exists $\delta>0$ with the following property: Let
$\tilde{f}:M\rightarrow\Sp^3$ be a {\sc{cmc}} cylinder with
non-negative mean curvature $\tilde{H}$ (with respect to normal
uniquely determined by the orientation of $M$) and with polynomial
Killing field $\tilde{\zeta}$ and marked points $\tilde{\lambda}_0$
and $\tilde{\lambda}_1$. If for all $p\in M$ there exists a
one-sided Alexandrov embedding $f:N\rightarrow\Sp^3$ with polynomial
Killing field $\zeta$ and marked points $\lambda_0$ and $\lambda_1$
obeying conditions ~\eqref{eq:killing bound 1} and \eqref{eq:killing
bound 2}, then $\tilde{f}$ extends to a one-sided Alexandrov
embedding $\tilde{f}:\tilde{N}\rightarrow\Sp^3$.
\end{lemma}
\begin{proof}
As in the proof of Lemma~\ref{thm:immersion bound} there exists
$\delta>0$ and a compact subset $\hat{\mathcal{K}}\subset\pksg$
containing all isospectral sets not disjoint from
$\cup_{\xi\in\mathcal{K}}B(\xi,\delta)$. Due to
Lemma~\ref{thm:killing bound} there exists a bound $\kappa\ind{max}$
on the absolute values of the corresponding principal curvatures.
From Lemma~\ref{thm:collar deformation 1} we get a $R>0$ and an
$\epsilon>0$ and from Lemma~\ref{thm:immersion bound} a $\delta>0$,
such that the immersions $f$ and $\tilde{f}$ obey
\eqref{eq:immersion bound} on $B(p,R)$ if the corresponding
polynomial Killing fields $\zeta$ and $\tilde{\zeta}$ and marked
points $(\lambda_0,\,\lambda_1)$ and
$(\tilde{\lambda}_0,\,\tilde{\lambda}_1)$ obey \eqref{eq:killing
bound 1} and \eqref{eq:killing bound 2}. Therefore these conditions
imply that all $p\in M$ have a neighbourhood, on which $\tilde{f}$
extends to a local one-sided Alexandrov embedding. We remark that
due to the uniform bound on the chord-arc ratio from
Lemma~\ref{thm:ec diameter bound} and the choice of $R>0$ in
Lemma~\ref{thm:collar deformation 1} these local one-sided
Alexandrov embeddings are not affected by the restriction of
$\tilde{f}$ to $M\setminus B(p,R)$. Due to Corollary~\ref{thm:unique
embedded 2} these local one-sided Alexandrov embeddings can be glued
together to an immersion $\tilde{f}:\tilde{N}\rightarrow\Sp^3$ with
boundary $\partial\tilde{N}=M$. A subsequence of any Cauchy sequence
in $\tilde{N}$ with Riemannian metric induced by $\tilde{f}$ is
contained in one of the local Alexandrov embeddings from $O$ to
$\Sp^3$. Hence $\tilde{N}$ with the Riemannian metric induced by
$\tilde{f}$ is complete. Since $\tilde{f}$ has non-negative mean
curvature, it is a one-sided Alexandrov embedding.
\end{proof}
\begin{definition}
For $g\in\mathbb{N}_0$ let $\pkcmc$ denote the subset of triples
$(\xi,\lambda_0,\lambda_1)\in\pksg\times\Sp^1\times\Sp^1$ containing
the initial values (not necessarily without roots) together with
both marked points of finite type {\sc{cmc}} cylinders
$f:M\rightarrow\Sp^3$. Let $\pkcmc^+$ denote the subset of triples
$(\xi,\lambda_0,\lambda_1)\in\pkcmc$ corresponding to finite type
{\sc{cmc}} cylinders $f:M\rightarrow\Sp^3$ with non-negative mean
curvature (with respect to the normal determined by the orientation
of $M$).
\end{definition}
\begin{theorem}\label{thm:spec_AE}
Let $\xi\in\pksg$ be the initial value and $(\lambda_0,\lambda_1)$
the marked points of a one-sided Alexandrov embedded {\sc{cmc}}
cylinder $f:N\rightarrow\Sp^3$. Then a {\sc{cmc}} cylinder of finite
type $\tilde{f}:M\rightarrow\Sp^3$ with non-negative mean curvature
and with initial value $\tilde{\xi}\in\pksg$ and marked points
$(\tilde{\lambda}_0,\tilde{\lambda}_1)$ extends to one-sided
Alexandrov embedding, if one of the following conditions are
satisfied:
\begin{enumerate}
  \item[\Con{M}] The polynomial $a(\lambda)=-\lambda\det(\xi)$ has
    $2g$ pairwise distinct roots and $\tilde{\xi}$ belongs to
    the same isospectral set $\mathcal{K}_a\subset\pksg$ as $\xi$
    with the same marked points.
  \item[\Con{N}] The initial value $\xi$ has only unimodular roots and
    $\tilde{\xi}$ belongs to the same isospectral set
    $\mathcal{K}_a\subset\pksg$ as $\xi$ with the
    same marked points.
  \item[\Con{O}] The polynomial $a(\lambda)=-\lambda\det(\xi)$ has
    $2g$ pairwise distinct roots and
    $(\tilde{\xi},\tilde{\lambda}_0,\tilde{\lambda}_1)$ belongs
    to the same connected component of $\pkcmc^+$ as
    $(\xi,\lambda_0,\lambda_1)$.
  \item[\Con{P}] The initial value $\xi$ has only unimodular roots and
    $(\tilde{\xi},\tilde{\lambda}_0,\tilde{\lambda}_1)$ belongs
    to the same connected component of $\pkcmc^+$ as
    $(\xi,\lambda_0,\lambda_1)$.
\end{enumerate}
In particular, real continuous deformations of spectral curves with
non-negative mean curvature and decreasing $G$ \eqref{eq:G} preserve
one-sided Alexandrov embeddedness.
\end{theorem}
\begin{proof}
The initial values $\xi$ and $\tilde{\xi}$ extend to unique
polynomial Killing fields $\zeta$ and $\tilde{\zeta}$, respectively.

We first show \Con{M}. We obtain the map $f:M \to \Sp^3$ by
composing an affine immersion of $M \cong \Sp^1 \times \R$ into the
real part of the Jacobian with the Sym-Bobenko formula. In case that
$a$ has pairwise distinct roots, the real part of the Jacobian is a
real $g$ dimensional compact torus, and acts freely and transitively
on the corresponding isospectral set $\mathcal{K}_a$. Due to
Lemma~\ref{thm:compact} the isospectral sets are compact. A
continuous action of a finite dimensional Lie group on a compact
space is uniformly continuous. Hence for all $\tilde{f}$ in a small
neighbourhood of $f$ in the isospectral set
\eqref{eq:immersion bound} is satisfied on $q\in M$ instead of
$q\in B(p,R)$.  Due to Lemma~\ref{thm:killing bound} and
Rosenberg's Lemma~\ref{th:Rosenberg_Lemma} the cut locus function
$c_f$ is uniformly bounded from below.
Corollary~\ref{thm:collar deformation 2} implies that the set of
one-sided Alexandrov embedded {\sc{cmc}} cylinders in the isospectral
set is open and closed. This proves \Con{M}.

Now we prove \Con{N}. If the polynomial $a$ of the initial value
$\xi$ has roots of higher order, then the
corresponding isospectral set $\mathcal{K}_a\subset\pksg$ has a
stratification into strata, on which the real part of the Jacobian
acts transitively. Due to Lemma~\ref{thm:dense stratum} the initial
value $\xi$ belongs to the stratum of highest dimension, which is
dense in the whole isospectral set. In this case the arguments above
concerning the case~\Con{M} carry over and show that whenever a
stratum contains a one-sided Alexandrov embedded {\sc{cmc}} cylinder,
then all elements of this stratum correspond to one-sided Alexandrov
embedded {\sc{cmc}} cylinders. Lemma~\ref{thm:convergency} now implies \Con{N}.

Now we prove part~\Con{O}. Again we show that the set of all spectral
curves of arithmetic genus $g$
(i.e.\ the set of polynomials $a(\lambda)=-\lambda\det(\xi)$ of
degree $2g$), whose isospectral sets contain one-sided
Alexandrov embedded {\sc{cmc}} cylinders, is open and closed. Since
the isospectral sets are compact, the Hausdorff distance between the
isospectral sets defines a metric on the space
of spectral data of given arithmetic
genus. Lemma~\ref{thm:convergency} implies that the
isospectral sets of one-sided Alexandrov embeddings are closed and
open in $\pkcmc^+$. This proves \Con{O}.

Due to \Con{N} the arguments of the proof of \Con{O}
carry over to \Con{P}, and concludes the proof.
\end{proof}
%
%
\section{One-sided Alexandrov embedded {\sc{cmc}} cylinders in $\Sp^3$
of finite type} \label{sec:AE}
In this concluding chapter we classify one-sided Alexandrov embedded
{\sc{cmc}} cylinders of finite type in the 3-sphere. We first show
in Theorem \ref{thm:onesided revolution} that there is a 1-parameter
family of flat one-sided Alexandrov embedded {\sc{cmc}} cylinders,
and a 2-parameter family of one-sided Alexandrov embedded rotational
{\sc{cmc}} cylinders of spectral genus $g=1$. We then prove (Lemma
\ref{thm:continuous genus 1}) that an arbitrary one-sided Alexandrov
embedded {\sc{cmc}} cylinder of finite type in the 3-sphere can be
continuously deformed into a flat one-sided Alexandrov embedded
{\sc{cmc}} cylinder, while preserving the condition of one-sided
Alexandrov embeddedness throughout the whole deformation. We next
prove (Lemma \ref{thm:not genus two}) that a surface with bubbletons
is not one-sided Alexandrov embedded. Putting the above results
together in Theorem \ref{thm:main1} gives us the following
classification: A one-sided Alexandrov embedded {\sc{cmc}}
cylinder of finite type in the 3-sphere is a surface of revolution.

Since an embedded {\sc{cmc}} torus in the 3-sphere is covered by a
one-sided Alexandrov embedded {\sc{cmc}} cylinder of finite type,
this result confirms the conjecture by Pinkall and Sterling
\cite{PinS} that the only embedded {\sc{cmc}} tori in the 3-sphere
are tori of revolution. In particular, since by a result by Hsiang
and Lawson \cite{HsiL}, the only embedded minimal torus of
revolution is the Clifford torus, this affirms the Lawson
Conjecture. We conclude the paper with a generalization of an
'unknottedness' result by Lawson \cite{Law:unknot}, and show that
for a one-sided Alexandrov embedded {\sc{cmc}} cylinder of finite
type in the 3-sphere, the 3-manifold is diffeomorphic to the
cartesian product $\overline{\mathbb{D}}\times\R$, where
$\overline{\mathbb{D}}$ is the closed unit disk.

We first turn our attention to the spectral data of one-sided
Alexandrov embedded rotational {\sc{cmc}} cylinders of spectral
genus $g \leq 1$.
\begin{theorem}\label{thm:onesided revolution}
There exists a family of spectral data of one-sided Alexandrov embedded
{\sc{cmc}} cylinders with spectral genus $g\leq 1$ parameterized
by the mean curvature $H \geq 0$ and $\alpha\in[0,1)$.
The corresponding spectral curves are given by
$a(\kappa) = \kappa^2 + \alpha$. The boundary of
the moduli $(H,\,\alpha) \in [0,\,\infty) \times [0,\,1)$ consists of
\begin{description}
  \item[flat cylinders in $\Sp^3$] $H \in [0,\,\infty ],\,\alpha = 0$,
  \item[minimal cylinders in $\Sp^3$] $H=0,\,\alpha \in [0,\,1)$,
\end{description}
In case of flat cylinders in $\Sp^3$ the spectral genus zero curves
have a real double point at $\kappa =0$ and no other real double
points. This case contains all spectral data of
one-sided Alexandrov embedded {\sc{cmc}} cylinders
of geometric genus zero. In case
$\alpha \neq 0$ the spectral curves have no real double points. As
$\alpha \to 0$ the two branch points coalesce into a double point at
$\kappa =0$.
\end{theorem}
\begin{proof}
We recall the spectral data for cylinders of spectral genus zero,
and the expressions for $\ln \mu$ of \eqref{eq:a_vac} and $d \ln
\mu$ of \eqref{eq:b_vac}. At $\kappa = \pm \kappa_0$ we require the
closing conditions that $\ln \mu \in \pi i \Z$. Thus for some
integers $m,\,n \in \Z$ we have
\begin{equation} \begin{split} \label{eq:b_eq}
    \left. \ln \mu \right|_{\kappa_0} \in \pi i \Z &\Longleftrightarrow
    4(b_0\kappa_0 - b_1)^2 = n^2(\kappa_0^2 +1)\,, \\
    \left. \ln \mu \right|_{-\kappa_0} \in \pi i \Z &\Longleftrightarrow
    4(b_0\kappa_0 + b_1)^2 = m^2(\kappa_0^2 +1)\,.
\end{split}
\end{equation}
We make the following claim: If a  cylinder is one-sided Alexandrov
embedded then $m = n =\pm 1$ (for this ensures that the surface is
simply wrapped with respect to the rotational period). To prove this
claim first note that any spectral genus zero cylinder is a covering
of a flat embedded torus. The complement of this flat embedded torus
with respect to $\Sp^3$ consists of two connected components
$\mathcal{D}_{\pm}$, both diffeomorphic to $\mathbb{D} \times
\Sp^1$. Assume the mean curvature vector points into
$\mathcal{D}_+$. For a one-sided Alexandrov embedded {\sc{cmc}}
cylinder the extension $f:N \to \Sp^3$ is a surjective immersion
onto the closure $\bar{\mathcal{D}}_+$ of $\mathcal{D}_+$. Hence this
map is a covering map. The fundamental group of
$\bar{\mathbb{D}}\times\Sp^1$ is isomorphic to
$\Z$. Now the covers of a topological space are in one-to-one
correspondence with subgroups of the fundamental group
\cite[\S14]{St}, and all non-trivial subgroups of
$\pi_1(\bar{\mathbb{D}}\times\Sp^1) \cong \Z$
correspond to compact covers. Hence the only non-compact cover of
$\bar{\mathbb{D}}\times\Sp^1$ is the universal cover
$\bar{\mathbb{D}}\times\R$. Therefore $f$ is the universal covering
map. In particular the period of the cylinder is a primitive period in
the kernel of
$$
H_1(\Sp^1 \times \Sp^1,\,\Z) \to
H_1(\bar{\mathbb{D}}\times\Sp^1,\,\Z )\,.
$$
This implies $m=n=\pm 1$ and proves the claim.

Returning to the proof of the theorem, then \eqref{eq:b_eq} with
$m^2=n^2=1$ simplify to $b_0b_1\kappa_0 = 0$. Thus if $\kappa_0 \neq
0$ then $b_0b_1=0$. Now $|H|$ attains all values in $[0,\,\infty)$
when either $\kappa_0$ ranges over values $\kappa_0 \in (0,\,1]$, or
when $\kappa_0$ ranges over values $\kappa_0 \in [1,\,\infty)$.
Hence we can pick a $\kappa_0$ either from $(0,\,1]$ or from
$[1,\,\infty)$. The M\"obius transformation~\eqref{eq:moebius}
$\kappa\mapsto-\frac{1}{\kappa}$ interchanges these two intervals and
the two cases $b_0=0$ and $b_1=0$. Hence we may assume without loss
of generality that $b_1=0$. Then $4\kappa_0^2 b_0^2 = \kappa_0^2+1$.

Now we claim that with this choice only $\kappa_0\in[1,\,\infty)$
correspond to one-sided Alexandrov embedded cylinders. In fact, we
have seen above that the period of the cylinder should correspond to
an element in the kernel of
$H_1(\partial\mathcal{D}_+,\mathbb{Z})\rightarrow
H_1(\bar{\mathcal{D}}_+,\mathbb{Z})$.
The M\"obius transformation~\eqref{eq:moebius}
$\kappa\mapsto-\frac{1}{\kappa}$ interchanges the cases
$\kappa_0\in(0,\,1]$ and $\kappa_0\in[1,\,\infty)$. Hence we may
consider both families of flat cylinders as two copies of one
family of tori considered as cylinders with respect to different
periods. In the limit $\kappa_0\rightarrow 0$ the length of the period
tends to infinity, and in the limit $\kappa_0\rightarrow\infty$ the
length of the period is bounded. Hence the period of the first family
are the rotation period of $\bar{\mathcal{D}}_-$, i.e.\ a primitive
element of the kernel of $H_1(\partial\mathcal{D}_-,\mathbb{Z})\rightarrow
H_1(\bar{\mathcal{D}}_-,\mathbb{Z})$ and therefore the
translation period of $\bar{\mathcal{D}}_+$. The period
of the second family are the rotation period of $\bar{\mathcal{D}}_+$.
Therefore only the family $\kappa_0\in[1,\infty)$ corresponds to
one-sided Alexandrov embedded cylinders.

For a double point $\kappa$ this means that
$(\kappa_0^2 +1)\,\kappa^2 = n^2(\kappa^2
+1)\,\kappa_0^2$ for some $n \in \Z$. We have
\begin{align*}
0\leq\frac{\kappa^2}{\kappa^2+1}&\leq 1\quad
\mbox{ for real }\kappa\mbox{ and }&
\frac{1}{2}&\leq\frac{\kappa_0^2}{1+\kappa_0^2}<1
\quad\mbox{ for }\kappa_0 \in[1,\,\infty).
\end{align*}
Hence for $\kappa_0\in[1,\,\infty)$
this family has only the real double points
$\kappa=\pm\kappa_0$ with $n=\pm 1$ and $\kappa=0$ with $n=0$. Due to
the closing condition~(ii) in Definition~\ref{thm:spec_bobenko} the
former has to be preserved. If we open the latter double points then
$\ln \mu$ remains a meromorphic function on the genus 1 spectral
curve. We thus altogether obtain families parameterized by
$\kappa_0\in[1,\,\infty)$ and $\alpha\in[0,1)$:
\begin{align*}
&a(\kappa) = \kappa^2 + \alpha\,,
\qquad \qquad \qquad \nu^2=(\kappa^2+1)(\kappa^2 +\alpha)\,,&\\
&\ln\mu = 2\pi i \,b_2\,\frac{\kappa^2 +\alpha}{\nu}\,,
\qquad d\ln\mu = 2\pi i \,b_2\,\frac{\kappa^2 +1-\alpha}{(\kappa^2 +1)\nu}\,.
\end{align*}
At $\pm \kappa_0$ the closing conditions $\ln\mu = \pm \pi i$ must hold
and thus $4(\kappa_0^2 +\alpha)^2b_2^2 = (\kappa_0^2 +\alpha)(\kappa_0^2 +1)$
giving $4(\kappa_0^2 +\alpha)\,b_2^2 = \kappa_0^2+1$.
Therefore double points have to satisfy
$$
\frac{\kappa^2+\alpha}{\kappa^2 +1} =
n^2\frac{\kappa_0^2+\alpha}{\kappa_0^2 +1}\,\mbox{ for some }n\in\Z\,.
$$
But for $\kappa_0 \in [1,\,\infty)$ and $\alpha\in(0,1)$ this equation
has no real solutions besides $n=\pm 1$ and $\kappa=\pm\kappa_0$.
Hence these families have no real double points.
\end{proof}
In addition to the two boundary components there exists two limiting
cases: When $H=\infty,\,\alpha \in [0,\,1]$ we obtain unduloidal
{\sc{cmc}} cylinders in $\R^3$; When $H\in[0,\infty),\, \alpha=1$,
the resulting surfaces are {\emph{chain of spheres}}. Note that as a
consequence of Theorem \ref{thm:onesided revolution}, when $\alpha
\neq 0$ there is no real branch point at $\Delta = \pm 2$. Hence it
is only possible to increase the genus by opening two conjugate
double points in this case.
\begin{lemma}\label{thm:2 real zeroes}
Every spectral data of a one-sided Alexandrov embedded {\sc{cmc}}
cylinder in $\Sp^3$ can be deformed into spectral data of a
one-sided Alexandrov embedded {\sc{cmc}} cylinder in $\Sp^3$
with arbitrary large mean curvature $H$ by a
continuous deformation that increases $H$ and decreases $G$
\eqref{eq:G}.
\end{lemma}
\begin{proof}
The real part is the fixed point set of the anti-linear involution
$\rho$. It is a compact one-dimensional submanifold of the domain of
$\Delta$, which is isomorphic to $\Sp^1$. Since the values of
$\Delta$ at both $\kappa_0,\,\kappa_1$ are equal to $\pm 2$ both
segments of the real part between $\kappa_0$ and $\kappa_1$ contain
an odd number of branch points of $\Delta$ (counted with
multiplicities). The simultaneous movements of $\kappa_0$ and
$\kappa_1$ directed inward to one of these two segments, increases
the mean curvature with respect to the inward normal (determined by
the orientation of $M$). This segment we call the {\bf{short
segment}}, and the other segment the {\bf{long segment}}. We use a
M\"obius transformation to ensure that $\kappa = \infty$ is
contained in the long segment.

We claim, that the movement of a simple real branch point of $\Delta$
within the short segment, which increases the value of $\Delta$ at
the simple branch point, if it is a local minimum, and decreases the
value of $\Delta$ otherwise, increases the mean curvature. In fact,
due to \eqref{eq:integrability_2} the signs of $\dot{\kappa}_i$ are
equal to the signs of $-c(\kappa_i)/b(\kappa_i)$. We may assume that
the simple branch point is a zero of $b$. Hence $c$ has the same
zeroes as $b$ with the exception of the simple branch point, whose
value of $\Delta$ is changed. Hence the sign of the function
$\kappa\mapsto-c(\kappa)/b(\kappa)$ on the short segment changes
the sign only at the simple branch point. The two points $\kappa_0$
and $\kappa_1$ sit on different sides of this simple branch point.
Hence the sign of $\dot{H}$ does not depend on the position of
$\kappa_0$ and $\kappa_1$, as long they stay on the corresponding
side of the simple branch point. If $\kappa_0$ and $\kappa_1$ sit
nearby the simple branch point, then due to \eqref{eq:dotDelta}
and \eqref{eq:integrability_2} the movement, which increases the value
of $\Delta$ at local minima and decreases the values of $\Delta$ at
local maxima, moves $\kappa_0$ and $\kappa_1$ towards each other.
This proves the claim.

Due to condition~\Con{F} the real branch point at $\Delta=\pm 2$
have odd orders. As a first step we shall move each such real
branch point of odd order in the short segment at $\Delta=\pm 2$ by a
small movement into $\Delta\in(-2,2)$. If we choose $a$ to have at all
real branch points at $\Delta=\pm 2$ roots of the same order as
$\Delta^2-4$, then $b$ has at all real branch points roots of the same
order as $\Delta'$. Moreover, the corresponding $\frac{c}{b}$ changes
the sign at the real branch point, and the arguments concerning simple
real branch points in the short segment carry over. This deformation
shortens the short segment, and therefore increases the mean
curvature. These movements do not change $G$ \eqref{eq:G} and
increases the geometric genus into $G$. All other deformations shall
decrease the geometric genus.

In a second step we increase the values of $\Delta$ of all simple
real branch points within the short segment, which are local minima,
and decrease the values of $\Delta$ of all simple real branch points
within the short segment, which are local maxima. This again
shortens the short segment and increases the mean curvature. All
higher order real branch points within the short segment we deform
into pairs of complex conjugate branch points close to the real part
and eventually one simple real branch point. Similar as in the
third step of the proof of Lemma~\ref{thm:connected} we can achieve
with increasing mean curvature a situation with only one simple real
branch point within the short segment. If we move this real branch
point arbitrarily close to the two marked points $\kappa_0$ and
$\kappa_1$, then the mean curvature becomes arbitrarily
large with decreasing $G$ \eqref{eq:G}.
\end{proof}

The arguments show that we can deform the spectral data of a
one-sided Alexandrov embedded {\sc{cmc}} cylinder in $\Sp^3$ into
the spectral data of a {\sc{cmc}} cylinder in $\R^3$. If the two
marked points $\kappa_0$ and $\kappa_1$ coalesce, then the
corresponding surface in $\Sp^3$ shrinks to a point. But if we
enlarge simultaneously all distances of $\Sp^3$ on an appropriate
scale, then the corresponding {\sc{cmc}} cylinders converge to a
{\sc{cmc}} cylinder in $\R^3$. In fact, the enlargement of the
distances of $\Sp^3$ corresponds to spheres in $\R^4$ of enlarged
radius and therefore also to three-dimensional space forms of
diminishing constant sectional curvature. Our arguments can be used
to show, that all spectral data of one-sided Alexandrov embedded
{\sc{cmc}} cylinders in $\Sp^3$ can be continuously deformed within
this class into spectral data of one-sided Alexandrov embedded
{\sc{cmc}} cylinders in $\R^3$. In \cite{KorKS} it is shown that
they have genus at most equal to one. Hence we could have used
\cite{KorKS} instead of the subsequent Lemmas to prove the main
Theorem~\ref{thm:main1}. Conversely, our arguments can be used to
show that all one-sided Alexandrov embedded {\sc{cmc}} cylinders in
$\R^3$ of finite type are Delaunay surfaces. We expect that our
arguments extend to all one-sided Alexandrov embedded
{\sc{cmc}} cylinders in $\R^3$ with constant Hopf differential
(compare \cite{MaOs,Sch,Ta}).

\begin{lemma}\label{thm:continuous genus 1}
All spectral data of one-sided Alexandrov embedded {\sc{cmc}}
cylinders in $\Sp^3$ can be continuously deformed within the class
of one sided Alexandrov embedded {\sc{cmc}} cylinders into spectral
data of flat cylinders in $\Sp^3$ described in
Theorem~\ref{thm:onesided revolution}.
\end{lemma}
\begin{proof}
In a first step, due to Lemma~\ref{thm:2 real zeroes}, we can deform
the spectral data of one-sided Alexandrov embedded
{\sc{cmc}} cylinders in $\Sp^3$ continuously with decreasing $G$
\eqref{eq:G} and increasing $H$ into spectral data with only one
simple real branch point in the short segment and arbitrary large mean
curvature.

In a second step we apply a similar procedure on the real branch
points in the long segment. We can control the mean curvature by
moving the only real branch point within the short segment closer
to the two marked points $\kappa_0$ and $\kappa_1$. Hence we obtain
spectral data with only two real branch points arbitrary close to
the two marked points $\kappa_0$ and $\kappa_1$ with arbitrary large
mean curvature.

In a third step we apply the deformation described in the fourth step
of the proof of Lemma~\ref{thm:connected}. We obtain spectral data of
genus at most equal to one with two simple real branch points very close
to the two branch points at $\kappa_0$ and $\kappa_1$.

In a fourth step we move the real branch point, which connects the
sheets with lower labels, away from the two real branch points at
$\kappa_0$ and $\kappa_1$ along $\Delta\in[-2,2]$ to and fro until
it reaches the place $\Delta=2$ on the sheets with labels $1^\pm$.
We claim that we can preserve positive mean curvature by moving, if
necessary, the other real branch point in the other direction along
$\Delta\in[-2,2]$ to and fro. The function $\Delta$ is equal to
$\mu+\mu^{-1}=2\cosh(\ln\mu)$. Hence it suffices to show that we
preserve the values of $\kappa$ at both marked points by changing the
values of $\ln\mu$ at the other real branch point. The polynomial $c$,
which corresponds to the deformation preserving the mean curvature is
proportional to $c(\kappa)=(\kappa-\kappa_0)(\kappa-\kappa_1)$. The
derivatives $\dot{\ln\mu}=\dot{\mu}\mu^{-1}$ at the real branch points,
that is the zeroes of the corresponding polynomial $b$, are equal to the
values of $\tfrac{c(\kappa)}{\nu}$ at the real branch points due
to \eqref{eq:def_c}. We remark that since $d\ln\mu$ vanishes at these real
branch points any change of the values of $\kappa$ at these real
branch points has no influence on the derivatives $\dot{\ln\mu}$ of
the values of $\ln\mu$ at these real branch points. As long as the
real branch point has not reached $\Delta=2$ on the sheets with labels
$1^\pm$, the function $\nu$ has no roots on the real part. The
two integrals of $d\ln\mu$ along the short and long segments are
preserved under the deformation. Both real roots of $b$ are local
extrema of $\tfrac{1}{2\pi i}\ln\mu$ on the real part. During the
deformation the value of $\tfrac{1}{2\pi i}\ln\mu$ is increased in
case of a local maximum, and decreased in case of a local
minimum. Hence during the deformation one of the real roots of $b$ stays in
the short segment, while the other real root stays in the long
segment. Therefore we can move, without changing the mean curvature,
the real branch point arbitrarily close to $\Delta=2$ on the sheets
with labels $1^\pm$. Due to Lemma~\ref{thm:H continuous} there exist
for any $\epsilon>0$ a $\delta >0$ such that the movement from
$\Delta=2-\delta$ to $\Delta =2$ changes the mean curvature by at most
$\epsilon$. This proves the claim.

Due to Theorem~\ref{thm:spec_AE} the
whole deformation preserves the one-sided Alexandrov embeddedness.
Hence  the final spectral data belong to the flat cylinders in
$\Sp^3$ described in Theorem~\ref{thm:onesided revolution}.
\end{proof}
\begin{lemma}\label{thm:not genus two}
The {\sc{cmc}} cylinders with bubbletons of finite type, whose
$\Delta$ \eqref{eq:covering} corresponds to one of the flat
cylinders described in Theorem~\ref{thm:onesided revolution} are not
one-sided Alexandrov embedded.
\end{lemma}
\begin{proof}
Due to Theorem~\ref{thm:spec_AE}~\Con{O} it is enough to show that
the {\sc{cmc}} cylinders, which are obtained from the continuous
deformation of a $\Delta$ corresponding to a flat cylinder described
in Theorem~\ref{thm:onesided revolution} by moving two simple
complex conjugate branch points from $\Delta=\pm 2$ into
$\Delta\in\C\setminus[-2,2]$ and increasing the genus by two, are
not one-sided Alexandrov embedded. Afterwards these two places at
$\Delta=\pm 2$ have no branch points. In the following first two steps
we can move the unique real branch point in the short segment
arbitrarily close to the two marked points and make the mean curvature
arbitrarily large.

In a first step we move these branch
points to the long segment of the real part between the two marked
points $\kappa_0$ and $\kappa_1$, which contains the real branch
point at $\Delta=2$. Afterwards this segment contains two additional
real branch points.

In a second step we move these real branch points
on the real part away from each other. Then between them two new sheets become
connected through a new segment of the real part. Hence we can move
one of these two real branch points several times along
$\Delta\in[-2,2]$ to and fro until it reaches one of the two places
at $\Delta=\pm 2$ without branch points. Hence the genus is reduced
by one.

In a third step we move the only real branch point in the short segment
away from the two marked points along $\Delta\in[-2,\,2]$ to and fro
until it reaches the other place $\Delta=\pm 2$ without other branch
points. In order to preserve positivity of the mean curvature, we move
the other real root of $b$ in the long segment away from the two
marked points along $\Delta \in [-2,\,2]$ to and fro. The arguments of
the fourth step in the proof of Lemma \ref{thm:continuous genus 1}
ensure this is always possible. The whole deformation has twice
deformed one zero of $b$, and two zeroes of $a$ into a
real double point, and thus arrive at a flat cylinder of spectral
genus zero. It cannot belong to the family of flat cylinders
described in Theorem~\ref{thm:onesided revolution}, since the
corresponding family of spectral genus one curves described in this
Theorem do not have real double points. Due to
Theorem~\ref{thm:spec_AE}~\Con{O}-\Con{P} these deformations
preserve one-sided Alexandrov embeddedness. Hence all spectral data
obtained by adding to the flat cylinders in $\Sp^3$ described in
Theorem~\ref{thm:onesided revolution} a bubbleton do not correspond
to one-sided Alexandrov embedded {\sc{cmc}} cylinders in $\Sp^3$.
\end{proof}

\begin{example}
As an example of the technique used in the proof of Lemma \ref{thm:not
  genus two} we demonstrate how to deform the simplest case of
  bubbletons at the branch point over $\Delta =2$ connecting the
  sheets with labels $2^\pm$ and $3^\pm$:
We describe the movement of the branch point connecting the sheets
which have labels with superscript $+$. The other branch point is
  moved along the conjugate path. First we move the branch point
  connecting the sheets with labels $2^+$ and $3^+$ counter clockwise
  around the interval $\Delta \in [-2,\,2]$ until it meets the origin
  $\Delta =0$ from the lower half-plane. Our branch point now
  meets the conjugate branch point at the real part and connects
sheets with labels $1^+$ and $4^+$.
Further, the branch
  point over $\Delta =-2$ which initially connected the sheets with
  labels $1^+$ and $2^+$, now connects the sheets with labels $4^+$
  and $2^+$, and the branch
  point over $\Delta =-2$ which initially connected the sheets with
  labels $3^+$ and $4^+$, now connects the sheets with labels $3^+$
  and $1^+$. Since the two coalescing branch points are now real, we
  can separate them along the real part. We move one of them
  connecting sheets with labels $1^+$ and $1^-$ to the left until
  $\Delta=-2$, and then move it further connecting the sheets with
  labels $3^+$ and $3^-$ to the right until $\Delta =2$. Now the genus
  is reduced to one. Next we simultaneously move the vertical cut
between the sheets
  with labels $1^+$ and $4^+$ and the vertical cut between the sheets
  with labels $1^-$ and $4^-$ along the same way until it reaches
  $\Delta =-2$ and continue on until it reaches the unique real branch
  point in the short segment. This has the effect that over $\Delta
  =-2$, the branch point which previously connected the sheets with
  labels $4^\pm$ and $2^\pm$, now connects the sheets with labels $1^\pm$
  and $2^\pm$, and the branch
  point over $\Delta =-2$ which previously connected the sheets with
  labels $3^\pm$ and $1^\pm$, now connects the sheets with labels $3^\pm$
  and $4^\pm$. The vertical cuts now connect the sheets with labels
  $3^\pm$ and $2^\pm$. Now we move the real branch point in the short
  segment to $\Delta =2$. To preserve positivity of the mean
  curvature, we move the other real branch point, which connects the
  sheets with labels $4^+$ and $4^-$, towards $\Delta =2$ and possibly
  further to and fro. The genus is thus reduced to zero.
  We arrive at a constellation described in the proof of
  Theorem~\ref{thm:onesided revolution} with $\kappa_0 \in
  (0,\,1)$, which we know not to be one-sided Alexandrov embedded.
\end{example}
\begin{theorem}\label{thm:main1}
A one-sided Alexandrov embedded {\sc{cmc}} cylinder of finite type
in the 3-sphere is a surface of revolution.
\end{theorem}
\begin{proof}
We will show that all spectral data of one-sided Alexandrov embedded
{\sc{cmc}} cylinders of finite type in $\Sp^3$ are described in
Theorem~\ref{thm:onesided revolution}. Due to
Lemma~\ref{thm:continuous genus 1} all spectral data corresponding
to a one--sided Alexandrov embedded {\sc{cmc}} cylinder can be
continuously deformed within the class of $\Delta$ corresponding to
one-sided Alexandrov embedded {\sc{cmc}} cylinders by a combination
of continuous deformations described in
Theorem~\ref{thm:spec_AE}~\Con{N}-\Con{P} into one of the spectral
data of the flat cylinders in $\Sp^3$ described in
Theorem~\ref{thm:onesided revolution}. Due to Lemma~\ref{thm:not
genus two} the $\Delta$ corresponding to these flat cylinders in
$\Sp^3$ can be continuously deformed within the class of one-sided
Alexandrov embedded {\sc{cmc}} cylinders only into the $\Delta$
described in Theorem~\ref{thm:onesided revolution}. All {\sc{cmc}}
cylinders with bubbletons of finite type, whose $\Delta$
\eqref{eq:covering} are described in Theorem~\ref{thm:onesided
revolution} can be deformed within the class of deformations
preserving the one-sided Alexandrov embeddedness into the {\sc{cmc}}
cylinders with bubbletons described in Lemma~\ref{thm:not genus
two}. Hence they are not one-sided Alexandrov embedded.
Consequently, all spectral data corresponding to one-sided
Alexandrov embedded {\sc{cmc}} cylinders in $\Sp^3$ are described in
Theorem~\ref{thm:onesided revolution}. They are all surfaces of
revolution!
\end{proof}
\begin{corollary} \label{th:Pinkall-Sterling}
All one-sided Alexandrov embedded {\sc{cmc}} tori in the 3-sphere
are tori of revolution. In particular, all embedded {\sc{cmc}} tori
in the 3-sphere are tori of revolution.
\end{corollary}
\begin{proof}
Let $f:N\rightarrow\Sp^3$ be a one-sided Alexandrov embedding of a
torus $M=\partial N\simeq\mathbb{T}^2$. We will show that there
exists a cover $\tilde{N}\rightarrow N$, whose composition with $f$
is a {\sc{cmc}} cylinder which is also one-sided Alexandrov
embedded. If the class in $\pi_1(M)$ of a smooth embedding
$\gamma:\Sp^1\hookrightarrow M$ belongs to the kernel of
$\pi_1(M)\rightarrow\pi_1(N)$, then, due to Dehn's Lemma, $\gamma$
is the boundary of a smooth embedding from the closed
two-dimensional disc $\bar{\mathbb{D}}\hookrightarrow N$. The
intersection number of another path $\Sp^1\rightarrow M$ with
$\gamma$ in $M$ is equal to the intersection number with the disc in
$N$. Hence the kernel of $\pi_1(M)\rightarrow\pi_1(N)$ is isotropic
with respect to the intersection form and therefore cyclic. In
particular, there exists a subgroup of $\Gamma\subset\pi_1(N)$,
whose pre-image under $\pi_1(M)\rightarrow \pi_1(N)$ is isomorphic
to $\mathbb{Z}$. Furthermore, we can choose the subgroup $\Gamma$
such that the restriction of $\pi_1(M)\rightarrow\pi_1(N)$ to the
subgroup $\mathbb{Z}\subset\pi_1(M)$ is surjective onto $\Gamma$.
This subgroup $\Gamma\subset\pi_1(N)$ corresponds to a covering
$\tilde{N}\rightarrow N$, whose boundary
$\tilde{M}=\partial\tilde{N}$ is a cylinder \cite[\S14]{St}. Hence
there exists a one-sided Alexandrov embedded cylinder
$\tilde{f}:\tilde{N}\rightarrow\Sp^3$, which is the composition of
$f$ with a covering map. Obviously the spectral data of $f$ and
$\tilde{f}$ coincide. Due to Pinkall and Sterling \cite{PinS} this
cylinder is of finite type, and by Theorem \ref{thm:main1}, it is a
surface of revolution.
\end{proof}
Hsiang and Lawson \cite{HsiL} prove that there are no embedded
minimal tori of cohomogeneity one. Hence the Clifford torus is the
only embedded minimal torus of revolution. By Corollary
\ref{th:Pinkall-Sterling} the only embedded {\sc{cmc}} tori are tori
of revolution, and we thus affirm Lawson's conjecture.
\begin{corollary} The Clifford torus is the only embedded minimal
torus in the 3-sphere.
\end{corollary}
Since all one-sided Alexandrov embedded {\sc{cmc}} cylinders in
$\Sp^3$ are surfaces of revolution around a closed geodesic, the
ambient 3-manifold is diffeomorphic to
$\overline{\mathbb{D}}\times\R$, where $\overline{\mathbb{D}}$
denotes the closed unit disk. We thus have the following
generalization of Lawson's 'unknottedness' result \cite{Law:unknot}.
\begin{corollary}
For all one-sided Alexandrov embedded {\sc{cmc}} cylinders of finite
type in the 3-sphere, the 3-manifold is diffeomorphic to the
cartesian product $\overline{\mathbb{D}}\times\R$.
\end{corollary}
%
%
\bibliographystyle{amsplain}

\begin{thebibliography}{10}

\bibitem{Abr}
U.~Abresch, \emph{Constant mean curvature tori in terms of elliptic
functions},
  J. Reine U. Angew Math. \textbf{374} (1987), 169--192.

\bibitem{Abr:twi}
\bysame, \emph{Old and new doubly periodic solutions of the
sinh-{G}ordon
  equation}, Seminar on new results in nonlinear partial differential equations
  (Bonn, 1984), Aspects Math., E10, Vieweg, Braunschweig, 1987, pp.~37--73.

\bibitem{AA}
R.~Alexander and S.~Alexander, \emph{Geodesics in {R}iemannian
  manifolds-with-boundary}, Indiana Univ. Math. J. \textbf{30} (1981), no.~4,
  481--488.

\bibitem{Ale0}
A.~D. Alexandrov, \emph{Uniqueness theorems for surfaces in the
large. {I}},
  Vestnik Leningrad. Univ. \textbf{11} (1956), no.~19, 5--17.

\bibitem{Bob:tor}
A.~I. Bobenko, \emph{All constant mean curvature tori in
{$\mathbb{R}^3$},
  {$\mathbb{S}^3$}, {$\mathbb{H}^3$} in terms of theta-functions}, Math. Ann.
  \textbf{290} (1991), 209--245.

\bibitem{Bob:cmc}
\bysame, \emph{Constant mean curvature surfaces and integrable
equations},
  Russian Math. Surveys \textbf{46} (1991), 1--45.

\bibitem{Bob:2x2}
\bysame, \emph{Surfaces in terms of 2 by 2 matrices. old and new
integrable
  cases}, Harmonic maps and integrable systems, Aspects of Mathematics, vol.
  E23, Vieweg, 1994.

\bibitem{Bur:iso}
F.~E. Burstall, \emph{Isothermic surfaces: conformal geometry,
{C}lifford
  algebras and integrable systems}, Integrable systems, geometry, and topology,
  AMS/IP Stud. Adv. Math., vol.~36, Amer. Math. Soc., Providence, RI, 2006,
  pp.~1--82.

\bibitem{BurFPP}
F.~E. Burstall, D.~Ferus, F.~Pedit, and U.~Pinkall, \emph{Harmonic
tori in
  symmetric spaces and commuting {H}amiltonian systems on loop algebras}, Ann.
  of Math. \textbf{138} (1993), 173--212.

\bibitem{BurK}
F.~E. Burstall and M.~Kilian, \emph{Equivariant harmonic cylinders},
Quart. J.
  Math. \textbf{57} (2006), 449--468.

\bibitem{BurP_adl}
F.~E. Burstall and F.~Pedit, \emph{Harmonic maps via
{A}dler-{K}ostant-{S}ymes
  theory}, Harmonic maps and integrable systems, Aspects of Mathematics, vol.
  E23, Vieweg, 1994.

\bibitem{BurP:dre}
\bysame, \emph{Dressing orbits of harmonic maps}, Duke Math. J.
\textbf{80}
  (1995), no.~2, 353--382.

\bibitem{Car:thesis}
E.~Carberry, \emph{On the existence of minimal tori in ${S}^3$ of
arbitrary
  spectral genus}, Ph.D. thesis, Princeton University, 2002.

\bibitem{Con2}
J.~B. Conway, \emph{Functions of one complex variable. {II}},
Graduate Texts in
  Mathematics, vol. 159, Springer-Verlag, New York, 1995.

\bibitem{Dav}
E.~B. Davies, \emph{{$L^1$} properties of second order elliptic
operators},
  Bull. London Math. Soc \textbf{17} (1985), 417--436.

\bibitem{DorH:per}
J.~Dorfmeister and G.~Haak, \emph{On constant mean curvature
surfaces with
  periodic metric}, Pacific J. Math. \textbf{182} (1998), 229--287.

\bibitem{DorPW}
J.~Dorfmeister, F.~Pedit, and H.~Wu, \emph{Weierstrass type
representation of
  harmonic maps into symmetric spaces}, Comm. Anal. Geom. \textbf{6} (1998),
  no.~4, 633--668.

\bibitem{ErcKT}
N.~M. Ercolani, H.~Kn{\"o}rrer, and E.~Trubowitz,
\emph{Hyperelliptic curves
  that generate constant mean curvature tori in {${\bf R}\sp 3$}}, Integrable
  systems (Luminy, 1991), Progr. Math., vol. 115, Birkh\"auser Boston, Boston,
  MA, 1993, pp.~81--114.

\bibitem{Esc}
J.-H. Eschenburg, \emph{Maximum principle for hypersurfaces},
Manuscripta Math.
  \textbf{64} (1989), no.~1, 55--75.

\bibitem{Fis-ColS}
D.~Fischer-Colbrie and R.~Schoen, \emph{The structure of complete
stable
  minimal surfaces in 3-manifolds of nonnegative scalar curvature}, CPAM
  \textbf{33} (1980), 199--211.

\bibitem{ForLR}
S.~Fornari, J.~H.~S. deLira, and J.~Ripoll, \emph{Geodesic graphs
with constant
  mean curvature in spheres}, Geom. Dedic. \textbf{90} (2002), 201--216.

\bibitem{GriS1}
P.~G. Grinevich and M.~U. Schmidt, \emph{Period preserving
nonisospectral flows
  and the moduli space of periodic solutions of soliton equations}, Phys. D
  \textbf{87} (1995), no.~1-4, 73--98.

\bibitem{GroKS:Tri}
K.~Gro{\ss}e-Brauckmann, R.~Kusner, and J.~M. Sullivan,
\emph{Triunduloids:
  embedded constant mean curvature surfaces with three ends and genus zero}, J.
  {R}eine {A}ngew. {M}ath. \textbf{564} (2003), 35--61.

\bibitem{Heb}
J.~Hebda, \emph{Cut loci of submanifolds in space forms and in the
geometries
  of {M}\"obius and {L}ie}, Geom. Dedicata \textbf{55} (1995), no.~1, 75--93.

\bibitem{Hit:tor}
N.~Hitchin, \emph{Harmonic maps from a 2-torus to the 3-sphere}, J.
  Differential Geom. \textbf{31} (1990), no.~3, 627--710.

\bibitem{Hop}
H.~Hopf, \emph{Differential geometry in the large}, Lecture Notes in
  Mathematics, vol. 1000, Springer-Verlag, 1983.

\bibitem{HsiL}
W-Y. Hsiang and H.~B. Lawson, Jr., \emph{Minimal submanifolds of low
  cohomogeneity}, J. Differential Geometry \textbf{5} (1971), 1--38.

\bibitem{Hur}
A.~Hurwitz, \emph{Ueber {R}iemann'sche {F}l\"achen mit gegebenen
  {V}erzweigungspunkten}, Math. Ann. \textbf{39} (1891), no.~1, 1--60.

\bibitem{Jag}
C.~Jaggy, \emph{On the classification of constant mean curvature
tori in {${\bf
  R}\sp 3$}}, Comment. Math. Helv. \textbf{69} (1994), no.~4, 640--658.

\bibitem{Kap1}
N.~Kapouleas, \emph{Complete constant mean curvature surfaces in
{E}uclidean
  three space}, Ann. Math. \textbf{131} (1990), 239--330.

\bibitem{Kap2}
\bysame, \emph{Compact constant mean curvature surfaces in
{E}uclidean three
  space}, J. Diff. Geom. \textbf{33} (1991), 683--715.

\bibitem{Kap3}
\bysame, \emph{Constant mean curvature surfaces by fusing {W}ente
tori}, Proc.
  Natl. Acad. Sci. USA \textbf{89} (1992), 5695--5698.

\bibitem{KapY}
N.~Kapouleas and S.-D. Yang, \emph{Minimal surfaces in the
three-{S}phere by
  doubling the {C}lifford {T}orus}, arXiv:math/0702565v1.

\bibitem{KarPS}
K.~Karcher, U.~Pinkall, and I.~Sterling, \emph{New minimal surfaces
in
  {$S\sp3$}}, J. Differential Geom. \textbf{28} (1988), 169--185.

\bibitem{Kil:del}
M.~Kilian, \emph{On the associated family of {D}elaunay surfaces},
Proc. Amer.
  Math. Soc \textbf{132} (2004), no.~10, 3075--3082.

\bibitem{KilMS}
M.~Kilian, I.~McIntosh, and N.~Schmitt, \emph{New constant mean
curvature
  surfaces}, Experiment. Math. \textbf{9} (2000), no.~4, 595--611.

\bibitem{KilSS}
M.~Kilian, N.~Schmitt, and I.~Sterling, \emph{Dressing {CMC}
n-{N}oids}, Math.
  Z. \textbf{246} (2004), no.~3, 501--519.

\bibitem{KorKR}
N.~Korevaar, R.~Kusner, and J.~Ratzkin, \emph{On the nondegeneracy
of constant
  mean curvature surfaces}, Geom. Funct. Anal. \textbf{16} (2006), no.~4,
  891--923.

\bibitem{KorKS}
N.~Korevaar, R.~Kusner, and B.~Solomon, \emph{The structure of
complete
  embedded surfaces with constant mean curvature}, J. Diff. Geom. \textbf{30}
  (1989), no.~2, 465--503.

\bibitem{Kri_77}
I.~M. Kri{\v{c}}ever, \emph{Methods of algebraic geometry in the
theory of
  nonlinear equations}, Uspehi Mat. Nauk \textbf{32} (1977), no.~6(198),
  183--208, 287, English translation: Russian Math. Surveys 32 (1977), no. 6,
  185--213.

\bibitem{KusMP}
R.~Kusner, R.~Mazzeo, and D.~Pollack, \emph{The moduli space of
complete
  embedded constant mean curvature surfaces}, Geom. Funct. Anal. \textbf{6}
  (1996), 120--137.

\bibitem{Law:compact}
H.~B. Lawson, Jr., \emph{Compact minimal surfaces in {$S\sp{3}$}},
Global
  Analysis (Proc. Sympos. Pure Math., Vol. XV, Berkeley, Calif., 1968), Amer.
  Math. Soc., Providence, R.I., 1970, pp.~275--282.

\bibitem{Law:S3}
\bysame, \emph{Complete minimal surfaces in {$S\sp{3}$}}, Ann. of
Math. (2)
  \textbf{92} (1970), 335--374.

\bibitem{Law:unknot}
\bysame, \emph{The unknottedness of minimal embeddings}, Invent.
Math.
  \textbf{11} (1970), 183--187.

\bibitem{MaOs}
V.~A. Mar{\v{c}}enko and I.~V. Ostrovs{\cprime}ki{\u\i}, \emph{A
  characterization of the spectrum of the {H}ill operator}, Mat. Sb. (N.S.)
  \textbf{97(139)} (1975), no.~4(8), 540--606, 633--634.

\bibitem{McI:tor}
I.~McIntosh, \emph{Harmonic tori and their spectral data},
math.DG/0407248. To
  appear in 'Surveys on {G}eometry and {I}ntegrable {S}ystems', Advanced
  Studies in Pure Mathematics, 2007.

\bibitem{Mee}
W.~H. {Meeks III}, \emph{The topology and geometry of embedded
surfaces of
  constant mean curvature}, J. Diff. Geom. \textbf{27} (1988), 539--552.

\bibitem{MeePR:limit}
W.~H. {Meeks III}, J.~Perez, and A.~Ros, \emph{Limit leaves of a cmc
lamination
  are stable}, arXive:0801.4345.

\bibitem{MonR:min}
S.~Montiel and A.~Ros, \emph{Minimal immersions of surfaces by the
first
  eigenfunctions and conformal area.}, Invent. Math. \textbf{83} (1986),
  153--166.

\bibitem{MonR:alex}
\bysame, \emph{Compact hypersurfaces: the {A}lexandrov theorem for
higher order
  mean curvatures}, Differential geometry, Pitman Monogr. Surveys Pure Appl.
  Math. \textbf{52} (1991), 279--296.

\bibitem{PinS}
U.~Pinkall and I.~Sterling, \emph{On the classification of constant
mean
  curvature tori}, Ann. Math. \textbf{130} (1989), 407--451.

\bibitem{Poh}
K.~Pohlmeyer, \emph{Integrable hamiltonian systems and interaction
through
  quadratic constraints}, Comm. Math. Phys. \textbf{46} (1976), 207--221.

\bibitem{PoeT}
J.~P{\"o}schel and E.~Trubowitz, \emph{Inverse spectral theory},
Pure and
  Applied Mathematics, vol. 130, Academic Press Inc., Boston, MA, 1987.

\bibitem{PreS}
A.~Pressley and G.~Segal, \emph{Loop groups}, Oxford Science
Monographs, Oxford
  Science Publications, 1988.

\bibitem{Rin}
W.~Rinow, \emph{Die innere {G}eometrie der metrischen {R}\"aume},
Die
  Grundlehren der mathematischen Wissenschaften, Bd. 105, Springer-Verlag,
  Berlin, 1961.

\bibitem{Ros}
A.~Ros, \emph{A two-piece property for compact minimal surfaces in a
  three-sphere.}, Indiana Univ. Math. J. \textbf{44} (1995), no.~3, 841--849.

\bibitem{Ros:com}
H.~Rosenberg, \emph{Private communication}, 27/03/2008.

\bibitem{Sch}
M.~U. Schmidt, \emph{Integrable systems and {R}iemann surfaces of
infinite
  genus}, Mem. Amer. Math. Soc. \textbf{122} (1996), no.~581, viii+111.

\bibitem{Sch:cmclab}
N.~Schmitt, \emph{cmclab}, http://www.gang.umass.edu/software.

\bibitem{SKKR}
N.~Schmitt, M.~Kilian, S.~Kobayashi, and W.~Rossman,
\emph{Unitarization of
  monodromy representations and constant mean curvature trinoids in three
  dimensional space forms}, J. London Math. Soc. \textbf{75} (2007), 563--581.

\bibitem{Sch:stable}
R.~Schoen, \emph{Stable minimal surfaces in three manifolds}, Annals
of
  Mathematics Studies \textbf{103} (1983), 111--126.

\bibitem{Ser}
J.~P. Serre, \emph{Algebraic goups and class fields}, Graduate Texts
in
  Mathematics, No. 117, Springer-Verlag, New York, 1988.

\bibitem{St}
N.~Steenrod, \emph{The {T}opology of {F}ibre {B}undles}, Princeton
Mathematical
  Series, vol. 14, Princeton University Press, Princeton, N. J., 1951.

\bibitem{SteW:bub}
I.~Sterling and H.~Wente, \emph{Existence and classification of
constant mean
  curvature multibubbletons of finite and infinite type}, Indiana Univ. Math.
  J. \textbf{42} (1993), no.~4, 1239--1266.

\bibitem{Symes_80}
W.~W. Symes, \emph{Systems of {T}oda type, inverse spectral
problems, and
  representation theory}, Invent. Math. \textbf{59} (1980), no.~1, 13--51.

\bibitem{TerU}
C.~Terng and K.~Uhlenbeck, \emph{B\"{a}cklund transformations and
loop group
  actions}, Comm. Pure and Appl. Math \textbf{LIII} (2000), 1--75.

\bibitem{Ta}
V.~Tkachenko, \emph{Spectra of non-selfadjoint {H}ill's operators
and a class
  of {R}iemann surfaces}, Ann. of Math. (2) \textbf{143} (1996), no.~2,
  181--231.

\bibitem{Uhl}
K.~Uhlenbeck, \emph{Harmonic maps into lie groups (classical
solutions of the
  chiral model)}, J. Diff. Geom. \textbf{30} (1989), 1--50.

\bibitem{UmeY:tori}
M.~Umehara and K.~Yamada, \emph{A deformation of tori with constant
mean
  curvature in {${\bf R}\sp 3$} to those in other space forms}, Trans. Amer.
  Math. Soc. \textbf{330} (1992), no.~2, 845--857.

\bibitem{Wal}
R.~Walter, \emph{{E}xplicit examples to the $h$-problem of {H}einz
{H}opf},
  Geom. Dedic. \textbf{23} (1987), 187--213.

\bibitem{Wen}
H.~C. Wente, \emph{Counterexample to a conjecture of {H}. {H}opf},
Pac. J.
  Math. \textbf{121} (1986), 193--243.

\bibitem{Wen:twi}
\bysame, \emph{Twisted tori of constant mean curvature in
{$\mathbb{R}^3$}},
  Seminar on new results in nonlinear partial differential equations (Bonn,
  1984), Aspects Math., E10, Vieweg, Braunschweig, 1987, pp.~1--36.

\bibitem{Wol:elliptic}
T.~Wolff, \emph{Recent work on sharp estimates in second order
elliptic unique
  continuation}, Jour. Geom. Ana. \textbf{3} (1993), no.~6, 621--650.

\end{thebibliography}
\def\cprime{$'$}
\providecommand{\bysame}{\leavevmode\hbox
to3em{\hrulefill}\thinspace}
\providecommand{\MR}{\relax\ifhmode\unskip\space\fi MR }
\providecommand{\MRhref}[2]{%
  \href{http://www.ams.org/mathscinet-getitem?mr=#1}{#2}
} \providecommand{\href}[2]{#2}

\end{document}